\documentclass[leqno]{amsart}
\usepackage{amsfonts}
\usepackage{hyperref}
\usepackage{amsmath}
\usepackage{amssymb}
\usepackage{cite}
\usepackage{caption}
\usepackage{tabularx}
\usepackage{color}
\usepackage{graphicx}
\usepackage{subfigure}
\usepackage{amsgen,amsthm}
\usepackage{epsfig,epstopdf,color,bm}
\usepackage{wrapfig}

\theoremstyle{plain}
\newtheorem{theorem}{Theorem} [section]

\newtheorem{proposition}[theorem]{Proposition}

\theoremstyle{remark}
\newtheorem{remark}[theorem]{Remark}

\allowdisplaybreaks

\newcommand{\p}{\partial}

\newcommand{\artanh}{\mathrm{artanh}}
\newcommand{\be}{\begin{equation}}
\newcommand{\ee}{\end{equation}}
\newcommand{\ba}{\begin{array}}
\newcommand{\ea}{\end{array}}
\newcommand{\bea}{\begin{eqnarray}}
\newcommand{\eea}{\end{eqnarray}}
\newcommand{\beas}{\begin{eqnarray*}}
\newcommand{\eeas}{\end{eqnarray*}}

\def\C{{\mathbb C}}
\def\R{{\mathbb R}}
\def\O{\mathcal O}
\def\F{\mathcal F}

\def\({\left(}
\def\){\right)}
\def\<{\left\langle}
\def\>{\right\rangle}
\def\le{\leqslant}
\def\ge{\geqslant}
\def\sech{\mathrm{sech}}
\def\d{{\partial}}
\def\eps{\varepsilon}
\def\om{\omega}
\def\si{{\sigma}}

\def\Eq#1#2{\mathop{\sim}\limits_{#1\rightarrow#2}}

\def\d{{\partial}}
\def\eps{\varepsilon}
\def\om{\omega}
\def\si{{\sigma}}

\DeclareMathOperator{\RE}{Re}

\numberwithin{equation}{section}

\title[LogNLS with repulsive harmonic
potential]{Numerical study of the logarithmic Schr\"odinger equation
  with repulsive harmonic potential}

\author[R. Carles]{R\'emi Carles}
\author[C. Su]{Chunmei Su}

\address{Univ Rennes, CNRS\\ IRMAR - UMR 6625\\ F-35000
  Rennes, France}
\email{Remi.Carles@math.cnrs.fr}

\address{Yau Mathematical Sciences Center\\Tsinghua University\\
  Beijing 100084, China}
\email{sucm@tsinghua.edu.cn}

\begin{document}

\begin{abstract}
We consider the Schr\"odinger equation with a
logarithmic nonlinearity and a repulsive harmonic potential. Depending
on the parameters of the equation, the solution may or may not be
dispersive. When dispersion occurs, it does with an exponential rate
in time. To control this, we change the unknown function through a
generalized lens transform. This approach neutralizes the possible
boundary effects, and could be used in the case of the Schr\"odinger
equation without potential. We then employ standard splitting methods on the new
equation via a nonuniform grid, after the logarithmic nonlinearity has been regularized.  We
also discuss the case of a power nonlinearity and give some results concerning the error estimates of the first-order Lie-Trotter splitting method for both cases of nonlinearities. Finally extensive numerical experiments are reported to investigate the dynamics of the equations.
\end{abstract}

\maketitle

\section{Introduction}
\label{sec:intro}

We consider the initial value problem of a time-dependent Schr\"odinger equation with potential $V$ and logarithmic nonlinearity:
\begin{equation}\label{RlogSE}\left\{
\begin{aligned}
&i\d_t u +\frac{1}{2}\Delta u=V u+\lambda
  u\ln\(|u|^2\),\quad x\in \R^d,\quad t>0,\\
&u(0,x)=u_0(x),\end{aligned}\right.
\end{equation}
where $u_0$ is the initial data.

In the absence of the logarithmic nonlinearity ($\lambda=0$), there is a large amount of literature on the so-called repulsive potential setting \cite{Vega, ReedSimon2}, i.e.,
\[(\partial_r V )_+(r):= \sup\limits_{x\in\R^d: |x|=r}(\partial_r V)_+(x)=0.\]
Particularly, solutions of the Schr\"odinger equation with repulsive $1/r^n$-type long-range potentials play an essential role in many physical processes such as collisions of similar
atoms in a radiation field \cite{burke1998, gallagher1989, orzel1998},
molecular spectra converging to thresholds where two fragments can
interact via resonant dipole-dipole interactions \cite{abraham1995,
  miller1993, napolitano}, atom-electron, atom-ion interactions
\cite{gao1999} and atom-surface interactions \cite{stwalley}. Compared
to the free-potential equation, the repulsive potential creates
acceleration of the field \cite{antoine2009, lorin2007}. In the case
of a repulsive harmonic potential, $V(x) = -|x|^2$, the acceleration
is exponential in time, as recalled below.

In the absence of potential ($V=0$), the logarithmic Schr\"odinger equation (logNLS) has been adopted in many physical models \cite{BEC,buljan,hansson,Hef85,KEB00,DMFGL03,Zlo10} since it was introduced in \cite{BiMy76}. For instance, as proposed in \cite{bouharia2015, Zloshchastiev}, the logarithmic model
may generalize the Gross-Pitaevskii equation, used in the case of two-body interaction, to the
case of multi-body interaction. A particular feature of the logarithmic
nonlinearity is that it leads to very special solitary waves, called
\emph{Gaussons} in \cite{BiMy76,BiMy79} when $\lambda<0$. These
solitary waves are orbitally stable \cite{Ar16,Caz83}. Furthermore, for $\lambda<0$,
no solution is dispersive (\cite[Proposition~4.3]{Caz83}), while for $\lambda>0$,
every solution is dispersive with an enhanced rate compared to the usual rate
of the free Schr\"odinger equation, and the modulus of the solution
converges to a universal Gaussian profile  \cite{CaGa18}.

When it comes to the logarithmic model with potential, a harmonic trapping
potential was considered in \cite{bouharia2015} to describe the
logarithmic Bose-Einstein Condensation:
\begin{equation}
  \label{eq:logNLSconf}
  i\d_t u +\frac{1}{2}\Delta u = \frac{\omega^2}{2}|x|^2u +\lambda
  u\ln\(|u|^2\),\quad x\in \R^d.
\end{equation}
Due to the presence of the potential, stationary solutions
(generalized Gaussons) are available and orbitally stable in both cases $\lambda<0$ \cite{ ACS20, bouharia2015} and $\lambda>0$ \cite{Carles21}.
For the logNLS with repulsive potential,
\begin{equation}
  \label{logNLSrep}
  i\d_t u +\frac{1}{2}\Delta u = -\frac{\omega^2}{2}|x|^2u +\lambda
  u\ln\(|u|^2\),\quad x\in \R^d,
\end{equation}
it was shown in \cite[Proposition~1.3]{Carles21} that for $\lambda\in
\R$ and any
\begin{align*}
  u_0\in \Sigma &:=H^1\cap \F(H^1)= \left\{ f\in H^1(\R^d),\quad x\mapsto |x|
    f(x)\in L^2(\R^d)\right\},\\
&\|f\|_\Sigma:= \|f\|_{L^2}+\|\nabla f\|_{L^2}+ \|xf\|_{L^2},
\end{align*}
there exists a unique solution $u \in L^\infty_{\rm
  loc}(\R;\Sigma)\cap C(\R;L^2(\R^d))$ to \eqref{logNLSrep}.
Furthermore, when $\lambda>0$, the solution shares the same
exponential dispersive
rate as in the linear case $\lambda=0$, and no solitary wave or
universal dynamics exists in this case. In the nondispersive
case ($\lambda<0$), the situation is different. It was shown in
\cite{ZZ20} that \eqref{logNLSrep} admits at least one positive bound
state under a suitable range of the coefficients. Recently, we proved
that if $-\lambda>\omega>0$, there exist two positive stationary
Gaussian solutions, which are
orbitally unstable \cite{carles2021}.

Along the numerical part, there have been few studies for the model with repulsive potential or logarithmic nonlinearity. For the Schr\"odinger equation with repulsive potential \eqref{logNLSrep} and without logarithmic interaction ($\lambda=0$), the strong dispersive effects made it almost impossible to simulate the dynamics on a truncated domain with naive homogeneous or periodic boundary conditions \cite{antoine2009}. To overcome this difficulty, some types of artificial or absorbing boundary conditions were proposed for one-dimensional Schr\"odinger equation with a general variable repulsive potential \cite{antoine2009} and for two-dimensional Schr\"odinger equation with a time and space varying exterior potential \cite{antoine2012, antoine2013}. On the other hand, for the logNLS without potential, the singularity of the logarithmic nonlinearity also makes it very challenging to
design and analyze numerical schemes. To overcome the singularity at
the origin, some numerical methods were proposed and analyzed for the
logarithmic Schr\"odinger equation (\eqref{logNLSrep} with $\omega=0$)
based on a global nonlinearity regularization model
\cite{bao2019error, bao2019} or a local energy regularization
approximation \cite{bao2021}. Note that even though it is rather
natural to regularize the nonlinearity at least in the case of splitting methods
(otherwise the solution of the ODE is singular), it may not be
necessary to do so, as proved in \cite{Zouraris-p}, in the case of the
Crank-Nicolson method.

To our knowledge, the model
\eqref{logNLSrep} is not motivated by physics: we consider it because it cumulates several difficulties in terms of computation and numerical analysis. Most importantly, the strong dispersion caused by the repulsive harmonic
potential requires to be handled very carefully as simulations are usually
performed on a bounded domain. In the case of a logarithmic
nonlinearity, the solution to \eqref{logNLSrep} may or may not be
dispersive, as recalled above. Moreover, the singularity of the
logarithm at the origin makes the use of splitting methods
delicate. In this paper, we propose a numerical strategy  to
face these two features, and extend it to the case where the
logarithmic nonlinearity is replaced by a power nonlinearity. Our numerical methods are based on a generalized lens transform. This formulation is able to extract the main dispersion as well as oscillation and transforms the original equation into an equivalent one with non-dispersive solutions. This enables standard numerical methods to work successfully. The other technique of our methods is to use a nonuniform temporal grid for the derived equivalent equation according to some designed rule consistent with the transform in the first step. This is the basis to establish the error estimates and numerical experiments show that it is superior than that by using a uniform grid directly. As pointed out in Remark~\ref{rem:lens}, the approach presented here can be adapted to the case when one starts from a nonlinear Schr\"odinger
equation without potential.

The rest of the article is organized as
follows. Section~\ref{sec:features} is devoted to recalling some properties, particularly the
dispersion of the logNLS with repulsive potential
\eqref{logNLSrep}. We introduce the generalized lens transform and
present the splitting methods with a series of non-equidistant time steps in
Section~\ref{sec:splitting}. This approach is extended to the nonlinear Schr\"odinger
equation with repulsive potential and power nonlinearity in
Section~\ref{sec:power}. Some error estimates of the
time discretization are presented in Section~\ref{sec:error}. Extensive numerical experiments are displayed in Section~\ref{sec:num} to show the accuracy and feasibility of the methods and investigate the dispersion properties of the equations with various parameters.

\section{Some features of the logNLS under repulsive potential}\label{sec:features}
In this section, we recall some properties of the dynamics of the
logNLS under repulsive potential. A first unusual property associated
to this logarithmic nonlinearity is that the size of the initial data
plays no role (apart from a purely time dependent oscillation): if $u$
solves \eqref{RlogSE}, then for all
$k\in \C$, so does
\begin{equation*}
  k u(t,x) e^{-it\lambda\ln|k|^2}.
\end{equation*}
In the case of \eqref{logNLSrep} (as well as \eqref{eq:logNLSconf}), the
potential $V(x)$ is the sum of potentials depending only on one space
variable, $V(x) = \sum\limits_{j=1}^d V_j(x_j)$. The second property we emphasize
is a tensorization property, which initially motivated the
introduction of the logarithmic nonlinearity \cite{BiMy76}:
if the initial data is a tensor
product,
\[
  u_0(x) =\prod_{j=1}^d u_{0j}(x_j),
\]
then the solution is given by
\[  u(t,x) =\prod_{j=1}^d u_{j}(t,x_j),
\]
where each $u_j$ solves a one-dimensional equation,
\begin{equation*}
   i\d_t u_j +\frac{1}{2} \d_{x_j}^2 u_j =V_j (x_j)u_j+ \lambda
   \ln\(|u_j|^2\)u_j  ,\quad u_j(0, x_j) =u_{0j} .
 \end{equation*}
 Like the first property recalled above, this is reminiscent of
 \emph{linear} Schr\"odinger equations, but effects caused
 by the logarithmic nonlinearity actually modify the dynamics, as we will see.

The third property, which is extremely convenient for numerical
simulations, is that initial Gaussian data leads to a solution which
remains Gaussian for all time. In view of the second property, we
consider the case $d=1$, and
we seek the Gaussian solution with the form
\be\label{gauss}
u( t,x)=b(t)e^{-a(t)x^2/2},
\ee
particularly the initial data $u_0(x)=b(0)e^{-a(0)x^2/2}$($a(0)=\alpha-i\beta$, $\alpha>0$) is Gaussian,
for the equation \eqref{logNLSrep}. It can be easily checked that
\begin{equation}\label{eqab}
  i\dot b = \frac{1}{2}ab+\lambda b\ln|b|^2,\quad i\dot a=a^2+2\lambda \mathrm{Re}\, a+\omega^2,
\end{equation}
where $\dot{b}$ is the derivative with respect to time.
Writing $a$ as
\begin{equation}\label{eqa}
  a=\frac{1}{\mu^2}-i\frac{\dot \mu}{\mu},\quad \mu\in\mathbb{R}^+,
\end{equation}
leads to
\begin{equation}\label{eq:tau-general}
  \ddot \mu=\frac{2\lambda}{\mu}+\frac{1}{\mu^3}+\omega^2 \mu,\quad \mu(0)=\frac{1}{\sqrt{\alpha}},\quad \dot \mu(0)=\frac{\beta}{\sqrt{\alpha}}.
\end{equation}
Plugging \eqref{eqa} into \eqref{eqab} and integrating in time yields
\begin{equation}\label{eqb}
  b(t) =b(0) e^{i\theta(t)}\sqrt{\mu(0)/\mu(t)},\quad
  \theta(t)\in \R.
\end{equation}
Multiplying \eqref{eq:tau-general} by $\dot{\mu}$ and integrating, we
get
\begin{equation}\label{eq:tau-energy}
  (\dot \mu)^2 =C_0+ 4\lambda \ln \mu -\frac{1}{\mu^2}+\omega^2 \mu^2,
\end{equation}
where $C_0=\dot{\mu}(0)^2-4\lambda \ln\mu(0)+\frac{1}{\mu(0)^2}-\omega^2\mu(0)^2$ is related to the  initial data. Noticing that $F(q)=C_0+ 4\lambda \ln q -\frac{1}{q^2}+\omega^2 q^2\to -\infty$ when $q\to 0$, this implies that $\mu$ is bounded away from zero, i.e.,
\begin{equation*}
  \exists \delta>0,\quad \mu(t)\ge\delta,\quad\forall t\ge 0.
\end{equation*}
Combining \eqref{gauss}, \eqref{eqa} and \eqref{eqb}, one arrives at
\[|u(t,x)|= \frac{|b(0)|}{\sqrt{\mu(t)\alpha^{1/2}}}e^{-\frac{x^2}{2\mu(t)^2}}.\]
Hence it suffices to study the dynamics of $\mu(t)$, i.e., the ODE \eqref{eq:tau-general}. Concerning this, we have:
\begin{proposition}[Propagation of Gaussian data, \cite{carles2021}]\label{prop1}
  Let $d=1$, $\lambda<0<\omega$. \\
 $1.$ If $\lambda<-\omega$, then \eqref{eq:tau-general} has exactly two stationary solutions, $\mu_\pm=1/\sqrt{k_\pm}$ with $k_{\pm}=-\lambda\pm
  \sqrt{\lambda^2-\omega^2}$, which correspond to the positive stationary solutions of \eqref{logNLSrep}: $\phi_\pm(x)=\exp\left(-\frac{k_\pm}{4\lambda}-\frac{k_\pm x^2}{2}\right)$ and generate a continuous family of solitary
waves,
\[
  u_{\pm,\nu}(t,x) =\phi_{\pm,\nu}(x)e^{i\nu t},\quad \phi_{\pm,\nu}(x) =
  e^{-\frac{\nu}{2\lambda}} \phi_{\pm}(x),\quad \nu\in \R.
\]
The other solutions to \eqref{eq:tau-general}  are either periodic or unbounded, which correspond to time-periodic or dispersive Gaussian solutions to
\eqref{logNLSrep}, respectively.\\
$2.$ If $\lambda=-\omega$, then \eqref{eq:tau-general} has exactly
one stationary solution, $\mu_0= 1/\sqrt\omega$. All the
  other solutions are unbounded. In other words, any Gaussian solution
 to  \eqref{logNLSrep}  which is not of the form
  \begin{equation*}
    e^{i\theta}e^{(2\nu+\omega)/(4\omega)}e^{i\nu t} e^{-\omega x^2/2},\quad \nu
    \in \R,\quad \theta\in [0,2\pi],
  \end{equation*}
  is dispersive.\\
  $3.$ If $\lambda>-\omega$, then every solution to
  \eqref{eq:tau-general} is unbounded and $e^{\omega t}\lesssim \mu(t)\lesssim e^{\omega t}$. This implies every Gaussian solution to \eqref{logNLSrep} disperses exponentially in time.
\end{proposition}

While for $\lambda>0$, it can be easily observed from \eqref{eqab} that there are no stationary Gaussian solutions. Furthermore, it was shown in \cite[Proposition~1.8]{Carles21} that for general initial data $u_0\in\Sigma$, the solution of \eqref{logNLSrep} can be scaled as
\[u(t,x)=\frac{1}{\nu(t)^{1/2}} v\(t, \frac{x}{\nu(t)}\)\exp\(i\frac{\dot{\nu}(t)}{\nu(t)}\frac{x^2}{2}\),\]
where $v$ is bounded and non-dispersive, and
\[\nu(t)\sim \nu_\infty e^{\omega t},\quad
\dot{\nu}(t)\sim \omega \nu_\infty e^{\omega t},\quad t\rightarrow \infty.\]

\section{Splitting numerical methods based on a transformation}
\label{sec:splitting}
It can be seen from the above section that for $\lambda<0$ and
Gaussian initial data, the solution of \eqref{logNLSrep} disperses
very fast, except some special solutions, e.g., standing waves
or time-periodic solutions. While for $\lambda>0$ and general initial
data, each solution disperses exponentially in time. On the one hand,
we are curious whether the solution is dispersive in the case $\lambda<0$,
for general initial data. One the other hand, possible quick
dispersion makes commonly used naive truncation by homogeneous
Dirichlet or periodic boundary conditions infeasible in practical
computation. In this section we introduce a transformation which
extracts the main dispersion such that the simple truncation technique
works well. Then a classical time splitting discretization is presented.

For simplicity of notation, we consider $d=1$ and extensions to higher dimensions are straightforward. Inspired by the lens transform
\cite{niederer73}
and its generalization in the repulsive harmonic case \cite{carles2003}, we change the unknown
$u$ to $v$ via the identity
\be\label{uv}
u(t,x)=\frac{1}{\sqrt{\cosh(\omega t)}}v\(s,y\)e^{i\frac{\omega}{2}x^2\tanh(\omega t)},\quad s=\frac{\tanh(\omega t)}{\omega}, \quad
y=\frac{x}{\cosh(\omega t)}.\ee
Plugging this formula into the equation \eqref{logNLSrep} leads to
\begin{align*}
0&= i\p_t u+\frac{1}{2}\p_{xx}u+\frac{\omega^2}{2}|x|^2u-\lambda u\ln|u|^2\\
&=\frac{e^{i\frac{\omega}{2}x^2\tanh(\omega t)}}{\sqrt{\cosh(\omega t)}}\(i\p_s v\,\sech^2(\omega t)+\frac{1}{2}\sech^2(\omega t)\p_{yy}v-\lambda v\ln|u|^2\)\\
&=\frac{e^{i\frac{\omega}{2}x^2\tanh(\omega t)}}{\cosh^{5/2}(\omega t)}\(i\p_s v+\frac{1}{2}\p_{yy}v-\lambda \cosh^2(\omega t)v\left[\ln|v|^2-\ln(\cosh(\omega t))\right]\).
\end{align*}
Noticing that $\cosh^2(\omega t)=\frac{1}{1-\tanh^2(\omega
  t)}=\frac{1}{1-\omega^2 s^2}$, we see that $v$ solves the equation
\begin{equation*} 
\left\{
\begin{aligned}
&i\partial_s v+\frac{1}{2}\p_{yy} v=\frac{\lambda v}{1-\omega^2
  s^2}\(\ln|v|^2+\frac{1}{2}\ln(1-\omega^2 s^2)\),\quad (s,y)\in (0,
1/\omega)\times\R,\\
&v(0,y)=u(0,y)=u_0(y).
\end{aligned}\right.
\end{equation*}
We emphasize that the time interval $t\in [0,\infty)$ has been
transformed into $s\in [0,1/\omega)$: on such bounded time interval,
dispersive effects will not be too strong, also because we have
filtered out the maximum (exponential) dispersive effects by changing
the space variable. On the other hand, the equation is now
non-autonomous, as a time dependent factor has appeared in front of
the nonlinearity. In the case of splitting methods, this will require
to compute various integrals in the time variable, when solving the
ODE part. The last term above involves a purely time-dependent potential,
and can be absorbed by a gauge transform.
Set
\be\label{vg}
v(s,y)=\kappa(s,y)e^{-i\lambda g(s)},\quad g(s)=\frac{1}{2}\int_0^s
\frac{\ln(1-\omega^2 p^2)}{1-\omega^2 p^2}dp.\ee
Then $\kappa$ satisfies the equation
\begin{equation}\label{weq}
\left\{
\begin{aligned}
&i\partial_s \kappa+\frac{1}{2}\p_{yy} \kappa=\frac{\lambda
  \kappa}{1-\omega^2 s^2}\ln|\kappa|^2,\quad (s,y)\in(0,
1/\omega)\times \R,\\
&\kappa(0,y)=u(0,y)=u_0(y).
\end{aligned}
\right.
\end{equation}
Hence investigating the large-time behavior of $u(t,x)$ reduces to studying the asymptotic dynamics of $\kappa$ when $s$ approaches $1/\omega$. Noticing \eqref{uv}, for each dispersive solution $u(t,x)$, the coefficient in the exponential function in space is of order $e^{-\omega t}$, hence it requires large computational bounded domain increasing exponentially in time. Moreover, due to the term $\exp\(i\frac{\omega}{2}x^2\tanh(\omega t)\)$, which is highly oscillatory in space, it requires very tiny mesh size in practical computation and brings difficulties for dynamics.

To simulate the dynamics of \eqref{weq}, firstly we regularize \eqref{weq} by introducing a small parameter $0<\eps\ll 1$ and truncate the problem on a bounded computational domain $\Omega =(a, b)$ with periodic boundary condition:
\be\label{rweq}\left\{
\begin{aligned}
&i\partial_s \kappa_\eps+\frac{1}{2}\p_{yy} \kappa_\eps=\frac{2\lambda \kappa_\eps}{1-\omega^2 s^2}\ln(|\kappa_\eps|+\eps),\quad (s,y)\in(0, 1/\omega)\times\Omega,\\
&\kappa_\eps(0, y)=u(0, y)=u_0(y).
\end{aligned}\right.
\ee
We use the splitting method to solve the equation \eqref{rweq}, based on the splitting
\[\p_s \phi=A(\phi)+B(s, \phi),\]
where
\[A(\phi)=\frac{i}{2}\p_{yy}\phi,\quad B(s, \phi)=-\frac{2i\lambda \phi}{1-\omega^2s^2}\ln(|\phi|+\eps),\]
and the solutions of the subproblems
\be\label{lin}\left\{
\begin{aligned}
&\p_s \phi(s_n+s, y) = A(\phi), \quad y\in\Omega, \quad s>0,\\
&\phi(s_n, y) = \phi^n(y),
\end{aligned}\right.
\ee
\be\label{nonl}\left\{
\begin{aligned}
&\p_s z(s_n+s, y)=B(s_n+s, z), \quad y\in\Omega, \quad s > 0,\\
&z(s_n, y)=z^n(y).
\end{aligned}\right.
\ee
The associated evolution operators are given by
\begin{align*}
\phi(s_n+s, \cdot)&=\Phi_A^s(\phi^n)=\exp\(\frac{is}{2}\p_{yy}\)\phi^n,\\
z(s_n+s, \cdot)&=\Phi_B^s(s_n, z_n)=z^n\exp\(-2i\lambda \ln(|z^n|+\eps)
\int_0^s\frac{1}{1-\omega^2(s_n+p)^2}dp\)\\
&=z^n (|z^n|+\eps)^{-\frac{i\lambda}{\omega}\ln\(\frac{(1-\omega s_n)(1+\omega(s_n+s))}{(1+\omega s_n)(1-\omega (s_n+s))}\)}.
\end{align*}

For fixed $T>0$ and the number of time steps $N\in \mathbb{N}^*$, we denote the time step $\tau=T/N>0$, $t_n=n\tau$, $s_n=\frac{\tanh(\omega t_n)}{\om}$ for $n\in \mathbb{N}^*$ and set \be\label{delta}
\delta_n =s_{n+1}-s_n.\ee
 We consider the Lie-Trotter splitting for solving \eqref{rweq} by using a series of non-equidistant time steps $\{\delta_n\}$:
\be\label{Lie1}
\kappa_\eps^{n+1}=\Phi_{\rm L}(\kappa_\eps^n)=\Phi_A^{\delta_n}(\Phi_B^{\delta_n}(s_n, \kappa_\eps^n)),\quad \kappa_\eps^0=u_0,\ee
and the Strang splitting scheme:
\be\label{scheme1}
\kappa_\eps^{n+1}=\Phi_{\rm S}(\kappa_\eps^n)=\Phi_A^{\delta_n/2}(\Phi_B^{\delta_n}(s_n, \Phi_A^{\delta_n/2}(\kappa_\eps^n))),\quad \kappa_\eps^0=u_0.\ee
With the value of $\kappa_\eps^n$, we can approximate $u(t_n, x)$ via \eqref{uv} by setting
\be\label{scheme2}
u_\eps^n(x)=\frac{1}{\sqrt{\cosh(\omega t_n)}}\exp\(i\frac{\omega}{2}x^2\tanh(\omega t_n)-i\lambda g(s_n)\)\kappa_\eps^n\(\frac{x}{\cosh(\omega t_n)}\),\ee
where $g$ is defined as \eqref{vg}.

The above scheme \eqref{scheme2} together with \eqref{Lie1} or \eqref{scheme1} displays a time integrator (semi-discretization) for solving \eqref{logNLSrep}. In practical computation, we combine it with the Fourier pseudospectral method for spatial discretization. The schemes are explicit and efficient thanks to the fast Fourier transform (FFT) for calculating \eqref{Lie1} or \eqref{scheme1} and the nonuniform fast Fourier transform (NUFFT) \cite{greengard,jiang2014} for computing \eqref{scheme2}.

\begin{remark}[Case without potential: lens transform]\label{rem:lens}
  In the case of a nonlinear Schr\"odinger equation without potential,
  \begin{equation}
    \label{eq:NLS}
    i\d_tu+\frac{1}{2}\Delta u = f\(|u|^2\)u,\quad (t,x)\in
    (0,\infty)\times \R^d,
  \end{equation}
one may introduce the lens transform
\begin{equation*}
  v(s,y) = \frac{1}{(\cos(\omega
    s))^{d/2}}u(t,x)e^{-i\frac{\omega}{2}|y|^2 \tan(\omega s)},\quad
    t=\frac{\tan (\omega s)}{\omega},\quad x=\frac{y}{\cos(\omega s)}.
\end{equation*}
Note that the relation between $u$ and $v$ is somehow reversed
compared to \eqref{uv}. The study of \eqref{eq:NLS} is equivalent to
the study of
\begin{equation*}
  i\d_s v+\frac{1}{2}\Delta v = \frac{\omega^2}{2}|y|^2v+\frac{1}{\cos^2(\omega s)} f\(
  |v|^2\cos^d(\omega s)\)v,\quad (s,y)\in
  (0,\frac{\pi}{2\omega})\times \R^d.
\end{equation*}
Like above, the time interval is now bounded, and the presence of the
(confining) harmonic potential makes it possible to avoid boundary
effects. This provides an alternative approach to compute the
scattering operator associated to \eqref{eq:NLS}, compared to
\cite{CGM3AS}. Note that in the case of an $L^2$-critical nonlinearity, $f(|z|^2)
=\lambda |z|^{4/d}$, the nonlinear term reduces to $\lambda |v|^{4/d}v$, and
the equation in $v$ is autonomous \cite{carles2002}.
\end{remark}

\section{Power nonlinearity case}
\label{sec:power}
For a comparison, next we consider the nonlinear Schr\"odinger
equation with repulsive potential and power nonlinearity. For
simplicity, we consider the one-dimensional case only (note that with
a power nonlinearity, the tensorization property recalled in
Section~\ref{sec:features} is lost):
\begin{equation}\label{eq:NLSrep}\left\{
  \begin{aligned}
  &i\d_t u +\frac{1}{2}\d_{xx} u = -\omega^2\frac{x^2}{2}u +\lambda
  |u|^{2\si} u ,\quad x\in\R,\\
  & u(0,x)=u_0(x),\end{aligned}\right.
\end{equation}
with $\sigma>0$. Note that in the case of a
power nonlinearity, initial Gaussian data do not propagate as
Gaussians (this can be seen by seeking solutions under the form of a
Gaussian function). The Cauchy problem and scattering theory have been
investigated in \cite{carles2003}. Among other results, it is proven
that in the case $\lambda>0$ (defocusing case, considered in the
simulations in Section 6), every power nonlinearity is short range, in the
sense that the (nonlinear) solution behaves like a solution to the
linear equation (\eqref{eq:NLSrep} with $\lambda=0$): for every
$u_0\in \Sigma$,
\begin{equation}\label{eq:scattering-power}
\exists u_+\in \Sigma,\quad  \left\| e^{-i\frac{t}{2}(\d_{xx}
      +\omega^2x^2)}u(t) - u_+\right\|_{\Sigma}=\O\(e^{-\frac{\sigma\omega}{\sigma+1}t}\)\quad \text{as }t\to+\infty,
\end{equation}
and, as can be shown for instance thanks to Mehler's formula
(e.g. \cite[Equation~(1.8)]{carles2003}),
\begin{equation*}
  e^{i\frac{t}{2}(\d_{xx}
      +\omega^2x^2)}u_+\Eq t{+\infty}
    e^{-i\pi/4}\sqrt{\frac{\omega}{\sinh (\omega t)}}
    \F \(u_+e^{i\omega |\cdot |^2/2}\)\(\frac{\omega x}{\sinh (\omega t)}\) e^{i\omega\coth
   (\omega  t)\frac{x^2}{2}},
\end{equation*}
where we normalize the Fourier transform as
\begin{equation*}
  \F(f)(\xi)=
  \widehat
  f(\xi)=\frac{1}{\sqrt{2\pi}} \int_{\R} e^{-ix  \xi
  }f(x)dx.
\end{equation*}
While for the focusing case $\lambda<0$, finite-time blow-up might occur under some assumptions on $\lambda$, $\omega$, $\sigma$ and $u_0$ (e.g. \cite[Theorem~1.1]{carles2003}), that is, there exists $T>0$ such that
\be\label{blowup}
\lim\limits_{t\rightarrow T}\|\nabla u(t)\|_{L^2}=+\infty.\ee

We study the dispersion properties of \eqref{eq:NLSrep} numerically. With the same formulation \eqref{uv}, we get the equation for $v$ as
\be\label{eq:v}\left\{
\begin{split}
&i\partial_s v+\frac{1}{2}\p_{yy} v=\lambda \(1-\omega^2
  s^2\)^{\frac{\si}{2}-1}|v|^{2\si}v,\quad (s,y)\in (0,
  1/\omega)\times \R,\\
&v(0,y)=u(0,y)=u_0(y).
\end{split}\right.
\ee
Denote $\Phi_B^s$ by the evolution operator of the nonlinear subproblem
\be\label{nonl1}\left\{
\begin{aligned}
&\p_s z(s_n+s,y)=-i\lambda (1-\omega^2(s_n+s)^2)^{\frac{\si}{2}-1}|z|^{2\si}z
,\quad s > 0, \quad y\in\R, \\
&z(s_n,y)=z^n(y).
\end{aligned}\right.
\ee
Then it can be written explicitly as
\[z(s_n+s,y)=\Phi_B^s(s_n, z^n)=z^n\exp\(-i\lambda |z^n|^{2\si}\int_{s_n}^{s_n+s}(1-\omega^2\zeta^2)^{\frac{\si}{2}-1} d\zeta\).\]
Using the same notations as before, i.e., $s_n=\frac{\tanh(\omega t_n)}{\omega}$ with $t_n=n\tau$, and $\delta_n=s_{n+1}-s_n$, we get
\be\label{zsol}
z(s_{n+1},y)=\Phi_B^{\delta_n}(s_n, z^n)=z^n\exp\(-i\lambda |z^n|^{2\si}\int_{t_n}^{t_{n+1}}\sech^\si(\omega \zeta) d\zeta\).\ee
Combining this with the evolution operator of the linear subequation $\Phi_A^s$ \eqref{lin}, we can get the first-order Lie-Trotter splitting and the second-order Strang splitting approximations as
\be\label{Lie}
v^{n+1}=\Phi_{\rm L}(v^n)=\Phi_A^{\delta_n}(\Phi_B^{\delta_n}(s_n, v^n)), \quad v^0=u_0;\ee
and
\be\label{Strang}
v^{n+1}=\Phi_{\rm S}(v^n)=\Phi_A^{\delta_n/2}(\Phi_B^{\delta_n}(s_n, \Phi_A^{\delta_n/2}(v^n))),\quad  v^0=u_0,\ee
respectively.
While the approximation $u^{n+1}$ can be recovered by
\be\label{recu}
u^{n+1}(x)=\frac{1}{\sqrt{\cosh(\omega t_{n+1})}}\exp\(i\frac{\omega}{2}x^2\tanh(\omega t_{n+1})\)v^{n+1}\(\frac{x}{\cosh(\omega t_{n+1})}\).\ee

In practical computation, we combine the splitting integrator with Fourier pseudospectral discretization in space and a highly accurate quadrature rule for approximating the integral in \eqref{zsol}.

\section{On error estimates for the splitting methods}
\label{sec:error}

In both cases, \eqref{logNLSrep} and \eqref{eq:NLSrep}, we have used
the generalized lens transform \eqref{uv} to turn the original
equation into a Schr\"odinger equation without potential, with a
nonautonomous nonlinearity, \eqref{weq} and \eqref{eq:v},
respectively. It is possible to rely on error estimates which have
been established for splitting methods applied to such equations. Three
differences must be taken into account though:
\begin{itemize}
\item A time dependent factor has appeared in front of the
  nonlinearity (expect in the very special $L^2$-critical case
  $\si=2/d$ in \eqref{eq:NLSrep}).
\item The time step is uniform in terms of the $t$ variable,
  hence not in terms of the $s$ variable.
\item The factor $\frac{1}{1-\omega^2 s^2}$ is singular, and not
    integrable, as $s\to 1/\omega$.
\end{itemize}

In the case of nonlinear Schr\"odinger equations without potential
(\eqref{weq} and \eqref{eq:v} without the time dependent factor in
front of the nonlinearity), error estimates for Lie-Trotter time splitting
schemes are available with initial data in $H^1$: see \cite{bao2019} for the
case of a (regularized) logarithmic nonlinearity, and
\cite{Ignat2011,ChoiKoh2021} for the power nonlinearity. Throughout
this section, we assume that the numerical solution is given by a
Lie-Trotter splitting scheme, as opposed to the Strang splitting
scheme considered so far: indeed error estimates for Strang splitting
scheme in the presence of a logarithmic nonlinearity are far less
satisfactory, see \cite[Remark~4]{bao2019}.

In view of \eqref{uv}, applying the gradient $\nabla_y$ to the
solution of \eqref{weq} and \eqref{eq:v} amounts to considering the
action of the vector-field $J$ on the initial unknown $u$, where $J$
is defined by
\begin{align*}
  J(t) &= \omega x \sinh(\omega t)+i\cosh (\omega
         t)\nabla=e^{-itH}\(i\nabla\) e^{-itH}\\
  &=
  i \cosh(\omega t)e^{i\omega \tanh (\omega t)\frac{|x|^2}{2} }
    \nabla  \left( e^{- i\omega  \tanh (\omega
    t)\frac{|x|^2}{2}}\ \cdot \right),
\end{align*}
where
\begin{equation*}
  H=-\frac{1}{2}\Delta -\omega^2\frac{|x|^2}{2}.
\end{equation*}
This vector-field plays a central role in the proofs in
\cite{carles2003}.

The splitting considered for $v$ (or, equivalently, $\kappa$), is
\emph{equivalent} to a splitting for $u$, thanks to the definition
of $s_n$,
\begin{equation*}
  s_n = \frac{\tanh(\omega t_n)}{\omega},\quad t_n= n\tau.
\end{equation*}
Indeed, the case $\lambda=0$ shows that \eqref{uv} maps a
solution $u$ of the linear Schr\"odinger equation with potential to a
solution $v$ of the linear Schr\"odinger equation without potential,
with $u_{\mid t=0}=v_{\mid s=0}$. We readily check that the same holds
regarding the ordinary differential equations. The last formulation
for $J$ shows that $\|J(t)u\|_{L^2(\R^d)}=\|\nabla
v(s)\|_{L^2(\R^d)}$, and we see that it is equivalent to consider
splitting methods,
\begin{itemize}
\item For $v$ with an analysis in $H^1$, that is, involving $v$ and
  $\nabla v$, with $e^{i\frac{s}{2}\Delta}$ as the linear operator;
\item For $u$ with an analysis involving $u$ and $Ju$, with $e^{-itH}$
  as the linear operator.
\end{itemize}
\subsection{Logarithmic nonlinearity}

In \cite{bao2019}, an error estimate is proven with initial data in
$H^1$ (in the 1D case -- data in $H^2$ if $d=2$ or $3$), on any
\emph{bounded} time interval: the constants given by the proof grow at
least exponentially in time. By resuming the same proof step by step,
following either of the two strategies described above, we can prove that for any  $T>0$, there exists $\eps_0>0$ such
that when $0<\eps\le \eps_0$
  and $0\le n\tau\le T$, we have:
  \begin{equation*}
    \|u^\eps(t_n)-
    u^n_\eps\|_{L^2(\R)}=\|\kappa^\eps(s_n)-\kappa^n_\eps\|_{L^2(\R)} \le
    C\(T,\|u_0\|_{\Sigma}\)\ln(\eps^{-1})\tau^{1/2},
  \end{equation*}
where $C(\cdot,\cdot)$ is independent of $\eps\in (0,\eps_0]$, and
$u^\eps$ is related to $\kappa^\eps$ via the same formula as the one
relating $u$ to $\kappa$.

To compare the numerical solution with the exact solution $u$, we then
have to compare $u$ with $u^\eps$, or, equivalently in $L^2$, $\kappa$
with $\kappa^\eps$. We first note that
\cite[Lemma~2.3]{bao2019error} remains unchanged in the presence of a
time dependent factor in front of the logarithm, but we have to stick
to time intervals where this factor remains bounded (or equivalently,
bounded time intervals for $u$): for $S\in [0,1/\omega)$, there exists
$C=C(S)$ such that for all $s\in [0,S]$,
\begin{equation*}
  \frac{d}{ds}\|\kappa^\eps(s)-\kappa(s)\|_{L^2(\R)}^2\le C\(
  \|\kappa^\eps(s)-\kappa(s)\|_{L^2(\R)}^2 +\eps
  \|\kappa^\eps(s)-\kappa(s)\|_{L^1(\R)}^2 \).
\end{equation*}
Note that the transform \eqref{uv} does not preserve the $L^1$-norm
in space ($L^2$ is the only Lebesgue norm which is preserved), but for bounded
 $t\in [0,T]$ (that is, for $s\in [0,S]$ and $S=\frac{\tanh(\omega
   T)}{\omega}\in [0,1/\omega)$), there exists $C(T)$ such that
 \begin{equation*}
 \frac{1}{C(T)}  \|\kappa^\eps(s)-\kappa(s)\|_{L^1(\R)}\le
 \|u^\eps(t)-u(t)\|_{L^1(\R)}\le C(T)
 \|\kappa^\eps(s)-\kappa(s)\|_{L^1(\R)},\forall t\in [0,T].
 \end{equation*}
Modifying e.g. the proof of \cite[Lemma~2.6]{bao2019error}, we readily
prove: there exists $C$ such that for any $v\in \Sigma$,
\begin{equation}\label{eq:GNdual}
  \|v\|_{L^1(\R)}\le C \|v\|_{L^2(\R)}^{1/2}\|v\|_{\Sigma}^{1/2}.
\end{equation}
In view of \cite[Proposition~1.3]{Carles21}, we infer that there
exists $C(T)$ such that
\begin{equation*}
  \|\kappa^\eps(s)-\kappa(s)\|_{\Sigma}\le C(T),\quad \forall s\in
  [0,S], \ S =\frac{\tanh(\omega
   T)}{\omega},
\end{equation*}
and so
\begin{equation*}
  \frac{d}{ds}\|\kappa^\eps(s)-\kappa(s)\|_{L^2(\R)}^2\le C(T)\(
  \|\kappa^\eps(s)-\kappa(s)\|_{L^2(\R)}^2 +\eps
  \|\kappa^\eps(s)-\kappa(s)\|_{L^2(\R)} \).
\end{equation*}
Using Gronwall Lemma, we conclude:
\smallbreak
\begin{proposition}\label{prop61}
  Suppose $d=1$, and consider the Lie-Trotter splitting method \eqref{Lie1} combined with \eqref{scheme2}.
 Let $T>0$ and $u_0\in \Sigma$. There exists $\eps_0>0$ such that when
 $0<\eps\le \eps_0$
  and $0\le n\tau\le T$, we have
  \begin{equation*}
    \|u(t_n)-
    u^n_\eps\|_{L^2(\R)}=\|\kappa(s_n)-\kappa^n_\eps\|_{L^2(\R)} \le
    C\(T,\|u_0\|_{\Sigma}\)\(\ln(\eps^{-1})\tau^{1/2} +\eps\),
  \end{equation*}
where $C(\cdot,\cdot)$ is independent of $\eps\in (0,\eps_0]$.
\end{proposition}

\subsection{Power nonlinearity}

In the case of a power nonlinearity, it is not necessary to regularize
the nonlinearity, but for $\si<1/2$, the nonlinearity in
\eqref{eq:NLSrep} is not $C^2$: working with an $H^2$ regularity in
space may be delicate (one way to overcome the lack of regularity
might be to consider one time derivative, and use the equation to
infer $H^2$ regularity in space, see \cite{CazCourant}). In
\cite{Ignat2011,ChoiKoh2021}, discrete Strichartz inequalities were
used in order to decrease the required regularity, to $H^1$ (still in
the case of Lie-Trotter method). However, as noted in
\cite{IgnatZuazua2006,IgnatZuazua2009}, in the absence of potential,
it is necessary to introduce a frequency cut-off in the free
propagator $e^{i\frac{t}{2}\Delta}$, and consider instead
$e^{i\frac{t}{2}\Delta} \Pi_\tau$, where $\Pi_\tau$ is a Fourier
multiplier of symbol $\chi(\tau^{1/2}\xi)$, where $\chi$ is
(sufficiently) smooth and compactly supported: this corresponds to a
modification of the operator $A$ as defined in Section~\ref{sec:splitting}. This implies that if
one wants to prove error estimates based on discrete Strichartz
inequalities like in \cite{Ignat2011,ChoiKoh2021}, in the case of
\eqref{eq:NLSrep}, one should consider the same frequency cut-off when
working with the unknown $v$. We note that the frequency cut-off
amounts to imposing the frequency localization
\begin{equation*}
  |\xi|\lesssim \tau^{-1/2}.
\end{equation*}
In the case of \eqref{eq:v}, the time step is not uniform, but
\begin{equation*}
  |s_{n+1}-s_n|=\left|\frac{\tanh(\omega(n+1)\tau)}{\omega} -
  \frac{\tanh(\omega\, n\tau)}{\omega}  \right|\le \tau,
\end{equation*}
and so $|\xi|\lesssim \tau^{-1/2}$ implies $|\xi|\lesssim
|s_{n+1}-s_n|^{-1/2}$ uniformly in $n$.

For $s\in[0,S]$ and $S\in [0,1/\omega)$, the time dependent factor in
front of the nonlinearity in \eqref{eq:v} is harmless, and the proofs
of \cite{ChoiKoh2021}  can be repeated, in order to obtain:
\begin{proposition}
  Let $d=1$, $\si>0$, $u_0\in \Sigma$ and $T>0$. In \eqref{Lie}, replace the operator
  $\Phi_A^s$ with
  \begin{equation*}
    \Phi_{\tilde A}^s = \Phi_A^s\Pi_\tau,\quad\text{where}\quad \widehat{\Pi_{\tau} \phi} (\xi)
    = \chi (\tau^{1/2}\xi) \widehat{\phi}(\xi) ,
  \end{equation*}
  and $\chi \in C^\infty(\R)$ is a cut-off function supported in
  $[-2,2]$ such that $\chi \equiv 1$ on $[-1,1]$. There exists
  $C=C(\si,T,\|u_0\|_\Sigma)$ such that
  \begin{equation*}
    \max_{0\le n\tau\le T}\|u^n- u(t_n)\|_{L^2(\R)} = \max_{0\le
      s_n\le S}\|v^n- v(s_n)\|_{L^2(\R)} \le C\tau^{1/2},
  \end{equation*}
where $   S=\frac{\tanh(\omega T)}{\omega}$.
\end{proposition}
\smallbreak

\begin{remark}
  In view of the analysis performed in \cite{CaSu-p}, it is likely
  that in the case $\lambda>0$, the above estimate is uniform in time,
  that is, $C$ can be
  chosen independent of $T$. Indeed, when working with the unknown
  $u$, the vector-field $J$ provides some exponential decay in time of
  various Lebesgue norms in space, as exploited in \cite{carles2003} to
  prove scattering results. The proof in \cite{CaSu-p} being already
  quite technical, we do not explore the details of the argument
  here.
\end{remark}
\smallbreak

\begin{remark}
  In the case $\sigma>1/2$, the nonlinearity in \eqref{eq:NLSrep} is
  $C^2$, and considering initial data in $H^2\cap \F(H^2)$, the above
  error rate $\tau^{1/2}$ can be improved to $\tau$, following the
  approach of \cite{BBD,Lu08}.
\end{remark}

\section{Numerical Results}
\label{sec:num}
In this section, we first test the order of accuracy of the proposed method \eqref{scheme2} combined with the Lie-Trotter splitting \eqref{Lie1} or the Strang splitting \eqref{scheme1}. Then we apply the Strang-splitting method to investigate some long time
dynamics of the Schr\"odinger equation with repulsive potential and logarithmic nonlinearity \eqref{logNLSrep} (or equivalently \eqref{rweq}) or power nonlinearity \eqref{eq:NLSrep} (or equivalently \eqref{nonl1}).

\smallbreak
\noindent{\bf Example 1}. Here, we set $\lambda=-3$, $\omega=2$ in the equation \eqref{logNLSrep} and $\eps=10^{-15}$ in the regularized model \eqref{rweq}. Choose the initial data as
$u_0(x)=Ae^{-\alpha x^2/2}$ with $\alpha =-\lambda-\sqrt{\lambda^2-\om^2}=3-\sqrt{5}$. Then according to Proposition~\ref{prop1}, it generates the solitary wave as
\[u(t,x)=A e^{i\nu t}e^{-\alpha x^2/2},\]
with $\nu=-(2\lambda \ln(A)+\alpha /2)$. Set $A=2$, then the exact solution is given by
\be\label{solit}
u( t,x)=2e^{-\alpha x^2/2}e^{-i(\alpha/2-6\ln2)t}.\ee
We set $\Omega=(-L, L)$ with $L=10$ in \eqref{rweq}.
To quantify the numerical error, we define the error function as
\[e(t_n)=u_\eps^n(x)-u(t_n,x),\]
where $u_\eps^n$ is obtained by \eqref{scheme2} combined with \eqref{Lie1} or \eqref{scheme1}. In practical computation, we compute until $T=2.5$ by using different values of $N$ (the number of time steps) and a fixed mesh size $h=1/2^9$ for spectral spatial discretization, which is small enough for neglecting the spatial error introduced. The time integral in $g(s)$ is approximated by Simpson rule with very fine mesh such that the error introduced by numerical integration is ignorable.

Before presenting the numerical results, we  recall the temporal grid introduced in Section 3. We notice that the computation is performed by utilizing a series of non-equidistant time steps $\delta_n$ defined by \eqref{delta}, which results as a series of numerical solutions $\kappa_\varepsilon^n$ as an approximation for $\kappa$ on a nonuniform grid $\{s_n=\frac{\tanh(\omega t_n)}{\omega}\}$ and a series of solutions $u^n$ via \eqref{scheme2} as an approximation for $u$ on a uniform grid $\{t_n=\frac{nT}{N}\}$. We denote the corresponding methods as Lie I (\eqref{Lie1}) and Strang I (\eqref{scheme1}), respectively, in the following. On the other hand, another intuitive approach is to solve the equation \eqref{rweq} directly by a standard temporal discretization: we set the step size $\delta_n\equiv\delta:=\frac{s_N}{N}$ with $s_N=\frac{\tanh(\omega t_N)}{\omega}$. Then we get a series of solutions $\kappa_\varepsilon^n$ as an approximation for $\kappa$ on a uniform grid $\{s_n=n\delta\}$ and a series of solutions $u^n$ via \eqref{scheme2} as an approximation for $u$ on a nonuniform grid $\{t_n=\frac{\artanh(\omega s_n)}{\omega}\}$. We denote the corresponding methods as Lie II (\eqref{Lie1}) and Strang II (\eqref{scheme1}) for convenience afterwards.

\begin{figure}[htbp]
\begin{center}
\includegraphics[width=2.45in,height=1.5in]{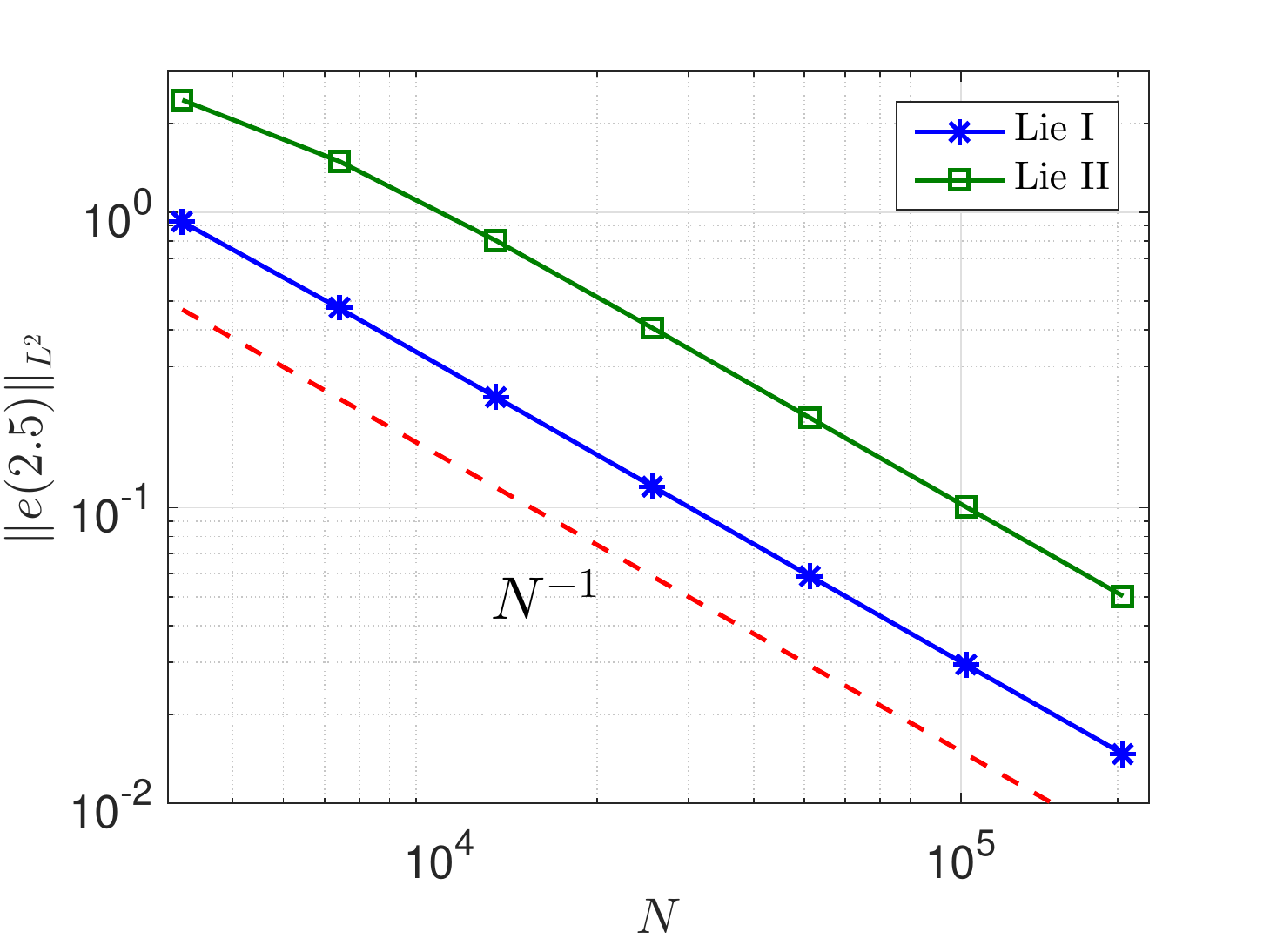}
\includegraphics[width=2.45in,height=1.5in]{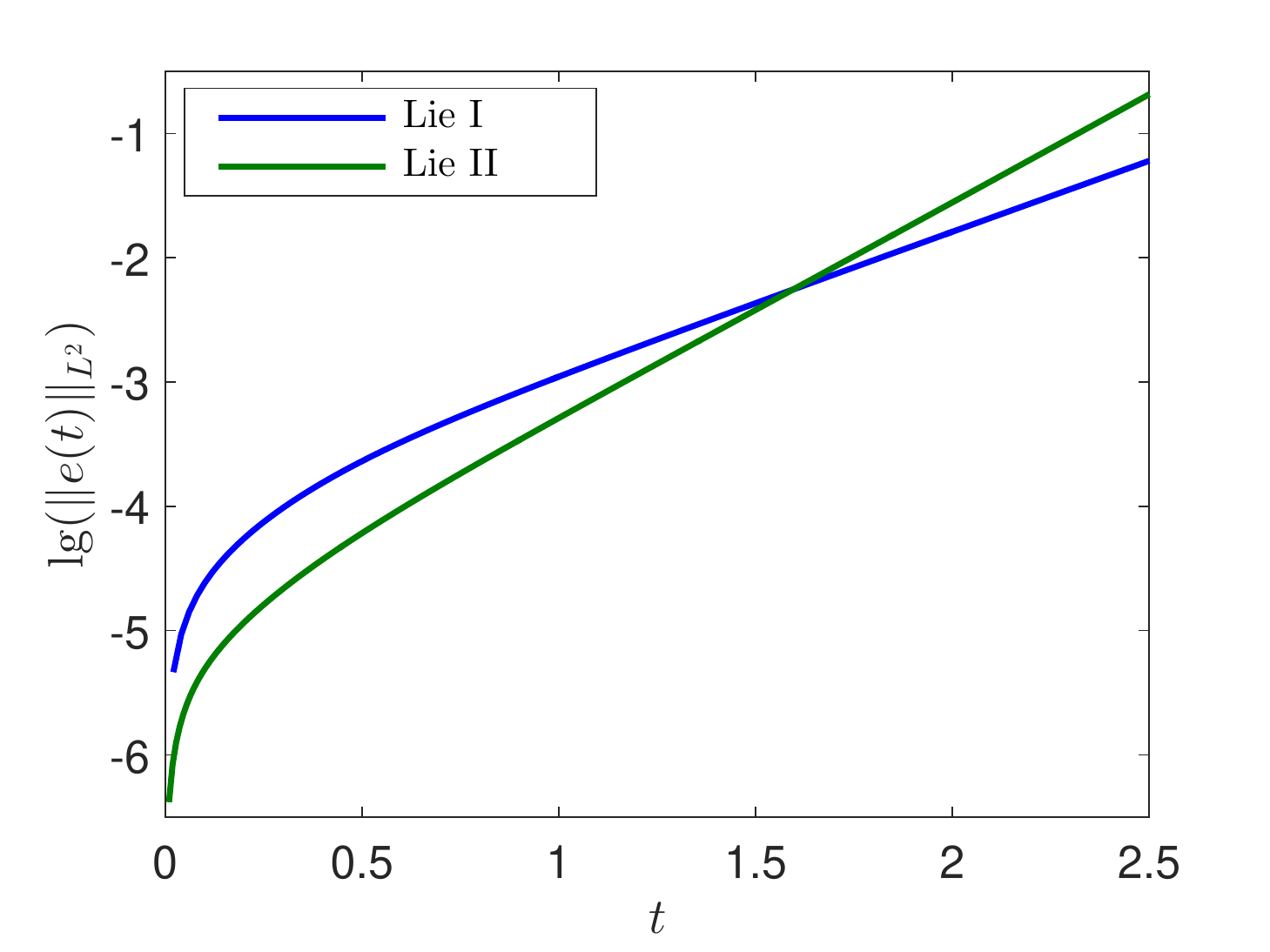}
\end{center}
\caption{Temporal errors of the Lie-Trotter splitting \eqref{Lie1} and \eqref{scheme2} for the solitary wave \eqref{solit}.}
\label{fig12}
\end{figure}

Fig. \ref{fig12} displays the temporal errors of the Lie-Trotter splitting method \eqref{Lie1} \& \eqref{scheme2} for the solitary wave. The left plot displays the errors at $t_N=T=2.5$ under different values of $N$. It can be clearly seen that via both cases of temporal grid, the Lie-Trotter splitting method converges at the first order in time. The right plot shows the evolution of the error when the step number $N$ is fixed as $N=25000$. We observe that the error increases exponentially in time for both cases and the error of Lie I increases a bit slower than that of Lie II, which suggests to choose Lie I for long time simulation.

\begin{figure}[htbp]
\begin{center}
\includegraphics[width=3.5in,height=1.8in]{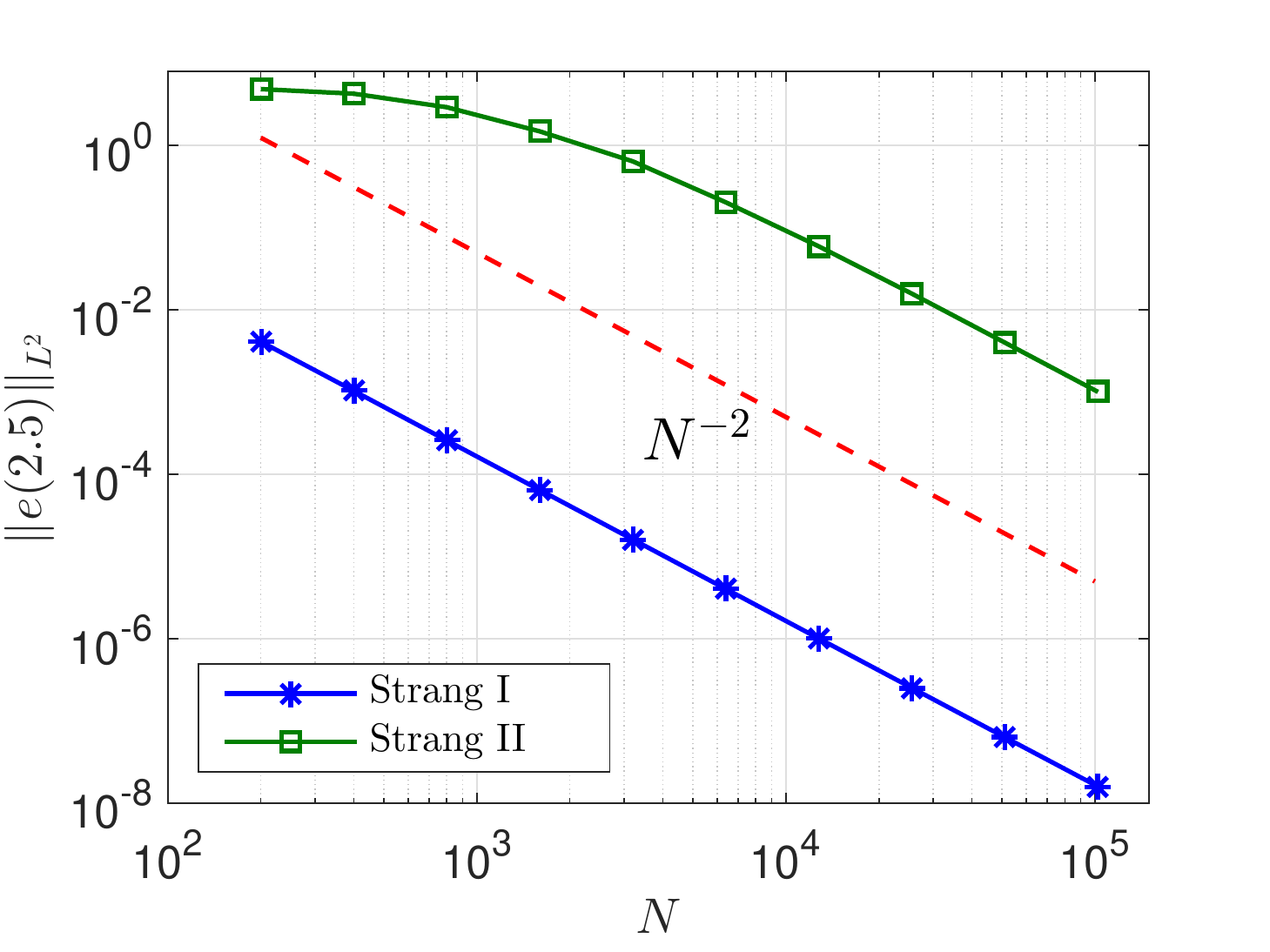}
\end{center}
\caption{Temporal errors of the Strang splitting scheme \eqref{scheme1}-\eqref{scheme2} for the solitary wave solution \eqref{solit} under different step numbers.}\label{fig1}
\end{figure}

\begin{figure}[htbp]
\begin{center}
\includegraphics[width=2.45in,height=1.5in]{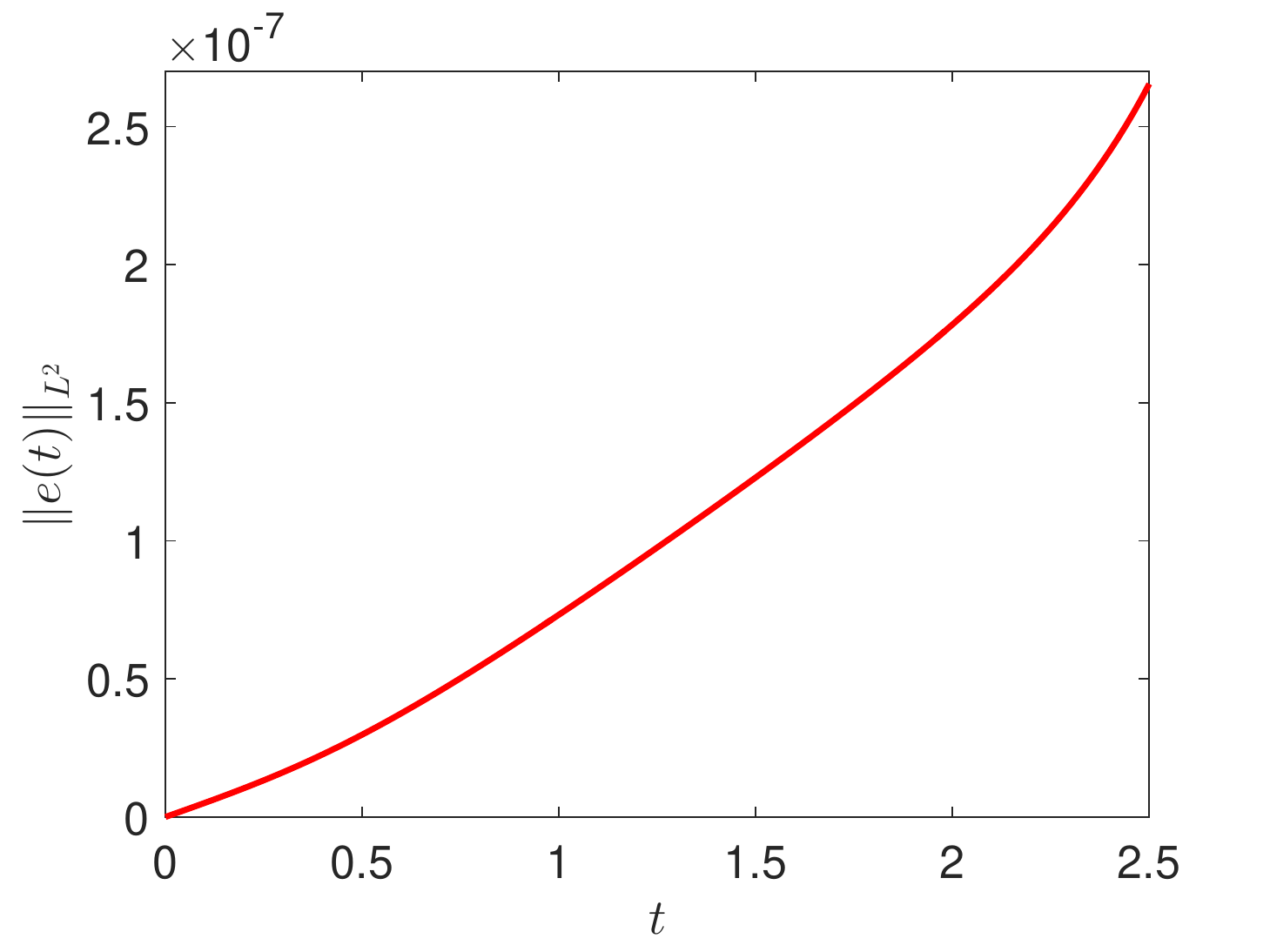}
\includegraphics[width=2.45in,height=1.5in]{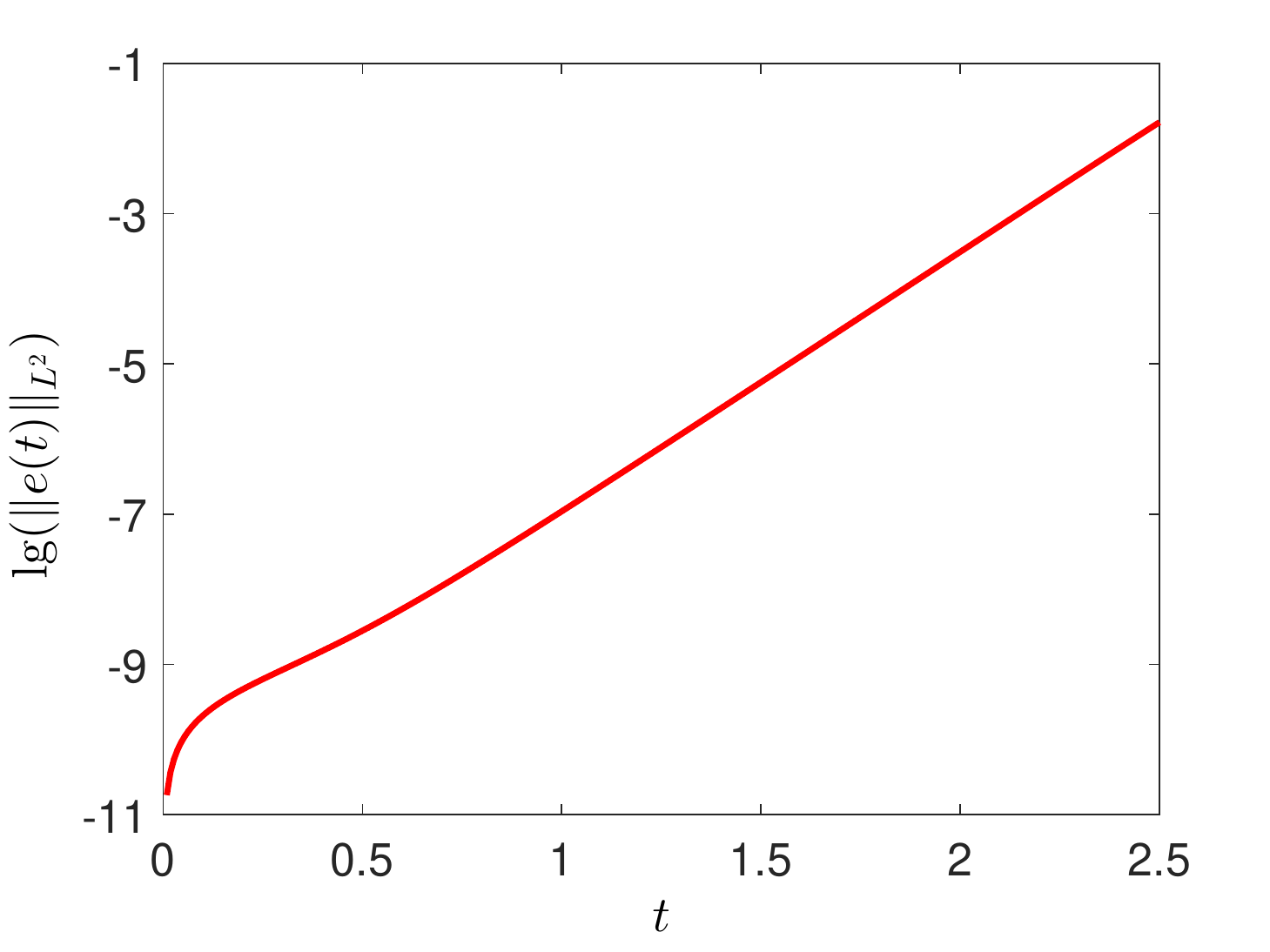}
\end{center}
\caption{Error evolution by using Strang I (left) and Strang II (right) for the solitary wave \eqref{solit}.}
\label{fig11}
\end{figure}

Fig. \ref{fig1} shows the error $\|e(2.5)\|_{L^2}$ for the methods
Strang I and Strang II. It can be clearly observed that the splitting
scheme converges at the second order in time and the error of Strang I
is far less than that of Strang II. Moreover, for the method Strang
II, it leads to a correct solution and is convergent quadratically
only when the step number $N$ is large enough $N\ge
N_T$. Fig. \ref{fig11} depicts the evolution of the error obtained via
Strang I and Strang II with fixed step number $N=25000$. We are
surprised to find that the error of Strang I method increases linearly
with respect to time while the error of Strang II increases
exponentially in time by noticing that the longitudinal axis
represents the error (left plot in Fig. \ref{fig11}) and the logarithm
of error (right plot in Fig. \ref{fig11}), respectively. The excellent
convergence behavior of Strang I is different from that of Lie I and
we can benefit a lot in practical computation, especially for long
time simulation. We postpone the analytical study of this method to a future work.

\begin{figure}[htbp]
\begin{center}
\includegraphics[width=2.45in,height=1.5in]{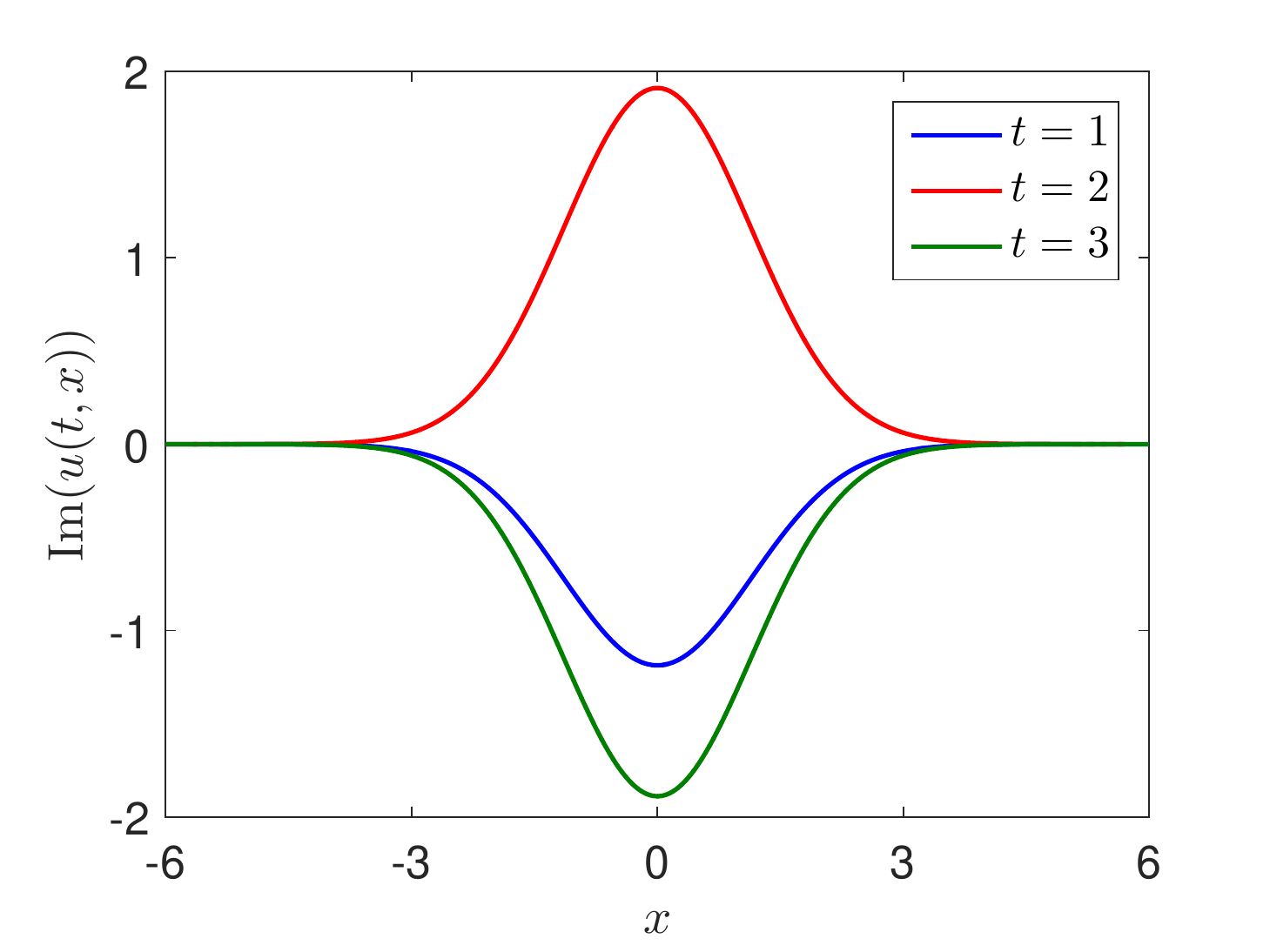}
\includegraphics[width=2.45in,height=1.5in]{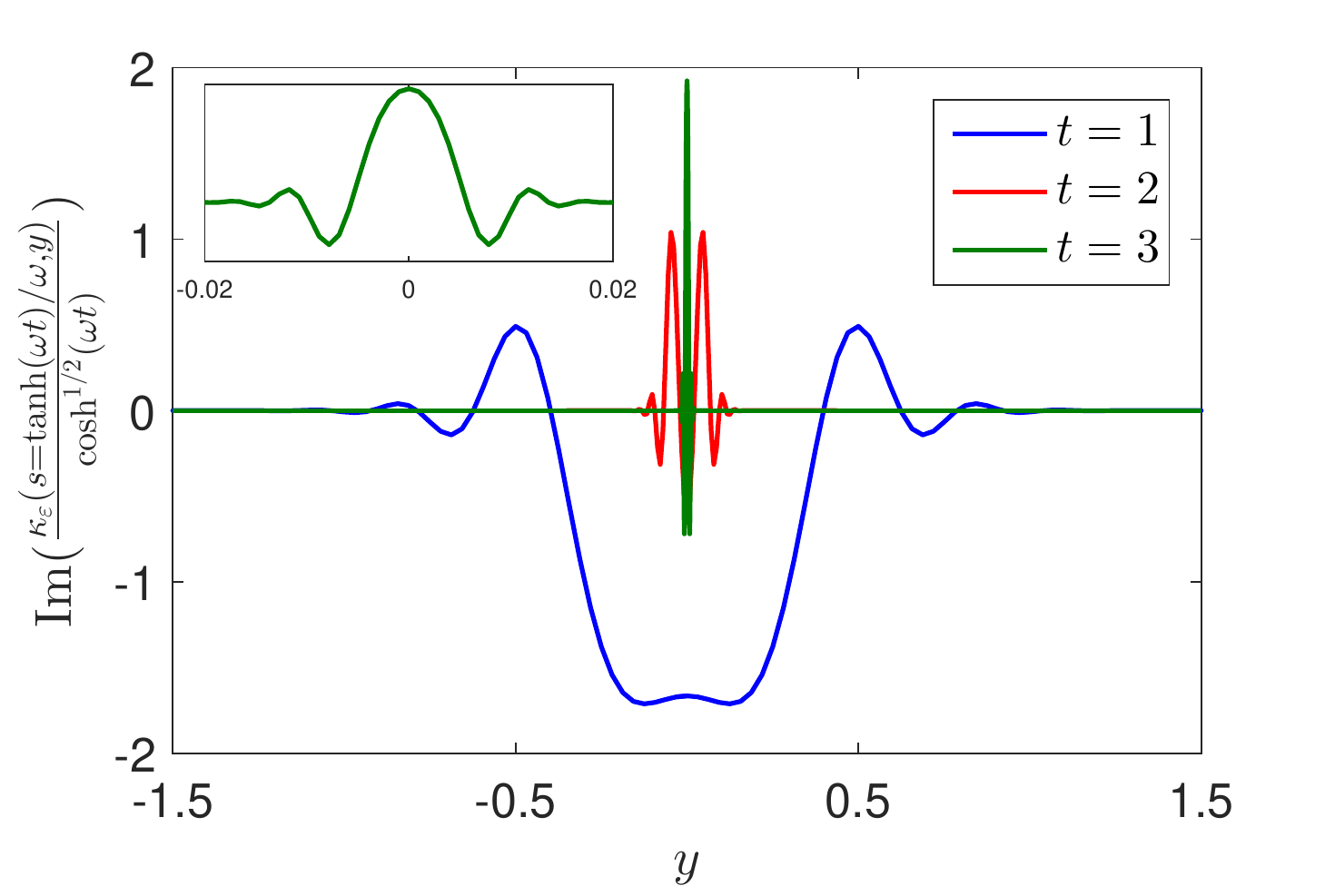}
\end{center}
\caption{Dynamics of $u(t,x)$ (left) and $\kappa(s,y)$ (right) for the solitary wave \eqref{solit}.}
\label{fig13}
\end{figure}

For the solitary wave, it is easily checked that direct classical
methods for $u$ work well since dispersion does not occur. On the contrary, the transform \eqref{uv} (and
\eqref{vg}) makes the ``essential" support of the solution $v$ (and
$\kappa$) shrink exponentially with respect to time, as is shown in
Fig. \ref{fig13}. This requires that the mesh size $h$ has to be set
tinier and tinier as time evolves in order to capture the correct
solution, which introduces huge costs if one keeps the computational
domain fixed. However, this difficulty can be overcome by shrinking
the computational region correspondingly in time, which enables the
computational cost to be comparable with the direct splitting methods
for \eqref{logNLSrep} (regularized version).

This suggests that a reasonable way to proceed may be to first simulate
\eqref{logNLSrep} directly: if the solution is not dispersive, then
nothing specific is needed (apart from the regularization of the
logarithm). If the solution is dispersive (which is always the case
when $\lambda>0$), then boundary effects
become strong, and it is more efficient to consider \eqref{rweq} and
use Strang I method via \eqref{scheme1}--\eqref{scheme2} on a nonuniform grid $\{s_n=\frac{\tanh(\omega t_n)}{\om}\}$.

\medskip
\noindent{\bf Example 2}.
We set $\lambda=-3$ and $\omega=2$ in the equation \eqref{logNLSrep}. We consider the following several cases:\\
(i). $u_0(x)=2e^{-\alpha x^2/2}$ for $\alpha=2$;\\
(ii). $u_0(x)=2e^{-\alpha x^2/2}$ for $\alpha=1/2$;\\
(iii). $u_0(x)=\sech(x^2/2)$.

\begin{figure}[h!]
\begin{center}
\includegraphics[width=3.5in,height=1.8in]{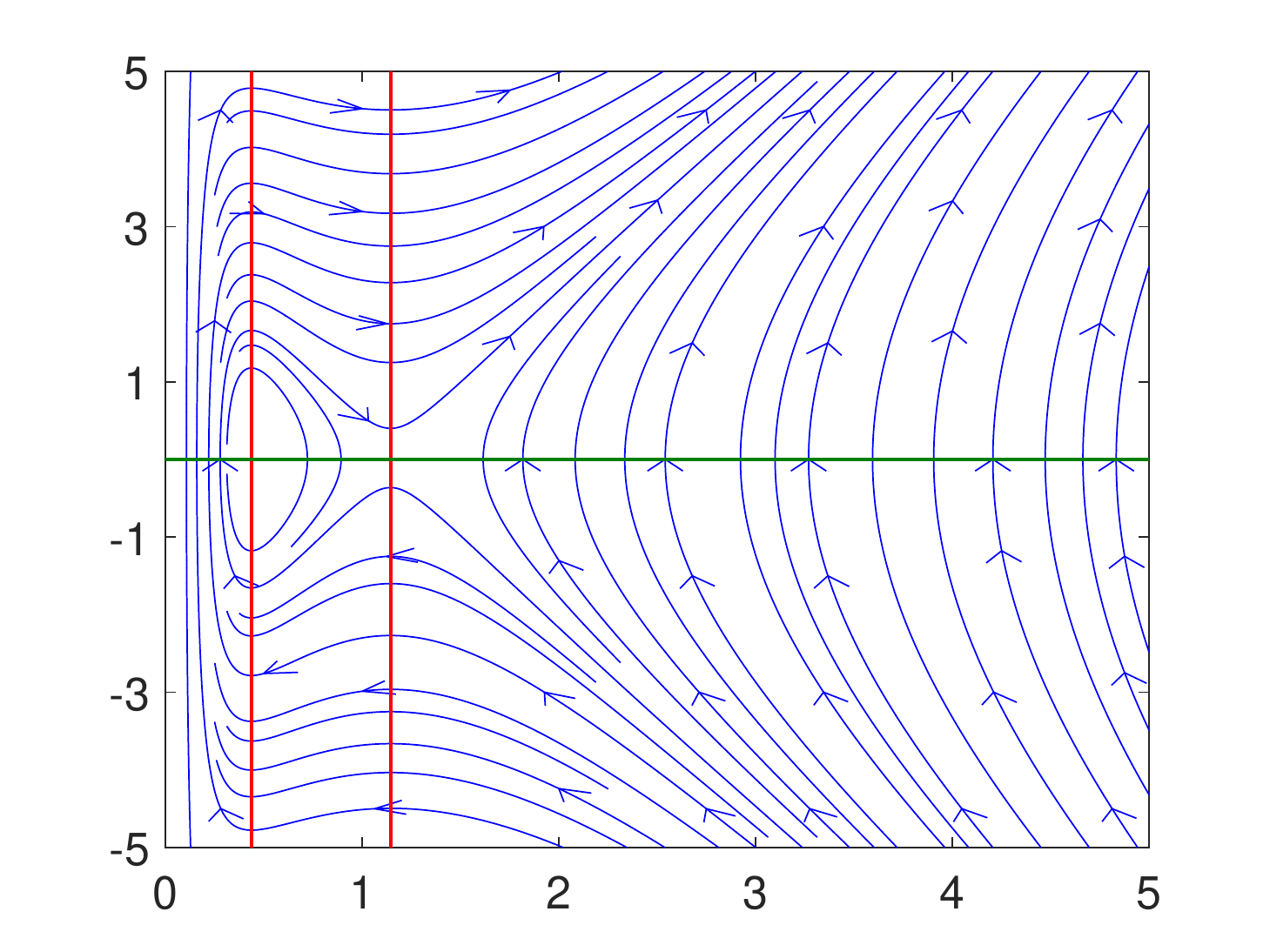}
\end{center}
\caption{Phase portraits for the ODE \eqref{eq:tau-general} with $\omega=2$ and $\lambda=-3$.}
\label{fig2}
\end{figure}

According to Proposition \ref{prop1}, for $\lambda=-3$, $\omega=2$, there exist two stationary solutions. The ODE \eqref{eq:tau-general}
\begin{equation*}
   \ddot \mu = \(\frac{1}{\mu^2}-k_-\)\(\frac{1}{\mu^2}-k_+\)\mu,\quad k_{\pm}=-\lambda\pm
  \sqrt{\lambda^2-\omega^2},
 \end{equation*}
 produces different trajectories due to different initial data for $\mu$, which is illustrated in Fig. \ref{fig2} displaying the phase portrait for the equation \eqref{eq:tau-general}.

\begin{figure}[htbp]
\begin{center}
\includegraphics[width=2.45in,height=1.5in]{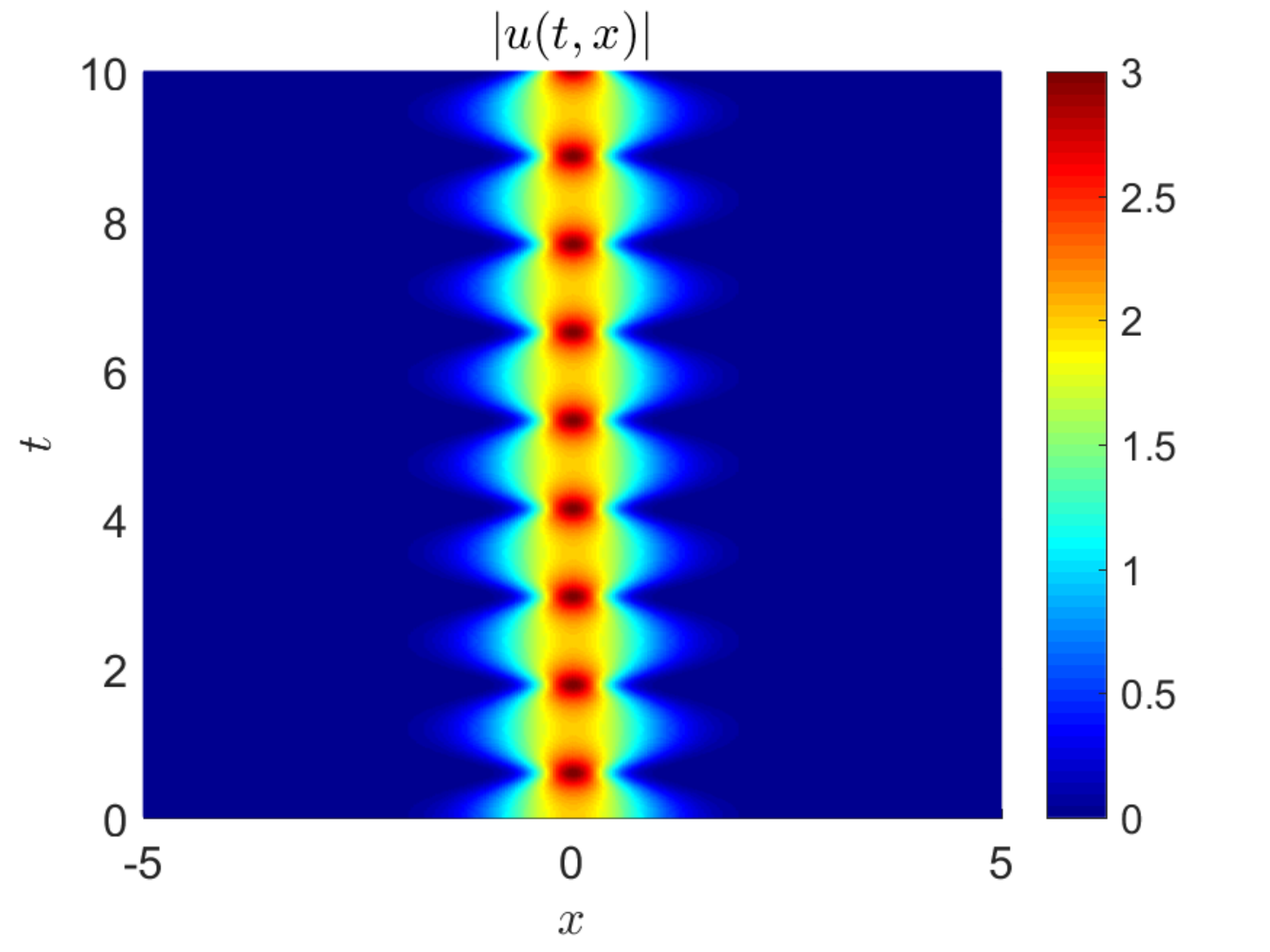}
\includegraphics[width=2.45in,height=1.5in]{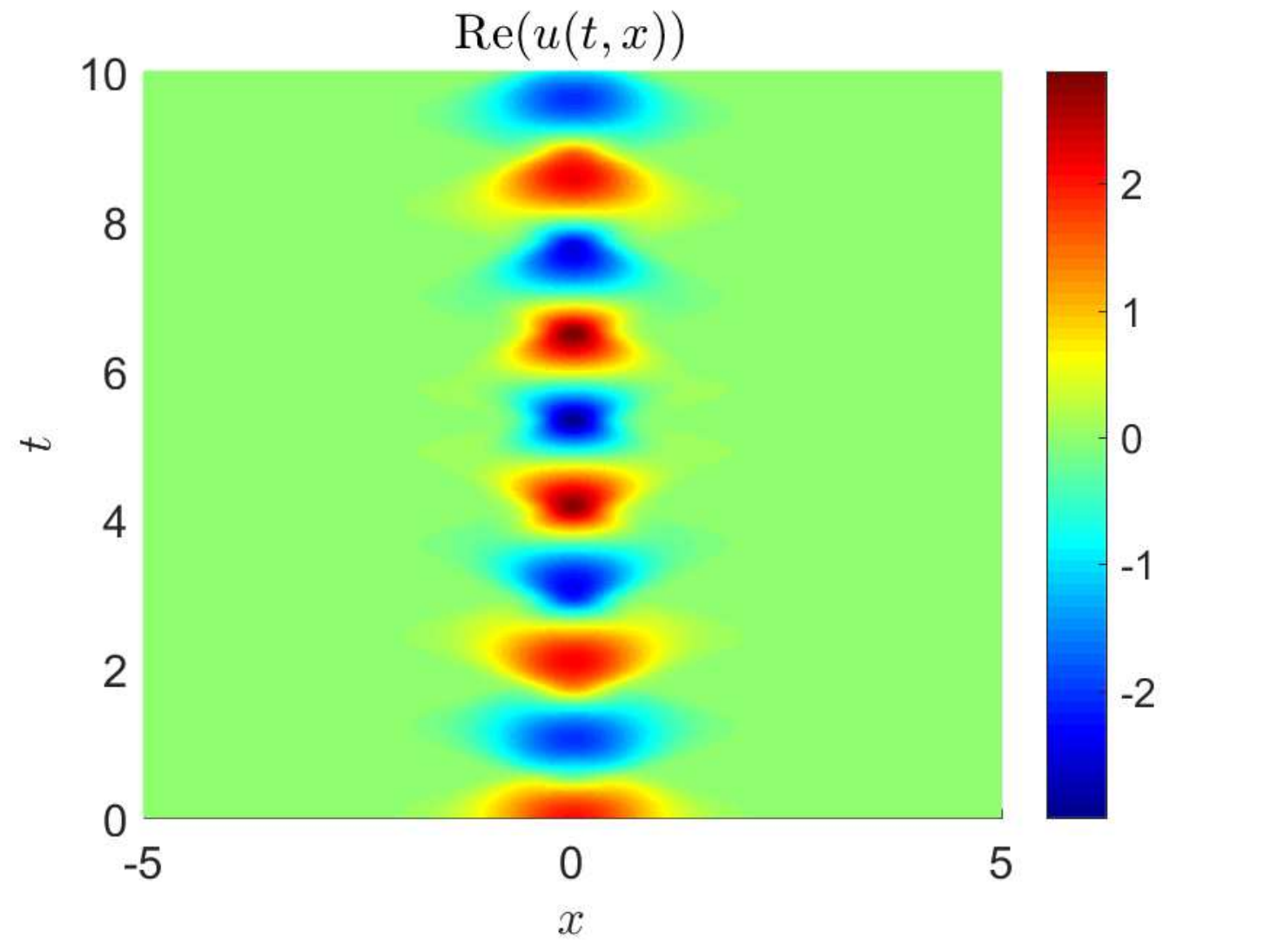}
\end{center}
\caption{Dynamics of $u$ for Case (i) initial data in Example 2.}
\label{fig3}
\end{figure}

For Case (i), $\alpha=2$ corresponds to $\mu(0)=1/\sqrt{\alpha}=1/\sqrt{2}$, $\dot{\mu}(0)=0$, which lies in the region producing time periodic solution $\mu(t)$, as seen from Fig. \ref{fig2}. This is confirmed by the dynamics shown in Fig. \ref{fig3}. We remark here that Fig. \ref{fig3} is obtained by the standard Strang splitting method for \eqref{logNLSrep} directly to avoid necessary treatments including adaptive mesh refinement and computational domain cut, since there is no dispersion for $u$.

\begin{figure}[htbp]
\begin{center}
\includegraphics[width=2.45in,height=1.5in]{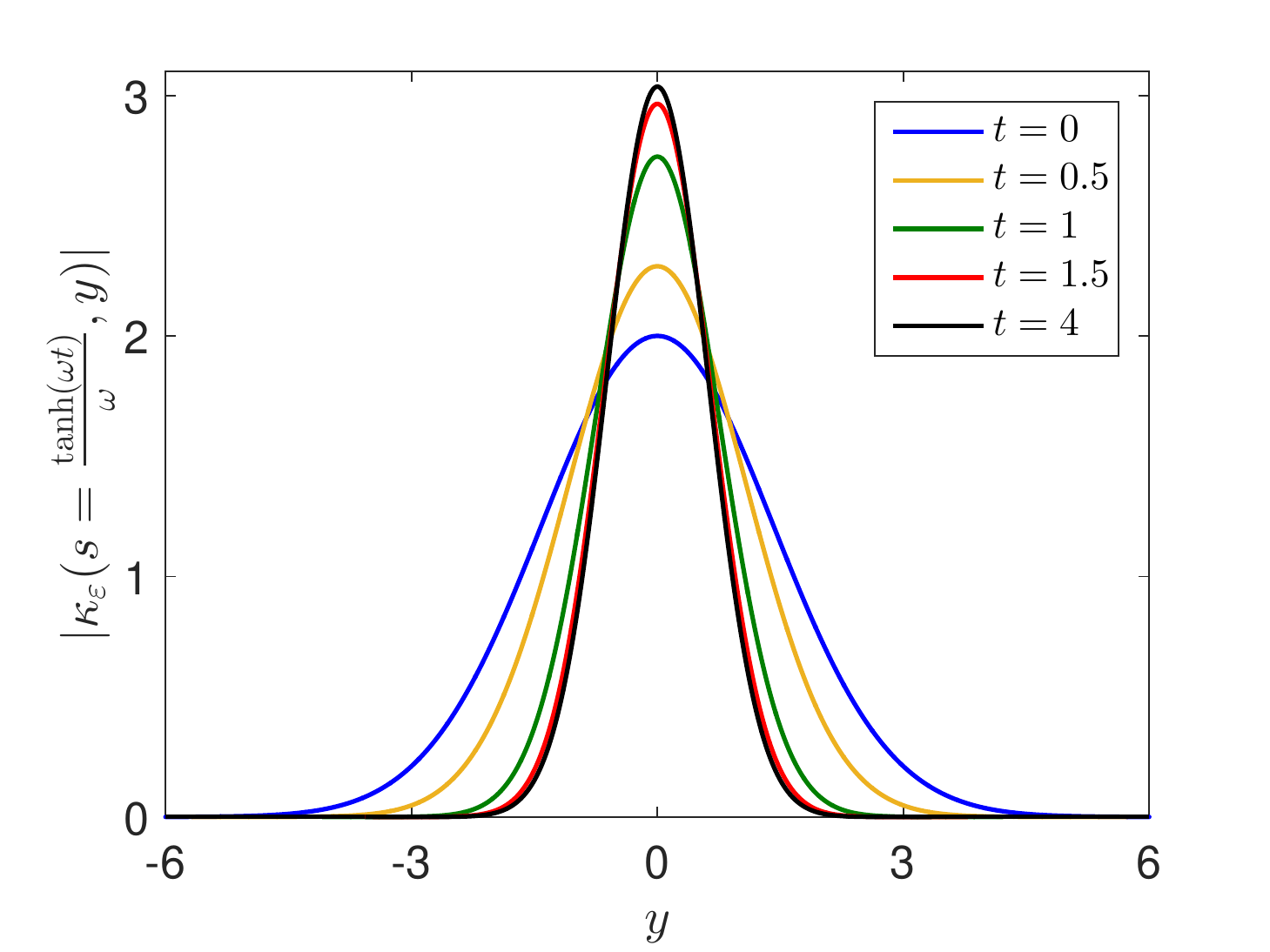}
\includegraphics[width=2.45in,height=1.5in]{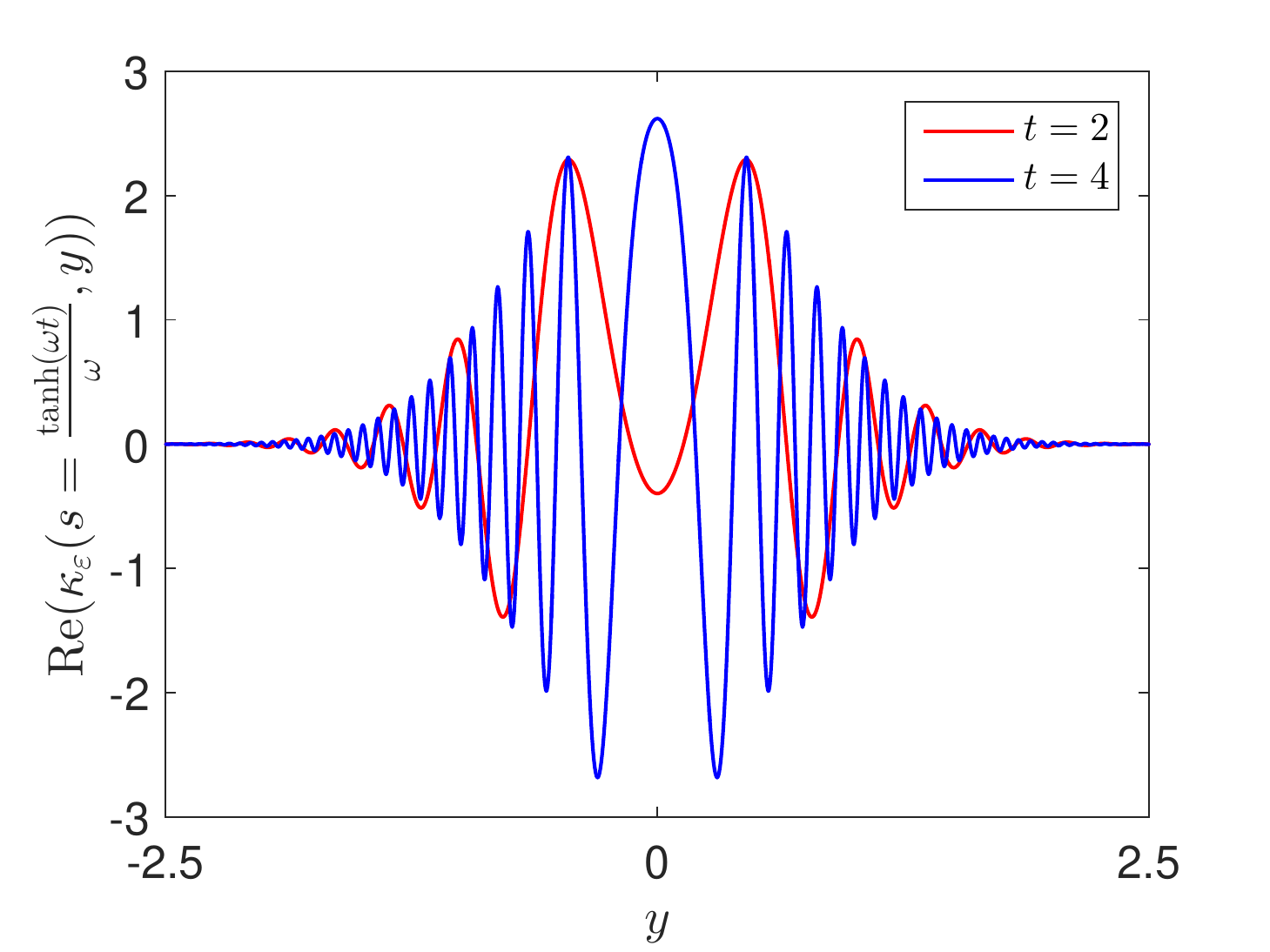}\\
\includegraphics[width=2.45in,height=1.5in]{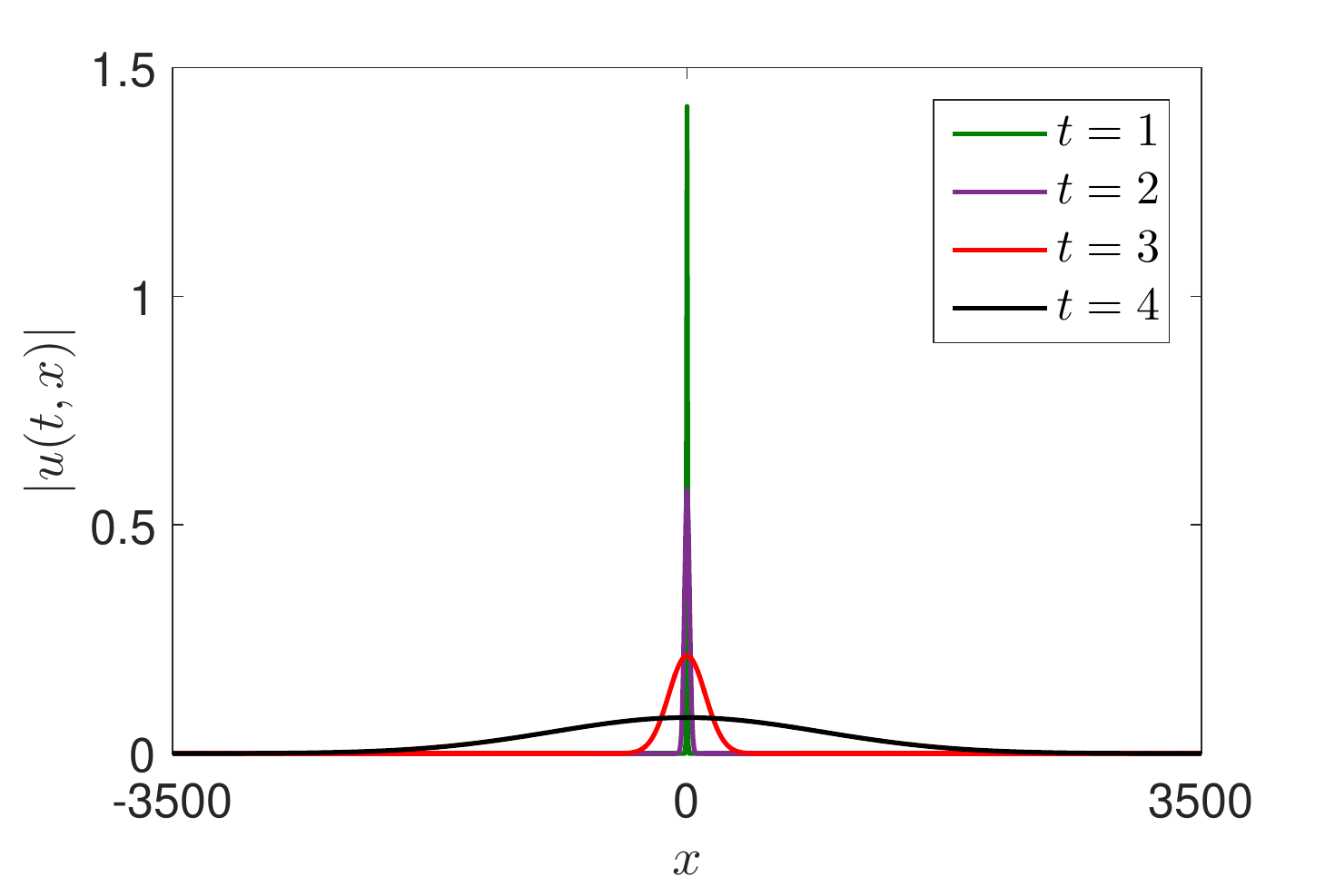}
\includegraphics[width=2.45in,height=1.5in]{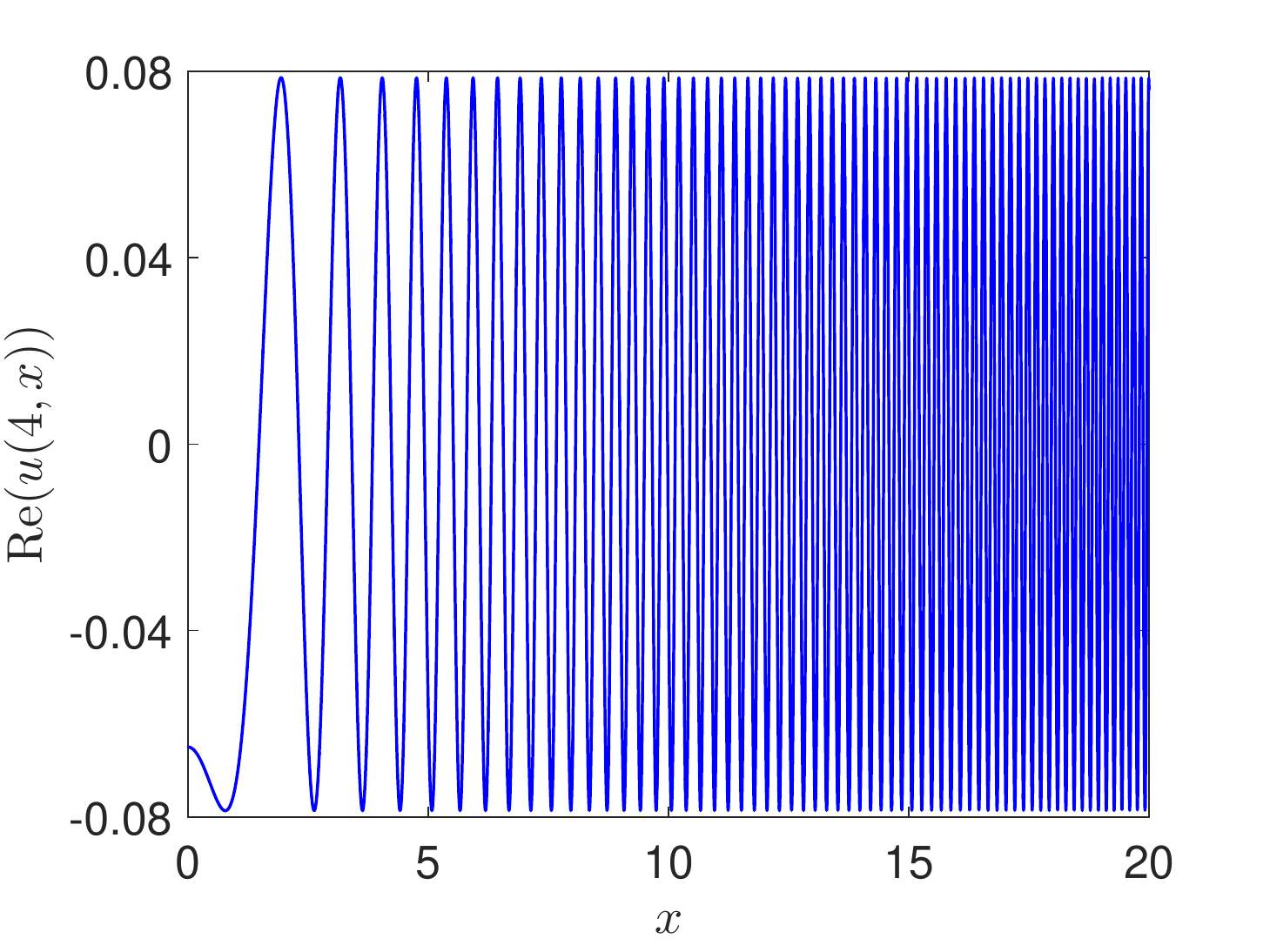}
\end{center}
\caption{Dynamics of $\kappa$ (top) and $u$ (bottom) for Case (ii) initial data in Example 2.}
\label{fig4}
\end{figure}

Fig. \ref{fig4} shows the dynamics of $\kappa_\varepsilon(s, y)$ for Case (ii) initial data by using Strang I method till $T=4$, or equivalently
$s=\tanh(4\omega)/\omega$  with time step $\tau=0.0001$ and mesh size
$h=1/2^{10}$. It can be observed that $\kappa_\varepsilon$, or accordingly $\kappa$, is well localized in
space and $|\kappa(s,y)|$ stays almost invariant after some time
(cf. the top left one in Fig. \ref{fig4}), however, the argument of
$\kappa$ varies in time (the top right one in Fig.
\ref{fig4}). Moreover, $\kappa$ oscillates in space. The bottom part
in Fig. \ref{fig4} shows the profile of $u(t,x)$ obtained by
\eqref{scheme2}, which shows the dispersion of $u$. The right one
displays $\RE(u(4,x))$, which shows that $u$ is highly oscillatory in
space. The oscillation is more and more drastic when observed
further away from the origin, which agrees with the transformation
\eqref{scheme2}, due to the term
$\exp\(i\frac{\omega}{2}x^2\tanh(\omega t_n)\)$.
Because of the rapid expansion and  drastic oscillation in space, it is a disaster if one solves
\eqref{logNLSrep} directly.

\begin{figure}[htbp]
\begin{center}
\includegraphics[width=2.45in,height=1.5in]{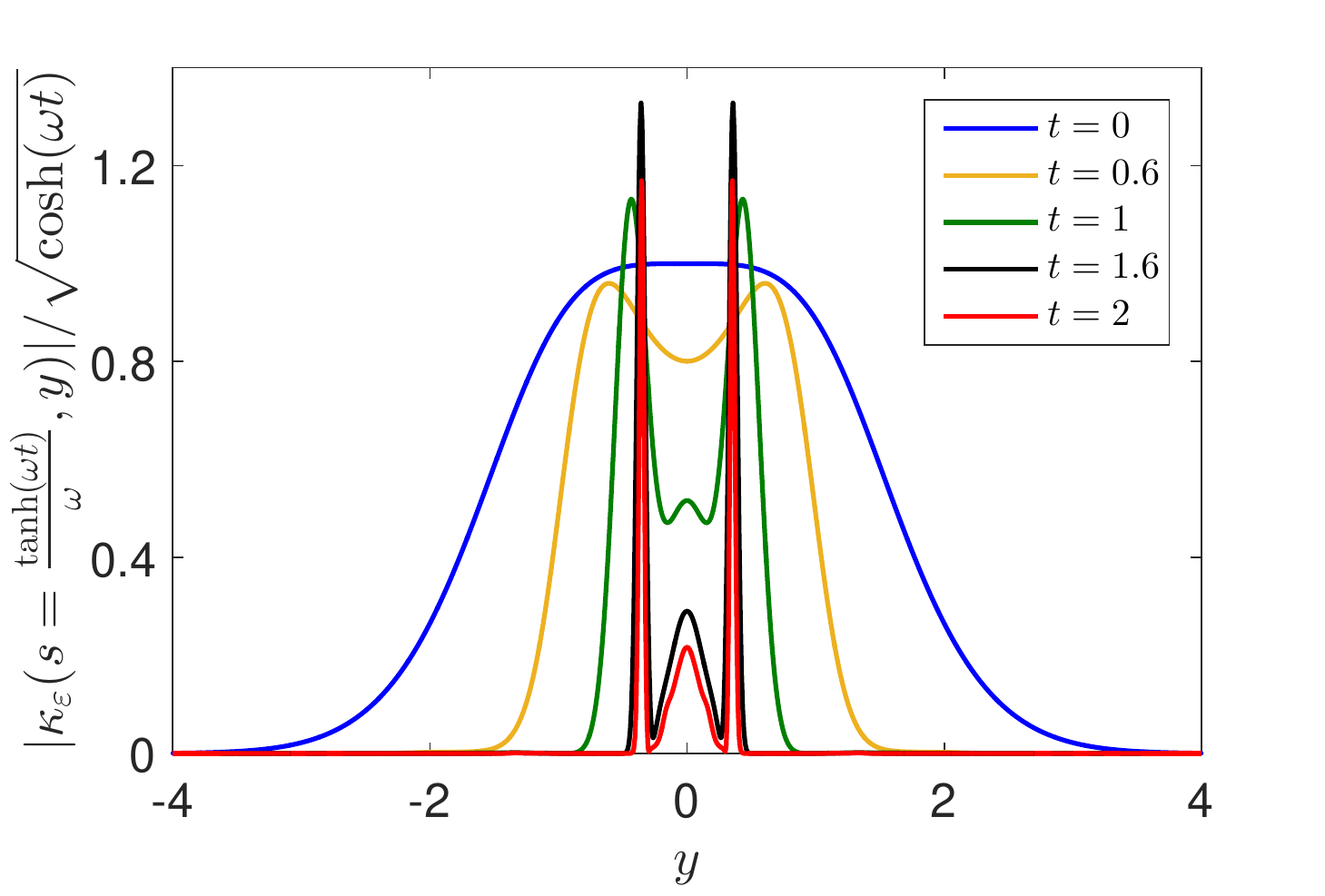}
\includegraphics[width=2.45in,height=1.5in]{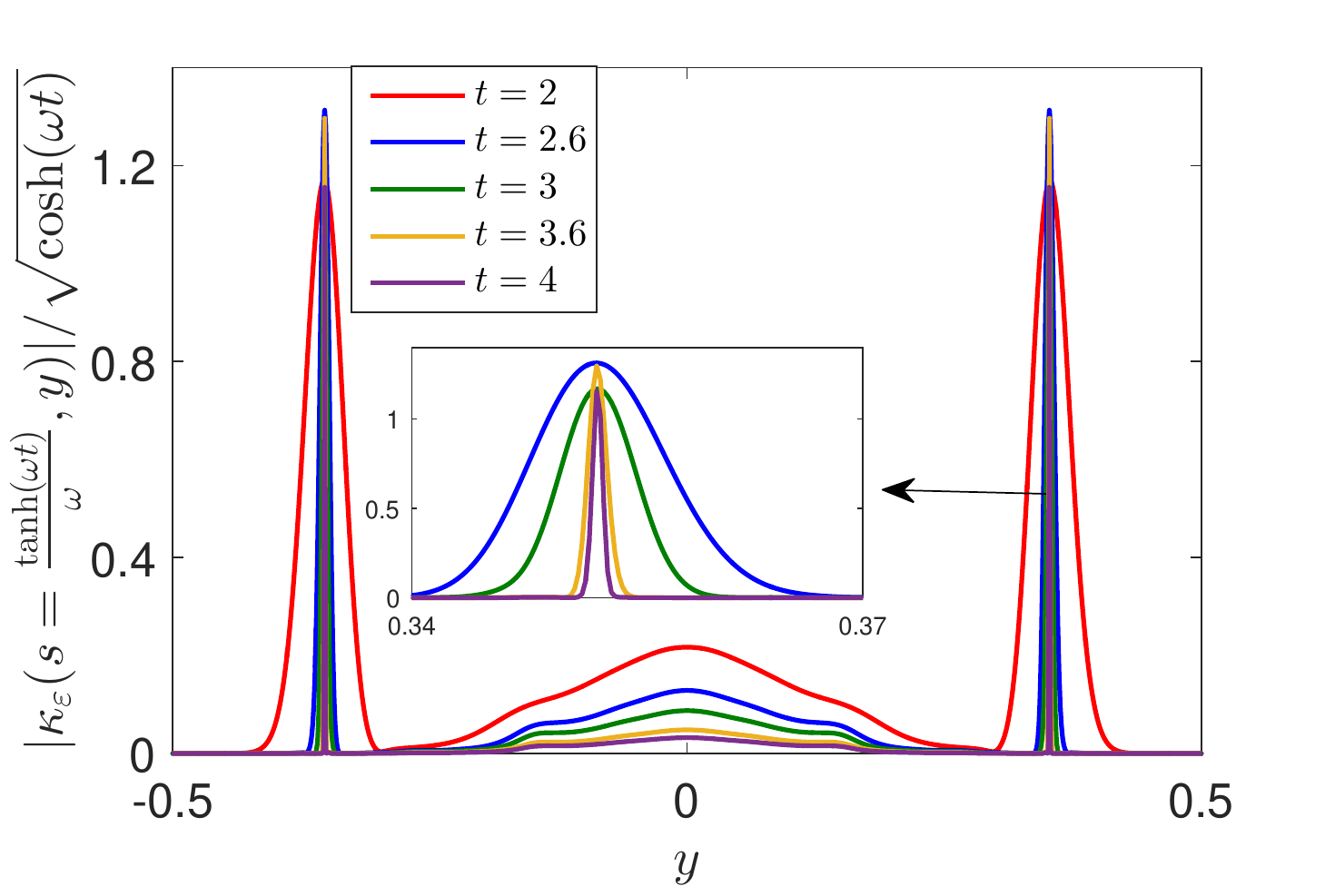}\\
\includegraphics[width=2.45in,height=1.5in]{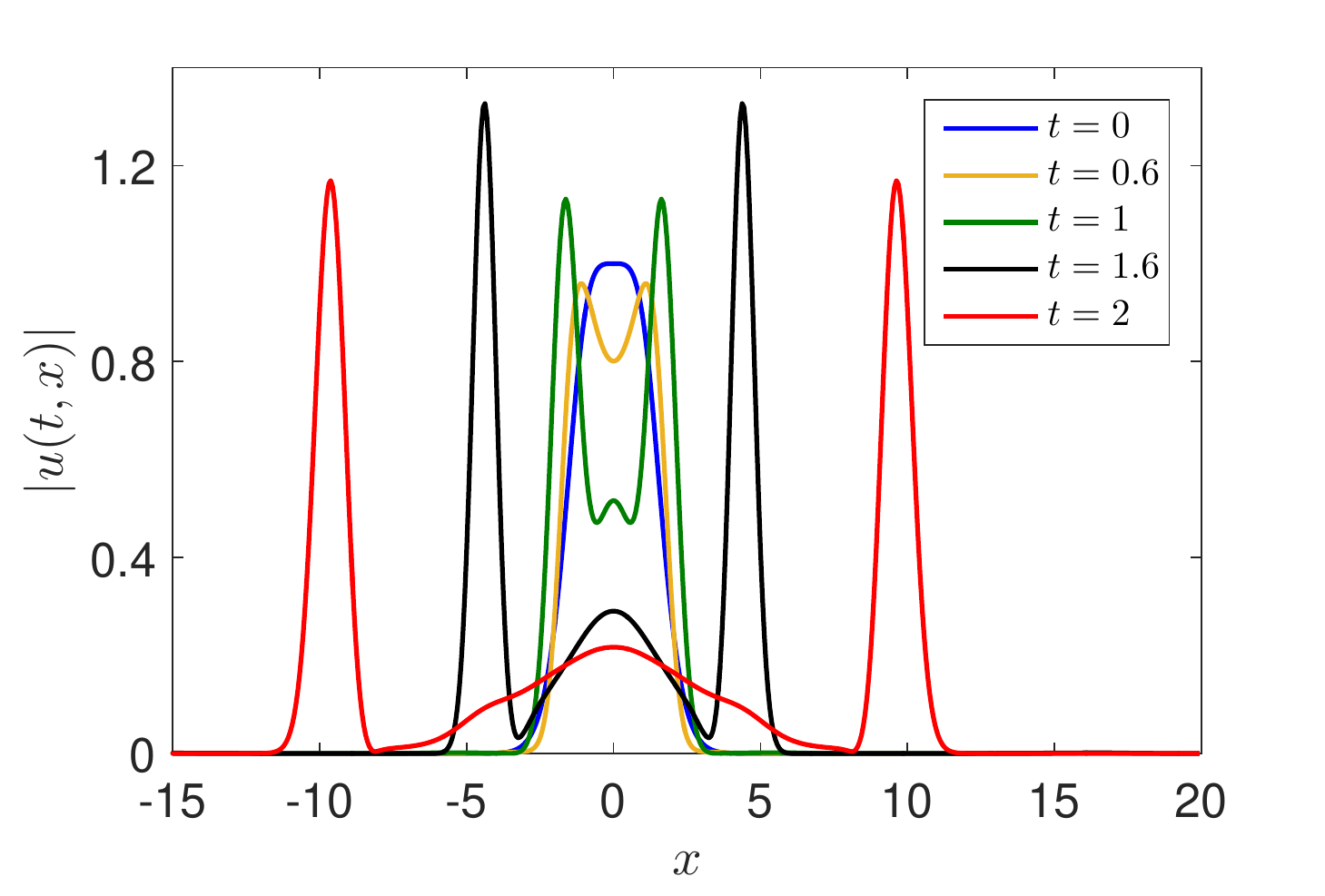}
\includegraphics[width=2.45in,height=1.5in]{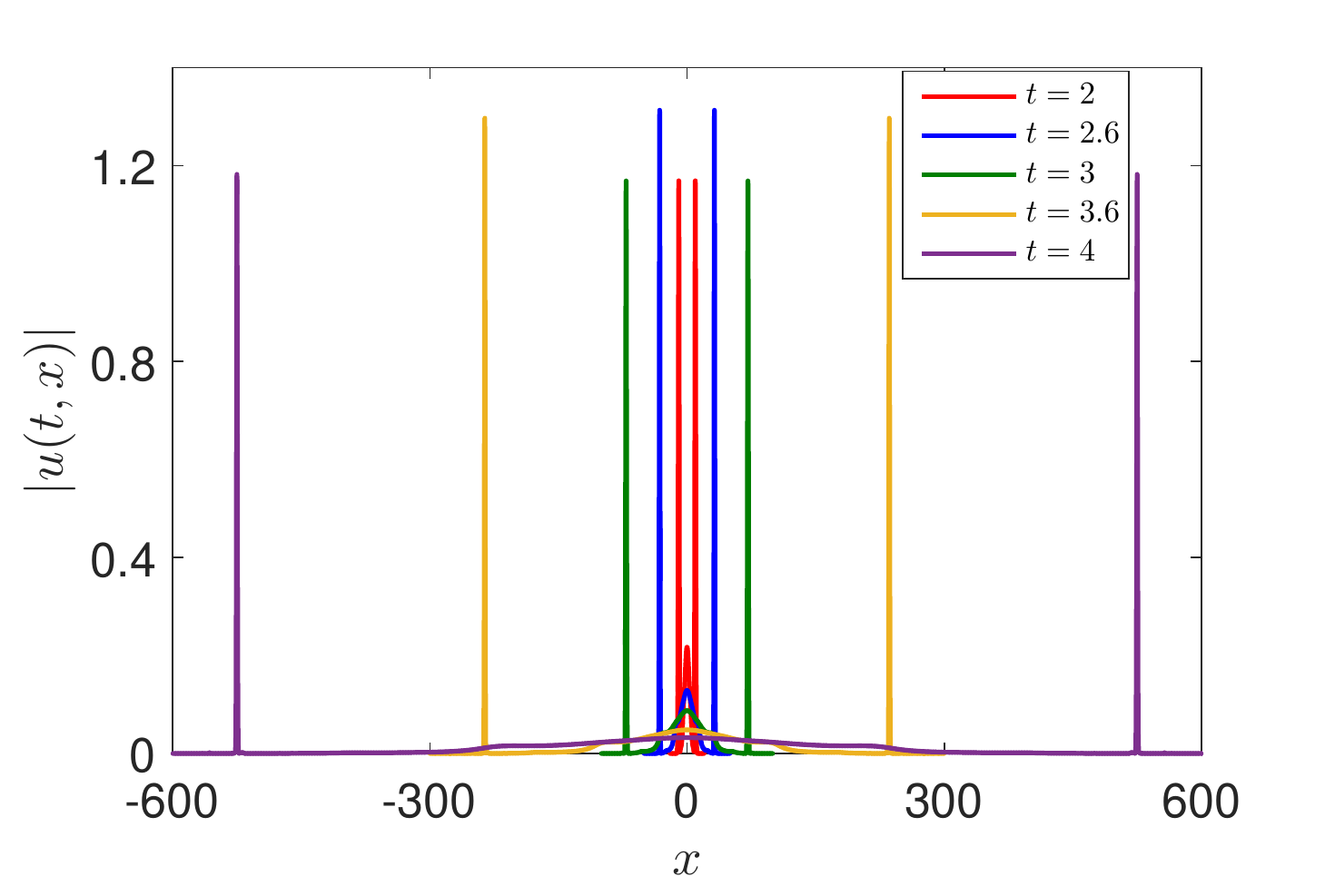}
\end{center}
\caption{Dynamics of $\kappa(s,y)$ (top) and $u(t,x)$ (bottom) for Case (iii) in Example 2.}
\label{fig23}
\end{figure}

\begin{figure}[htbp]
\begin{center}
\includegraphics[width=2.45in,height=1.5in]{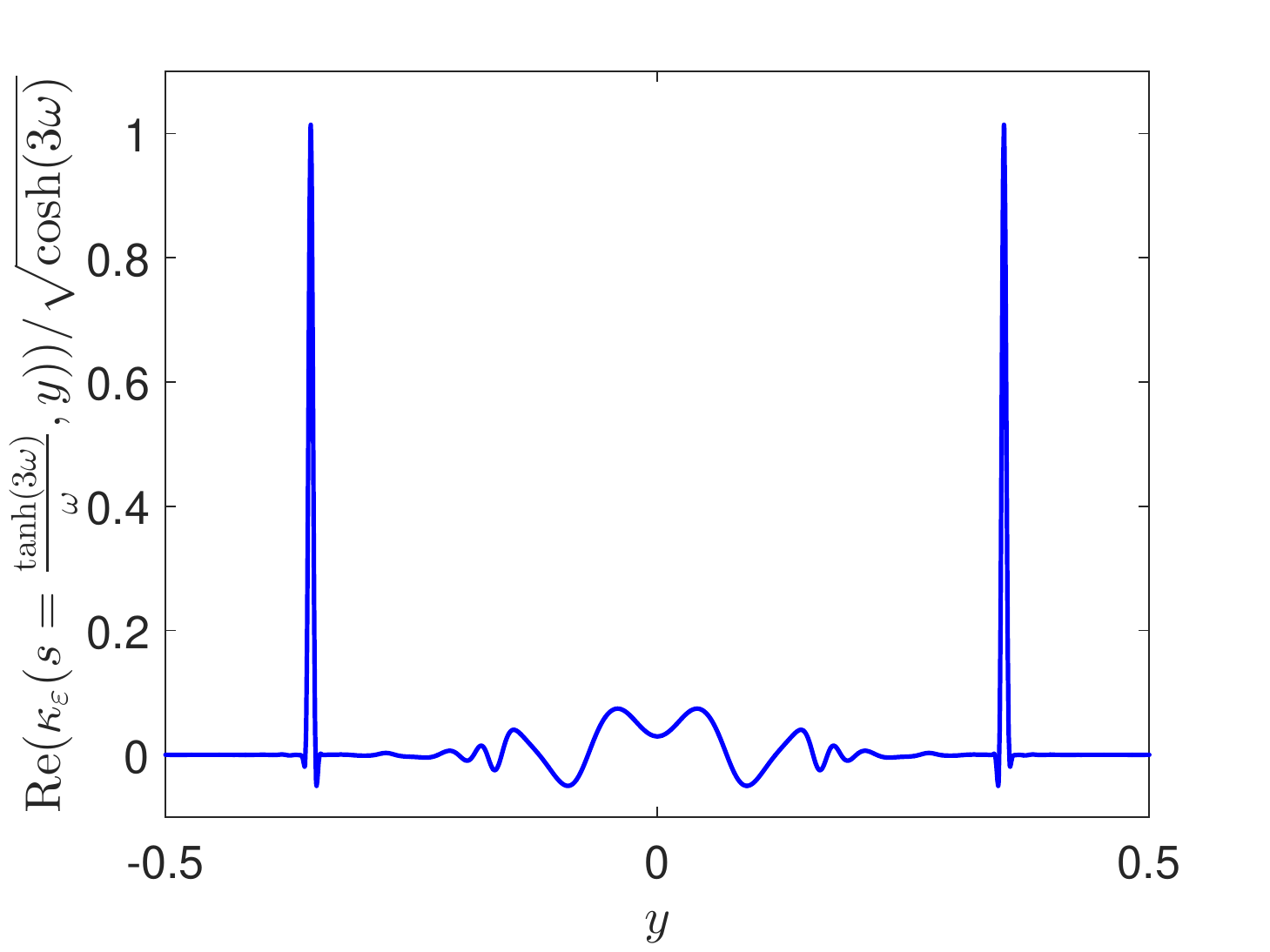}
\includegraphics[width=2.45in,height=1.5in]{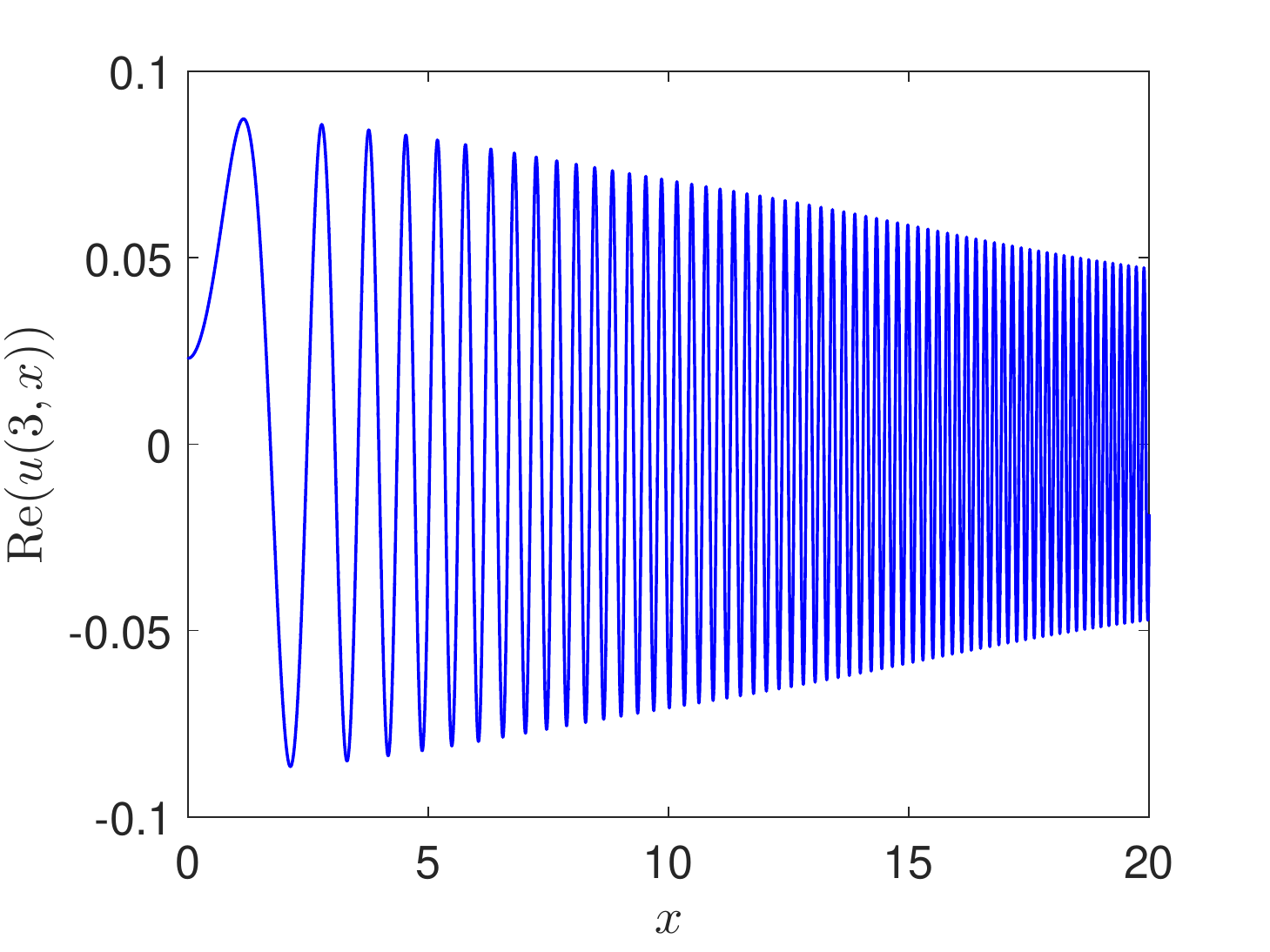}
\end{center}
\caption{The profile of
  ($\RE \kappa_\varepsilon(s=\frac{\tanh(3\omega)}{\omega},y)$) (left)
  and $\RE(u(3,x))$ (right) for Case (iii) in Example 2.}
\label{fig231}
\end{figure}

Fig. \ref{fig23} shows the profile of $\kappa(s,y)$ for Case (iii) initial data. Different from Case (ii) where the modulus of
$\kappa$ keeps invariant after some time, here the modulus keeps
increasing with respect to time, which is illustrated in the top
graphics. The wave is split into two symmetric wave packets and a wave in
between. As time evolves, the width of the two wave packets gets narrower
and narrower. With the formulation \eqref{scheme2}, we get the
dynamics of $u(t,x)$ in the bottom in Fig. \ref{fig23}. We observe
that the initial wave is split into two symmetric breathers which
spread out quickly (exponentially) in space with modulus varying
periodically in time and a dispersive wave in between. Fig.
\ref{fig231} displays the profile of $\RE (\kappa)$ and $\RE(u)$ at $t=3$, which suggests the feasibility of solving $u$ via the formulation \eqref{scheme2} and the virtual difficulty of solving $u$ directly.

\medskip

\begin{figure}[htbp]
\begin{center}
\includegraphics[width=3.5in,height=1.8in]{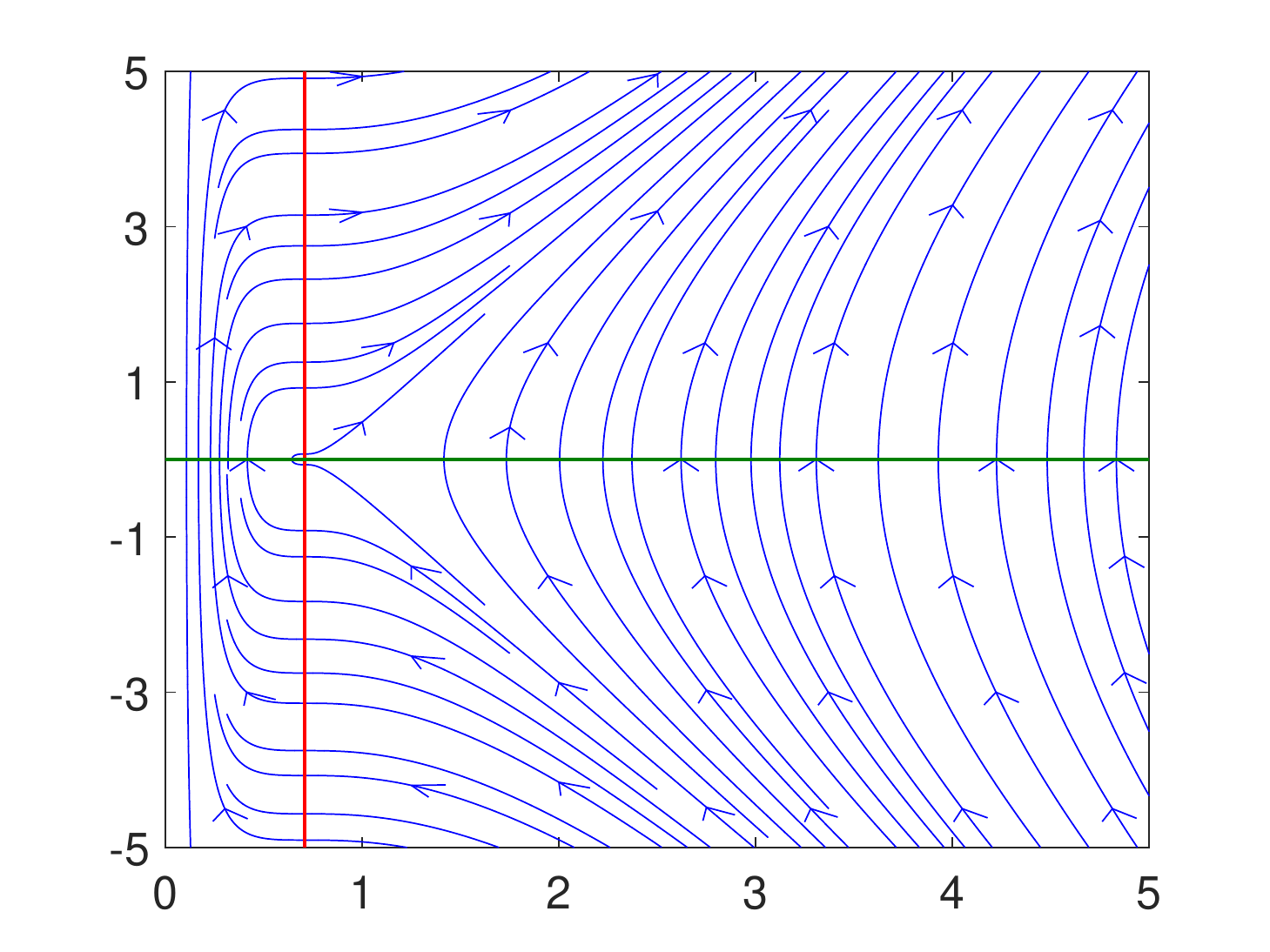}
\end{center}
\caption{Phase portraits for the ODE \eqref{eq:tau-general} with $\omega=2$ and $\lambda=-2$.}
\label{portrait2}
\end{figure}

\noindent{\bf Example 3}.
We set $\lambda=-2$ and $\omega=2$ in the equation \eqref{logNLSrep}, in which case there exists a continuous family of solitary wave when $u_0(x)=Ae^{-\alpha x^2/2}$ for $\alpha=2$ (or, equivalently $\mu(0)=1/\sqrt{2}$), according to Proposition \ref{prop1}, that can be illustrated from the phase portrait for the ODE \eqref{eq:tau-general} (cf. Fig. \ref{portrait2}). We consider the following three cases for the initial data $u_0$:\\
(i). $u_0(x)=2e^{-\alpha x^2/2}$ for $\alpha=3$;\\
(ii). $u_0(x)=2e^{-\alpha x^2/2}$ for $\alpha=1/2$;\\
(iii). $u_0(x)=\sech(x^2/2)$.

\begin{figure}[htbp]
\begin{center}
\includegraphics[width=2.4in,height=1.5in]{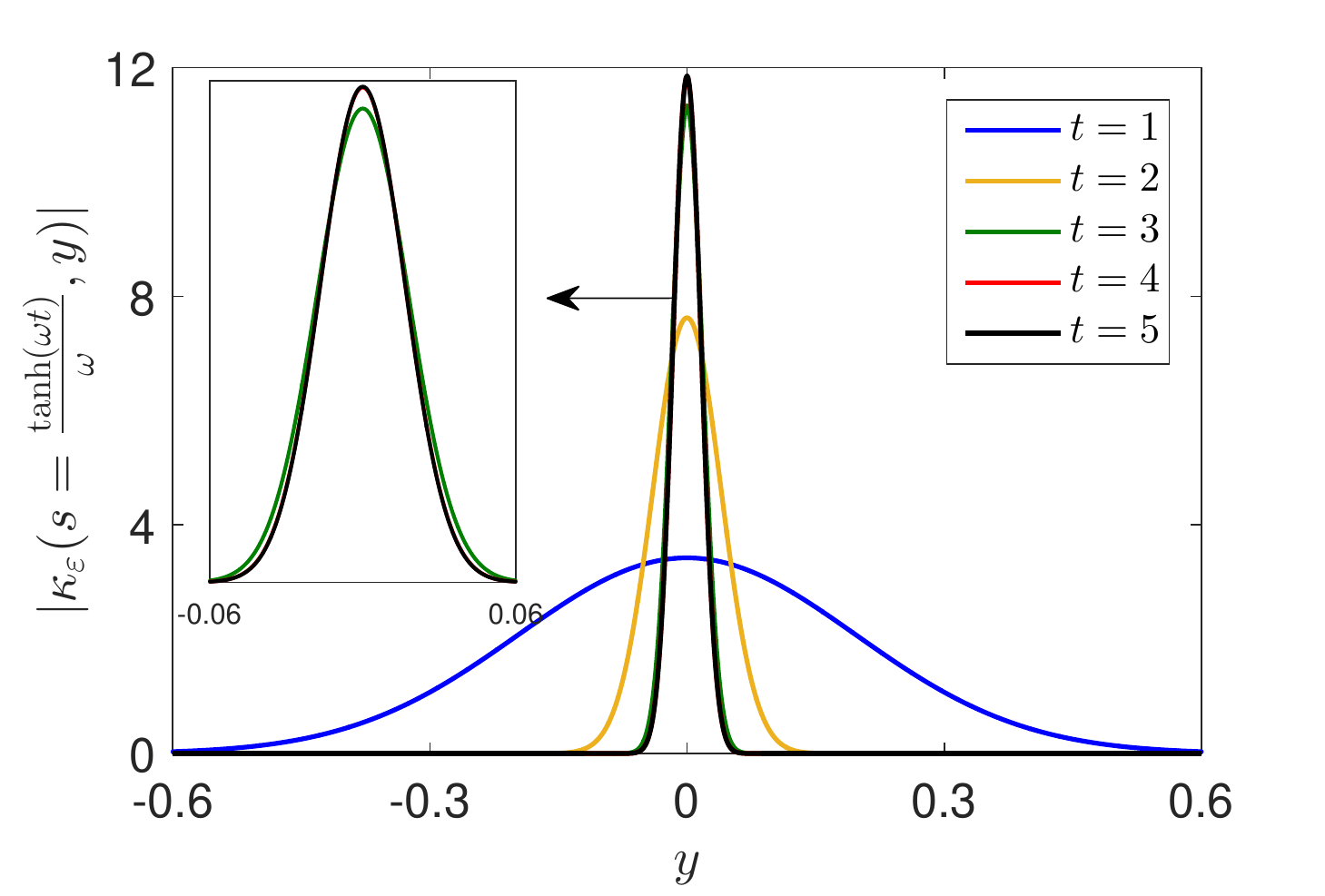}
\includegraphics[width=2.4in,height=1.5in]{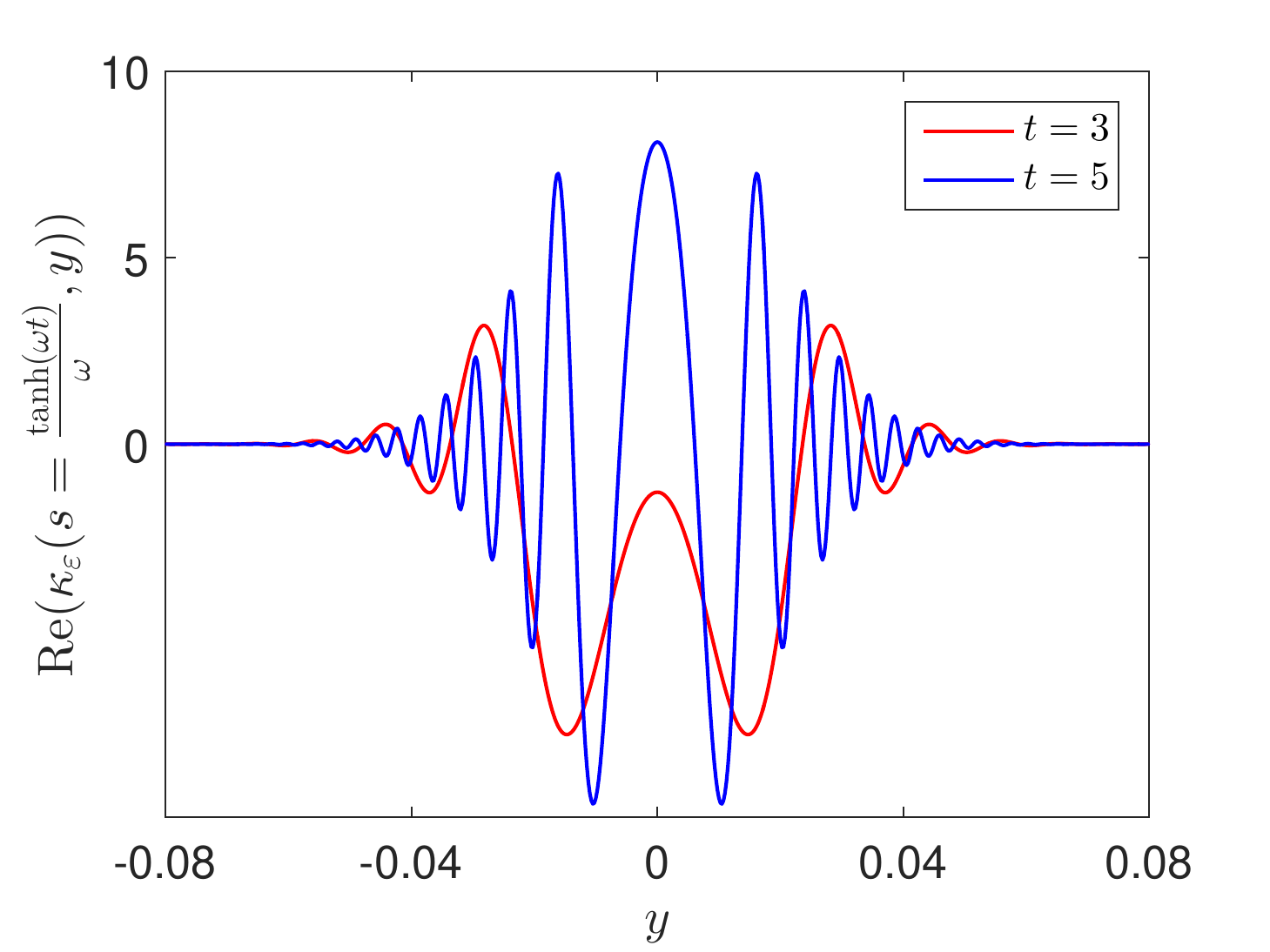}\\
\includegraphics[width=2.4in,height=1.5in]{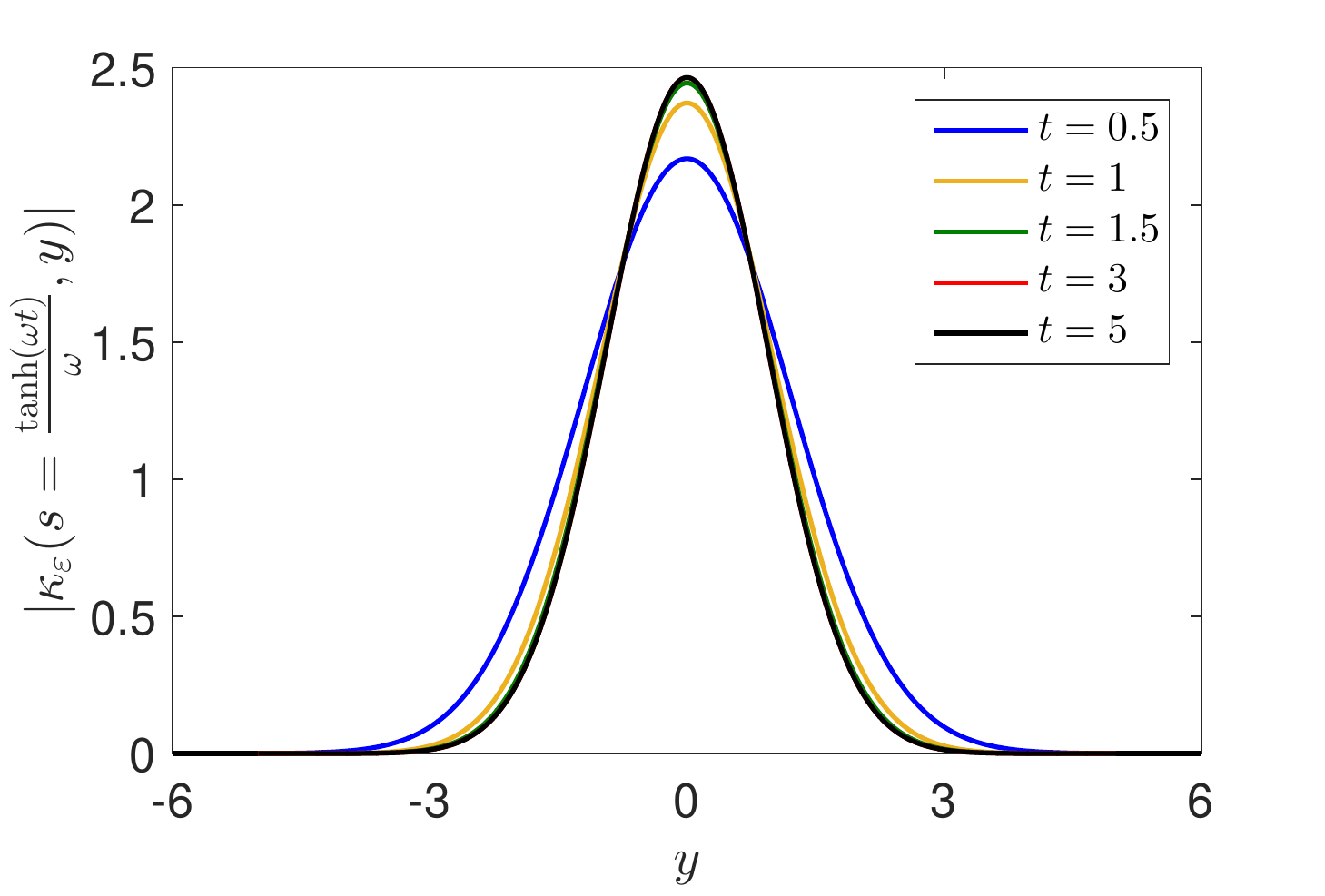}
\includegraphics[width=2.4in,height=1.5in]{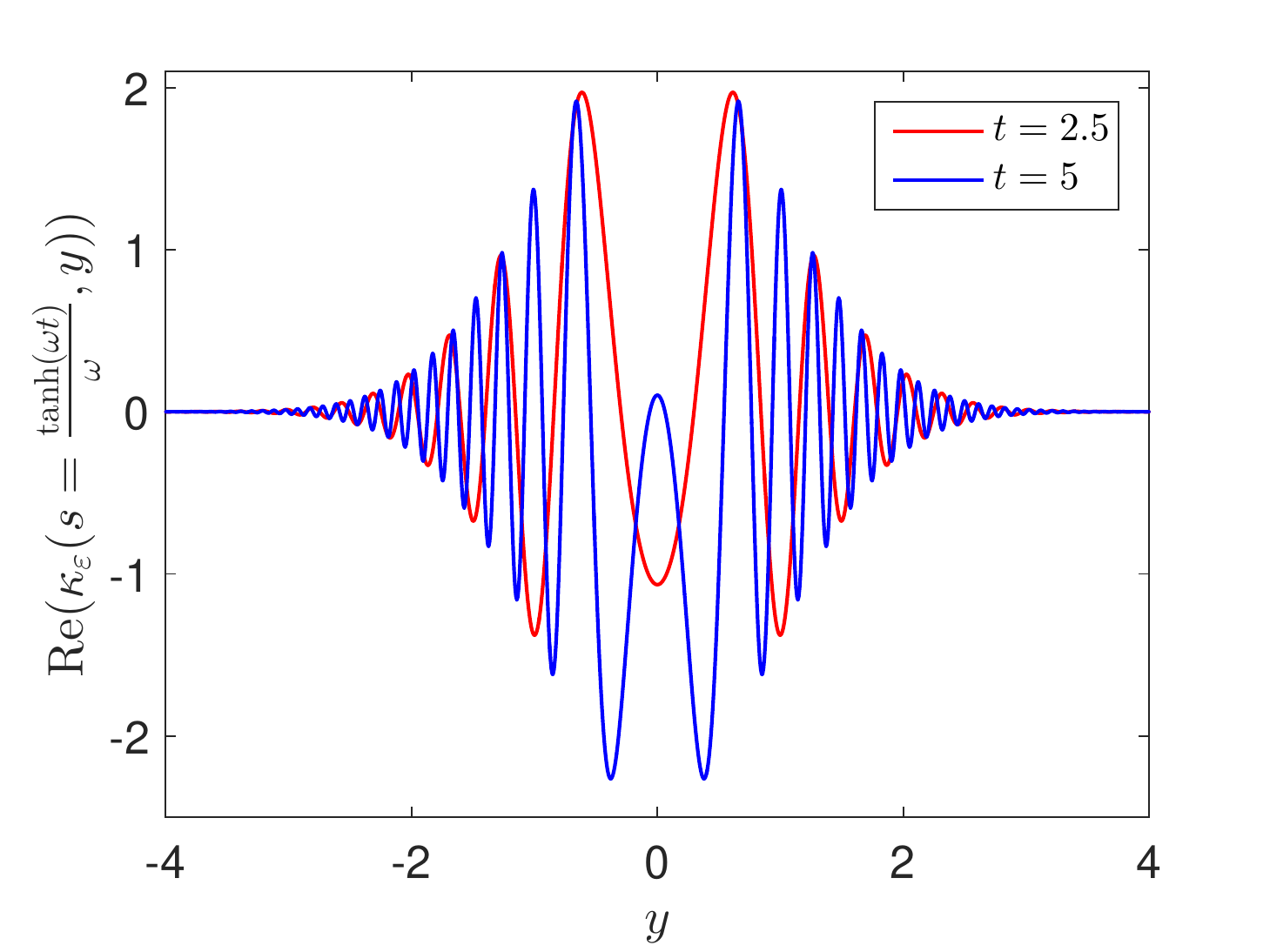}\\
\includegraphics[width=2.4in,height=1.5in]{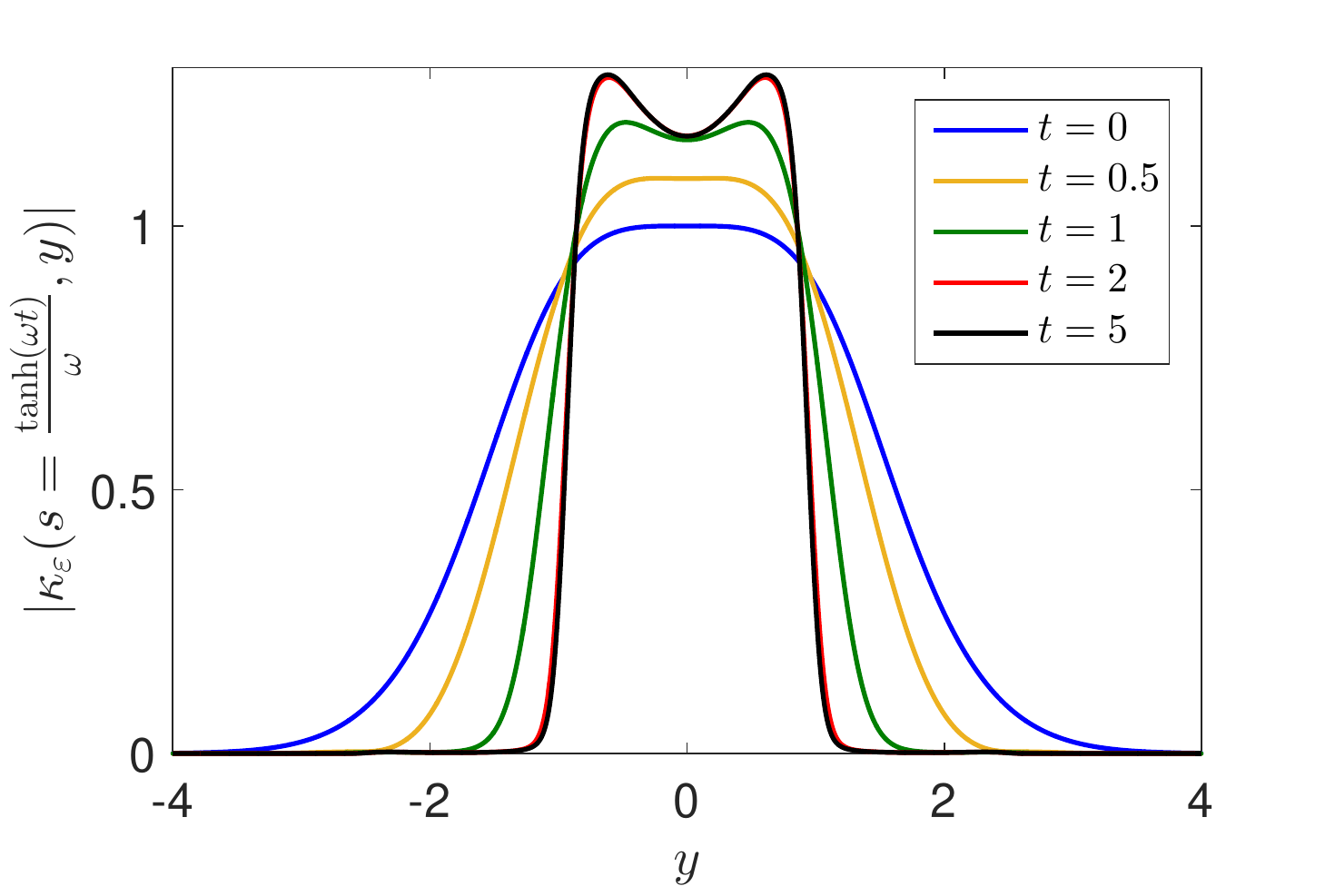}
\includegraphics[width=2.4in,height=1.5in]{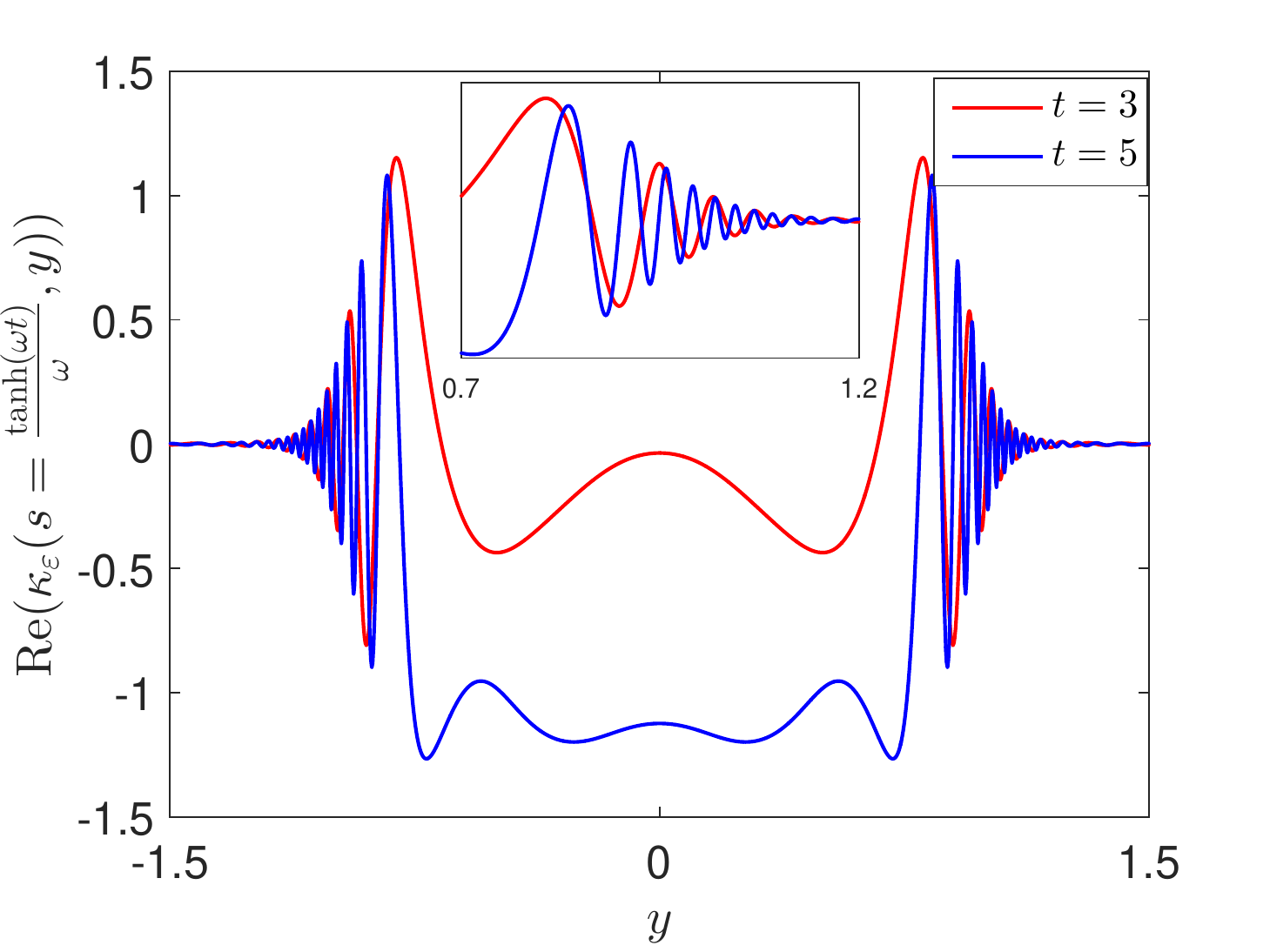}
\end{center}
\caption{Dynamics of $\kappa$ for Cases (i)-(iii) (up to down) in Example 3.}
\label{Ex3}
\end{figure}

Fig. \ref{Ex3} shows the dynamics of $\kappa$ until
$s=\frac{\tanh(\omega T)}{\omega}$ for $T=5$ by using the method Strang I
with $\tau=0.0002$ and $h=1/2^{12}$ on a computational
domain $\Omega=(-10, 10)$ with periodic boundary condition. We observe
that for all three cases, the modulus $|\kappa|$ remain invariant
after some time while the variation of the argument introduces the
oscillation of $\kappa$ in space. This suggests that $|u|$ decreases
exponentially in time, with rate $e^{-\omega t/2}$ in view of \eqref{scheme2}. Moreover, the argument of $\kappa$ introduces more and more oscillation in space as time involves (cf. the right column in Fig. \ref{Ex3}).

\bigskip

\noindent{\bf Example 4}.
We set $\lambda=1$ and $\omega=2$ in the equation \eqref{logNLSrep},
in which case there exists no stationary solution for equation
\eqref{eq:tau-general}. We consider the following two cases for the
initial data $u_0$:\\
(i). $u_0(x)=\sech(x^2/2)\sin(x)$;\\
(ii). $u_0(x)=2e^{-(x-3)^2/4}+2e^{-(x+3)^2}$.

\begin{figure}[h!]
\begin{center}
\includegraphics[width=2.4in,height=1.5in]{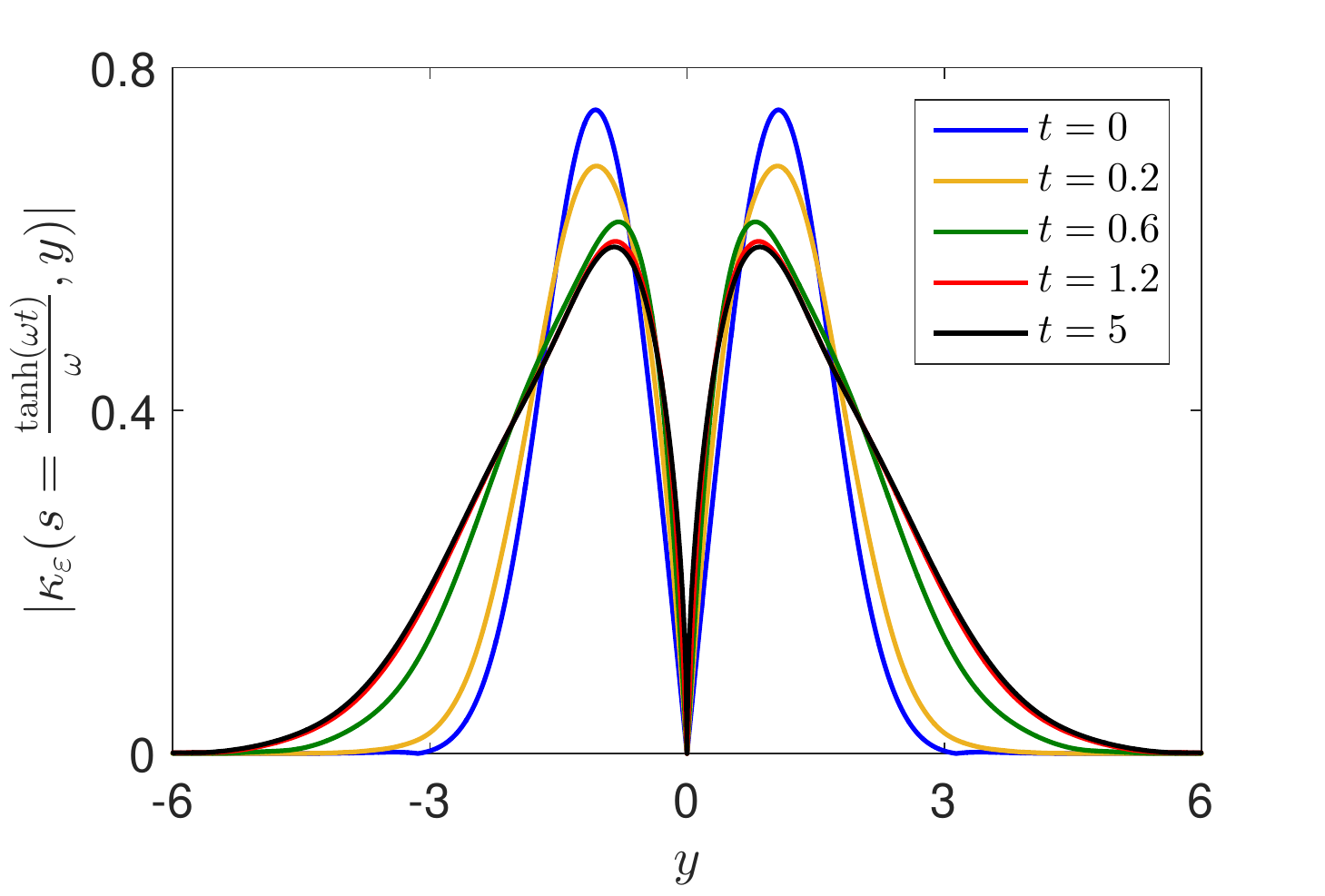}
\includegraphics[width=2.4in,height=1.5in]{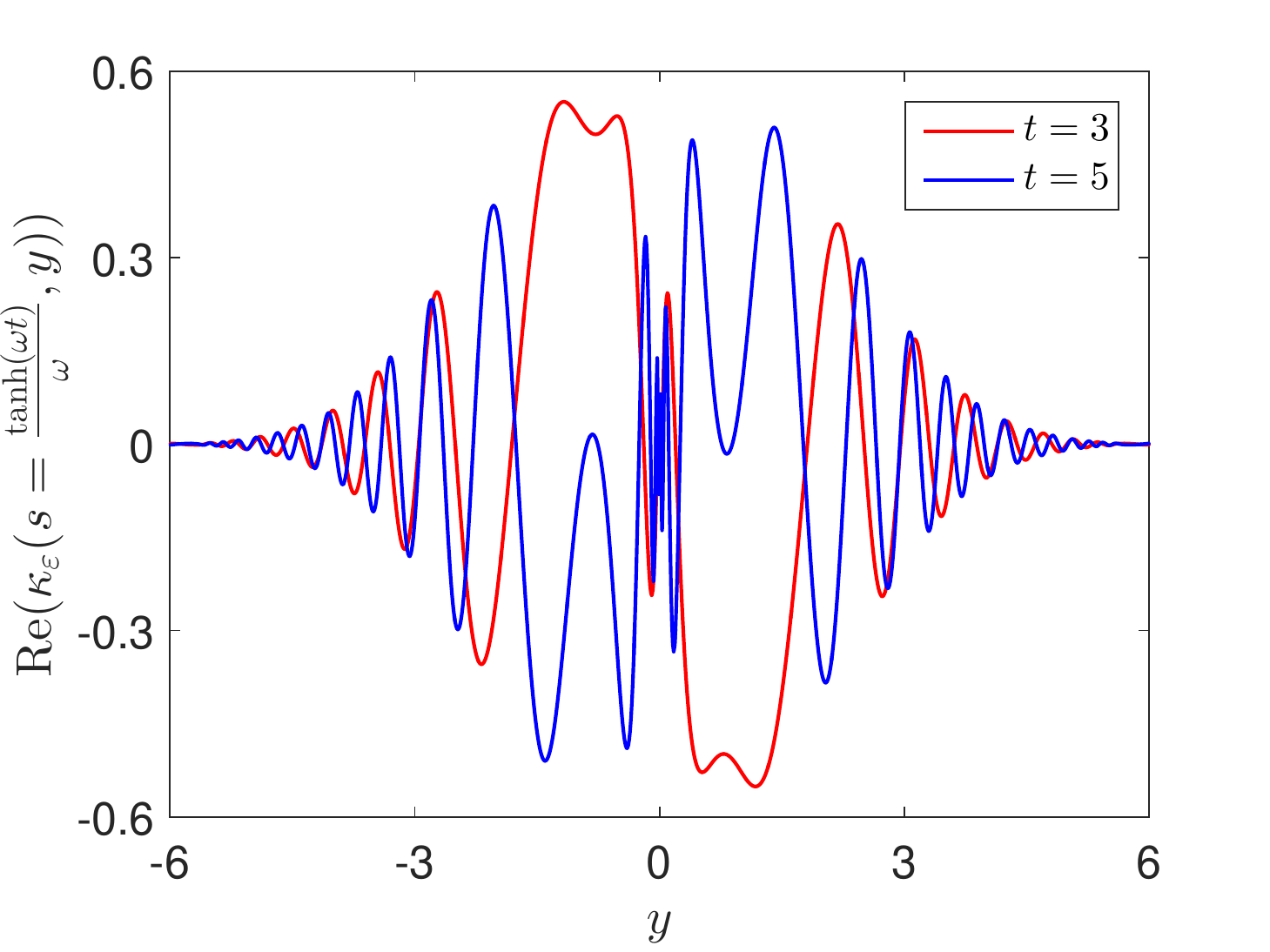}\\
\includegraphics[width=2.4in,height=1.5in]{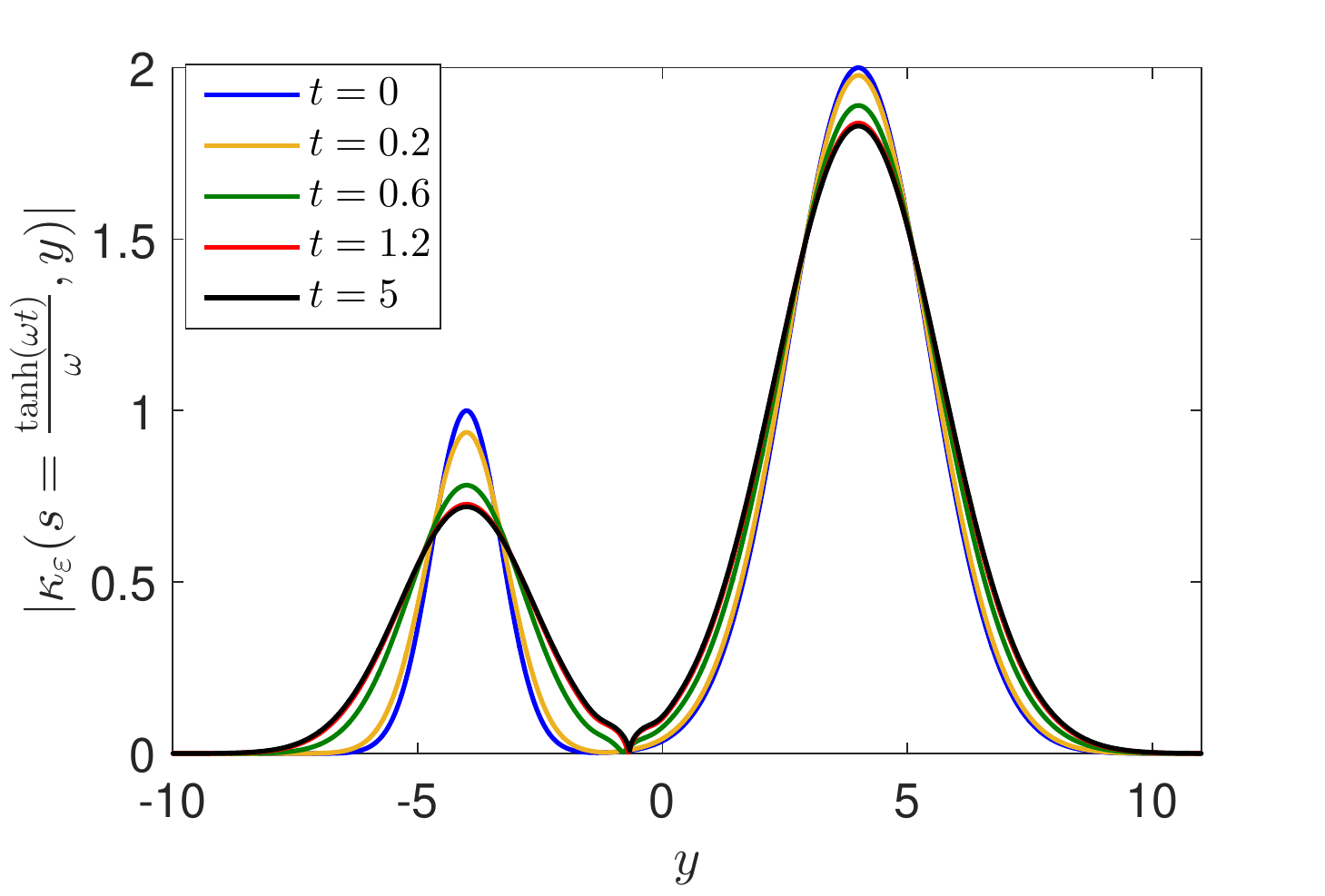}
\includegraphics[width=2.4in,height=1.5in]{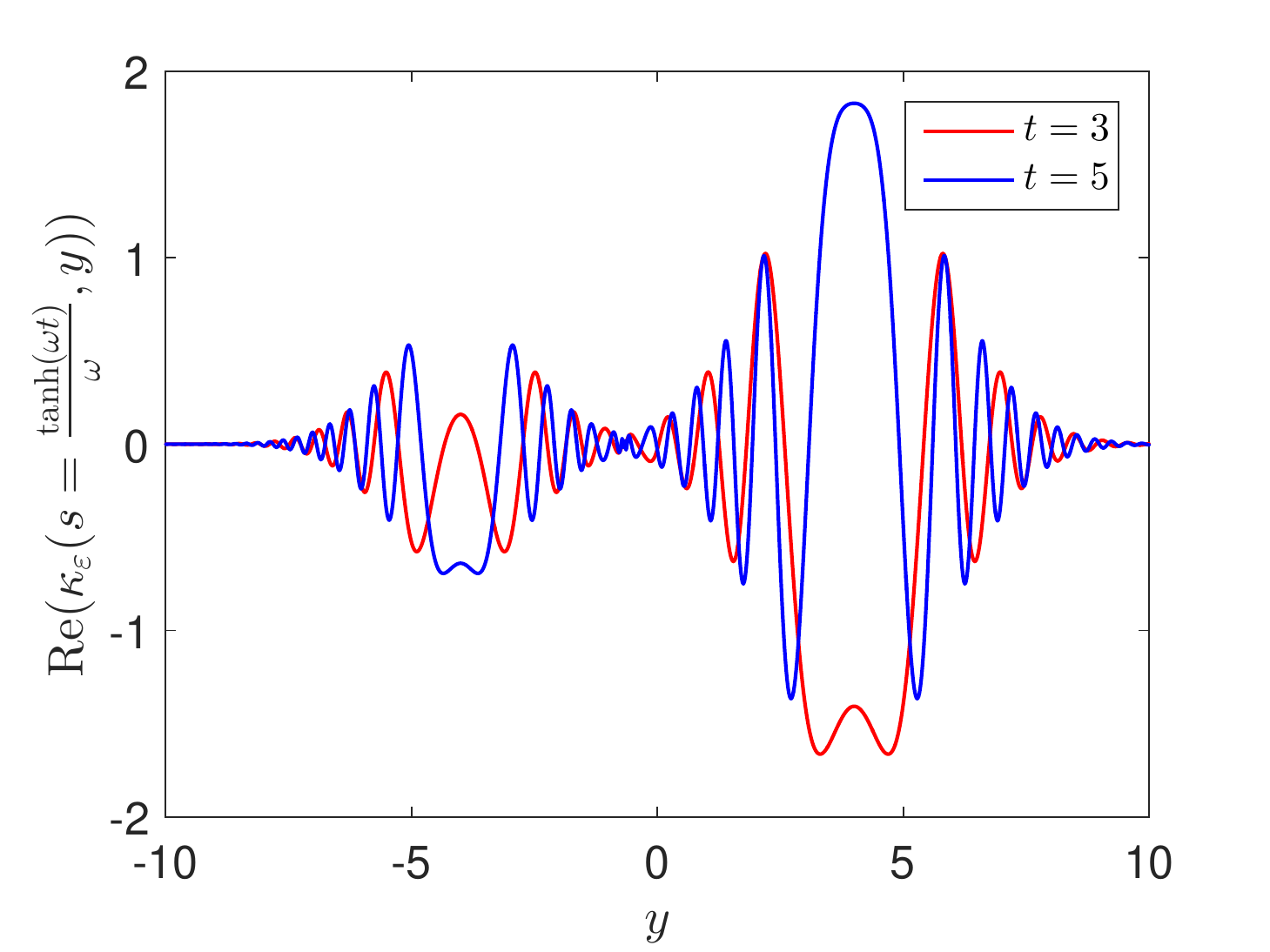}
\end{center}
\caption{Dynamics of $\kappa$ for Cases (i)-(ii) (up to down) in Example 4.}
\label{Ex4}
\end{figure}
Fig. \ref{Ex4} shows the dynamics of $\kappa$ until $s=\frac{\tanh(\omega T)}{\omega}$ for $T=5$ with $\tau=0.0002$ and $h=1/2^{12}$ on a computational domain $\Omega=(-20, 20)$ with periodic boundary condition. We observe that the behavior is similar to that in Example 3, the modulus can attain the invariant state soon and more and more mild oscillations are created as time evolves. With the formulation \eqref{uv}, the dispersion and oscillation properties of $u$ can be clearly revealed.

\medskip

Next we display some numerical experiments for dynamics of the solution $v$ for the equation \eqref{eq:v}, or accordingly the dynamics of $u$ (through \eqref{recu}) for the power nonlinearity case \eqref{eq:NLSrep}. We still use Strang I method for temporal discretization
with $\tau=0.0002$ and $h=1/2^{10}$ as the spatial mesh size on a computational domain $\Omega=(-20, 20)$ with periodic boundary condition.
\medskip

\noindent{\bf Example 5}.
We set $\lambda=1$, $\omega=2$ in the equation \eqref{eq:NLSrep} with different $\sigma$ and initial data $u_0$:\\
(i). $u_0(x)=2e^{-x^2}$;\\
(ii). $u_0(x)=\sech(x^2/2)\sin(x)$.

\begin{figure}[htbp]
\begin{center}
\includegraphics[width=2.4in,height=1.5in]{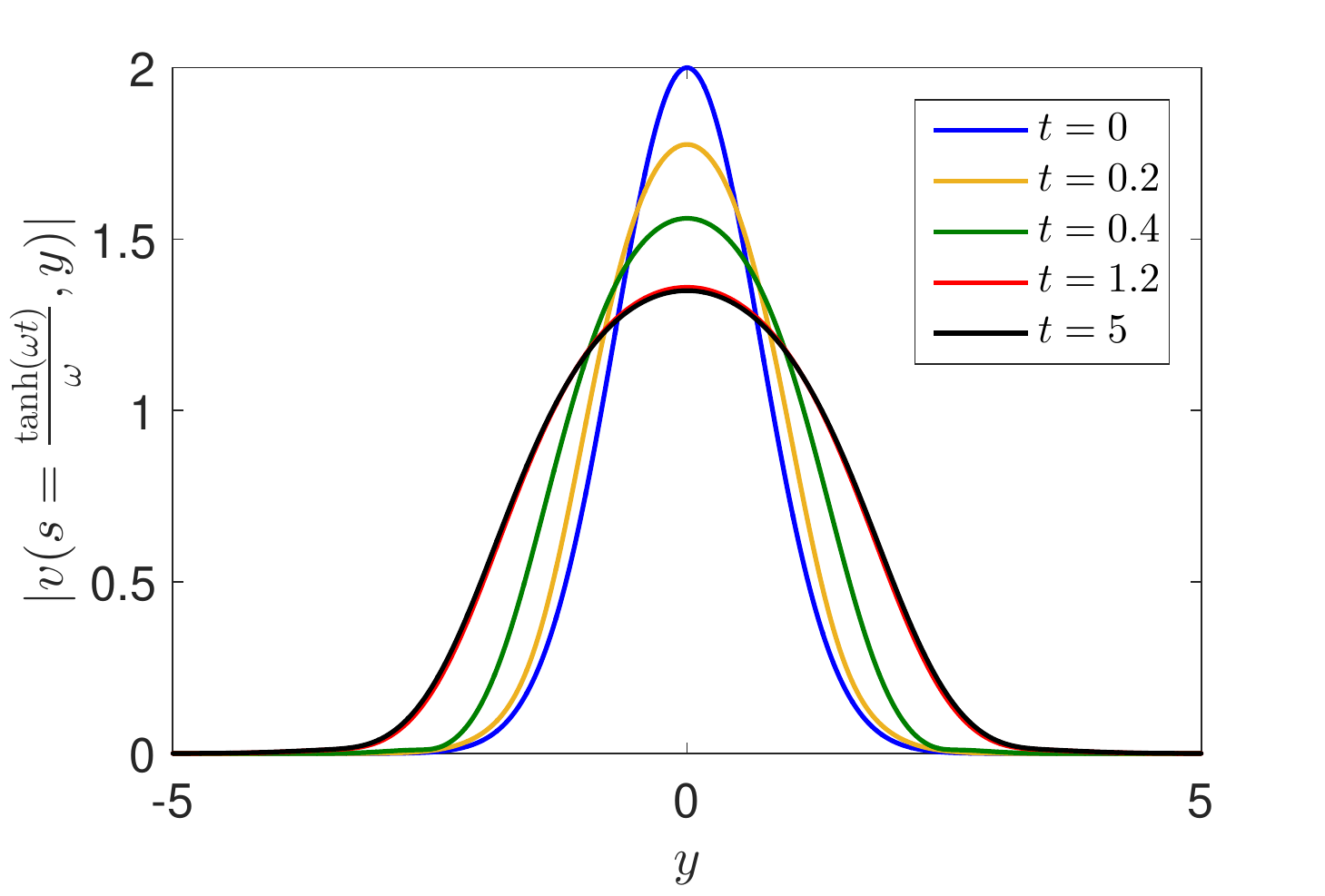}
\includegraphics[width=2.4in,height=1.5in]{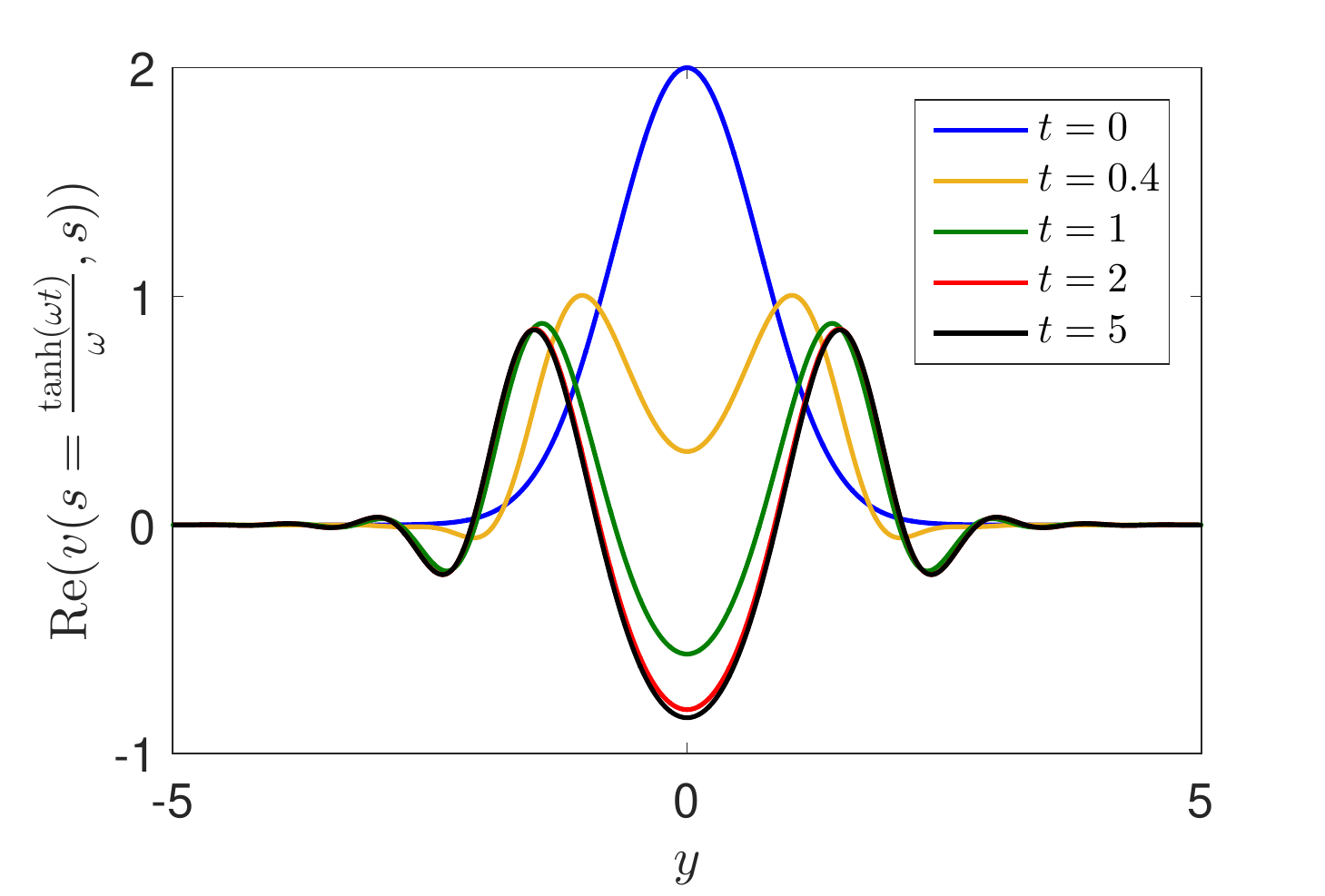}\\
\includegraphics[width=2.4in,height=1.5in]{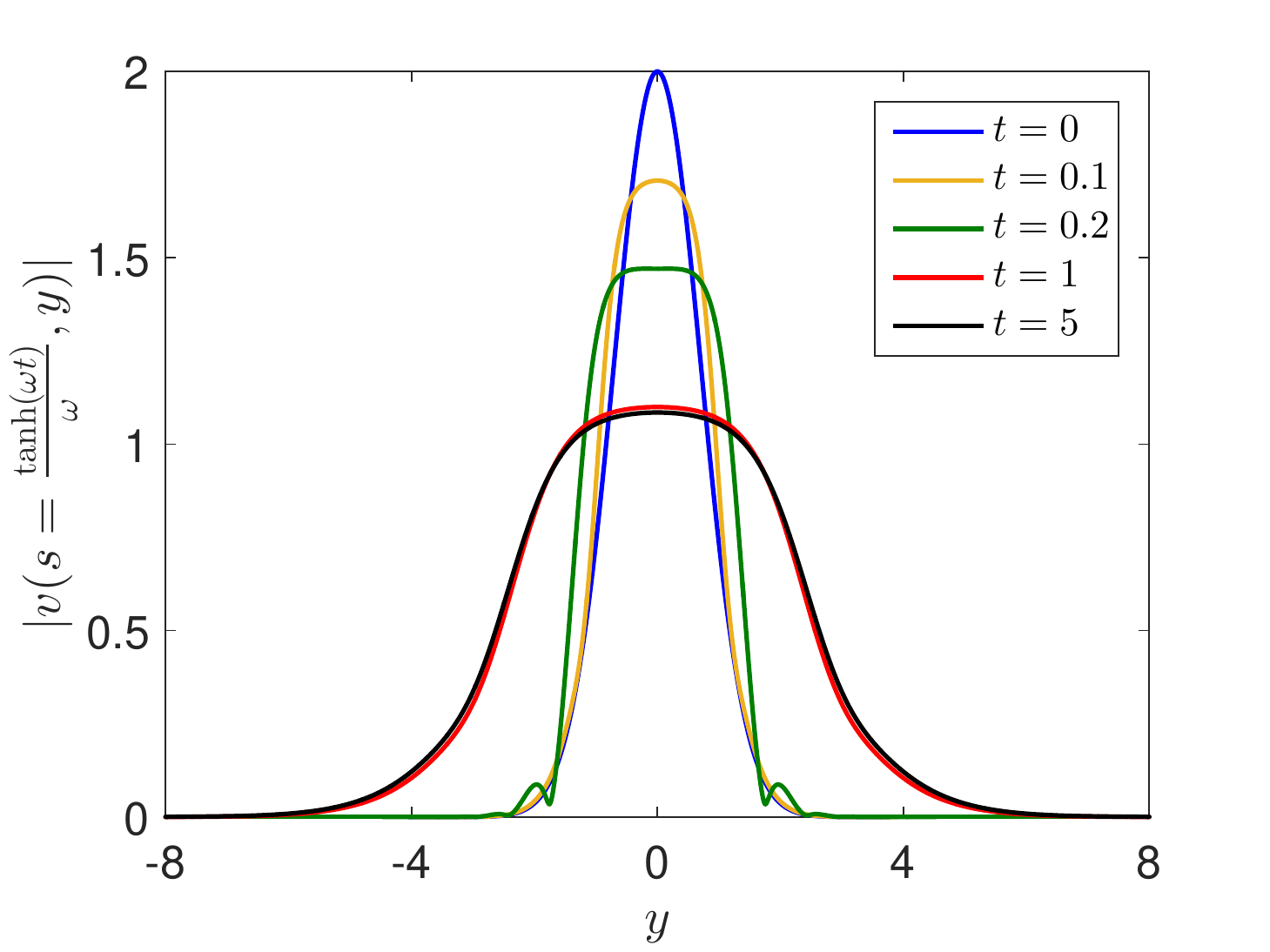}
\includegraphics[width=2.4in,height=1.5in]{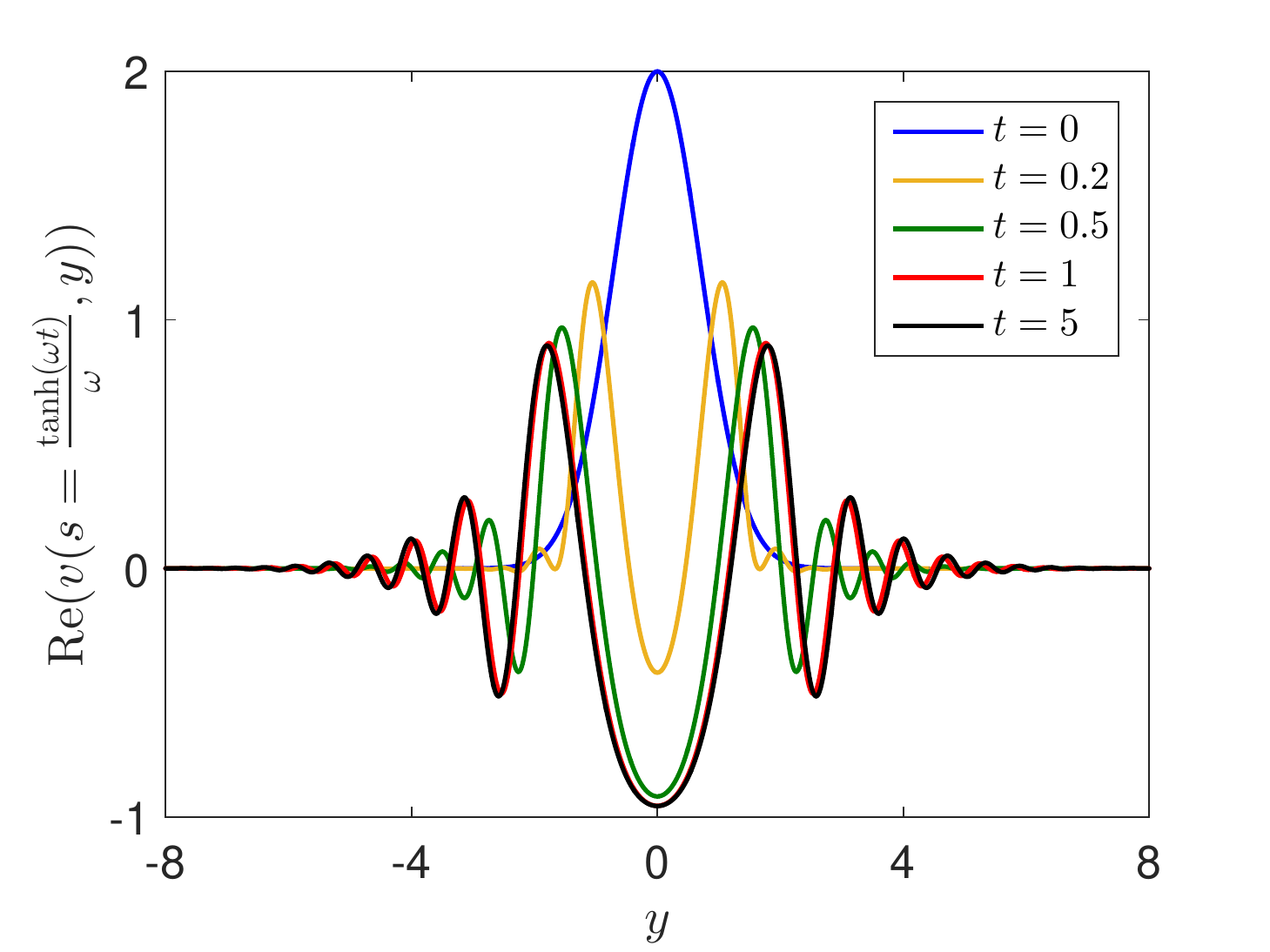}\\
\includegraphics[width=2.4in,height=1.5in]{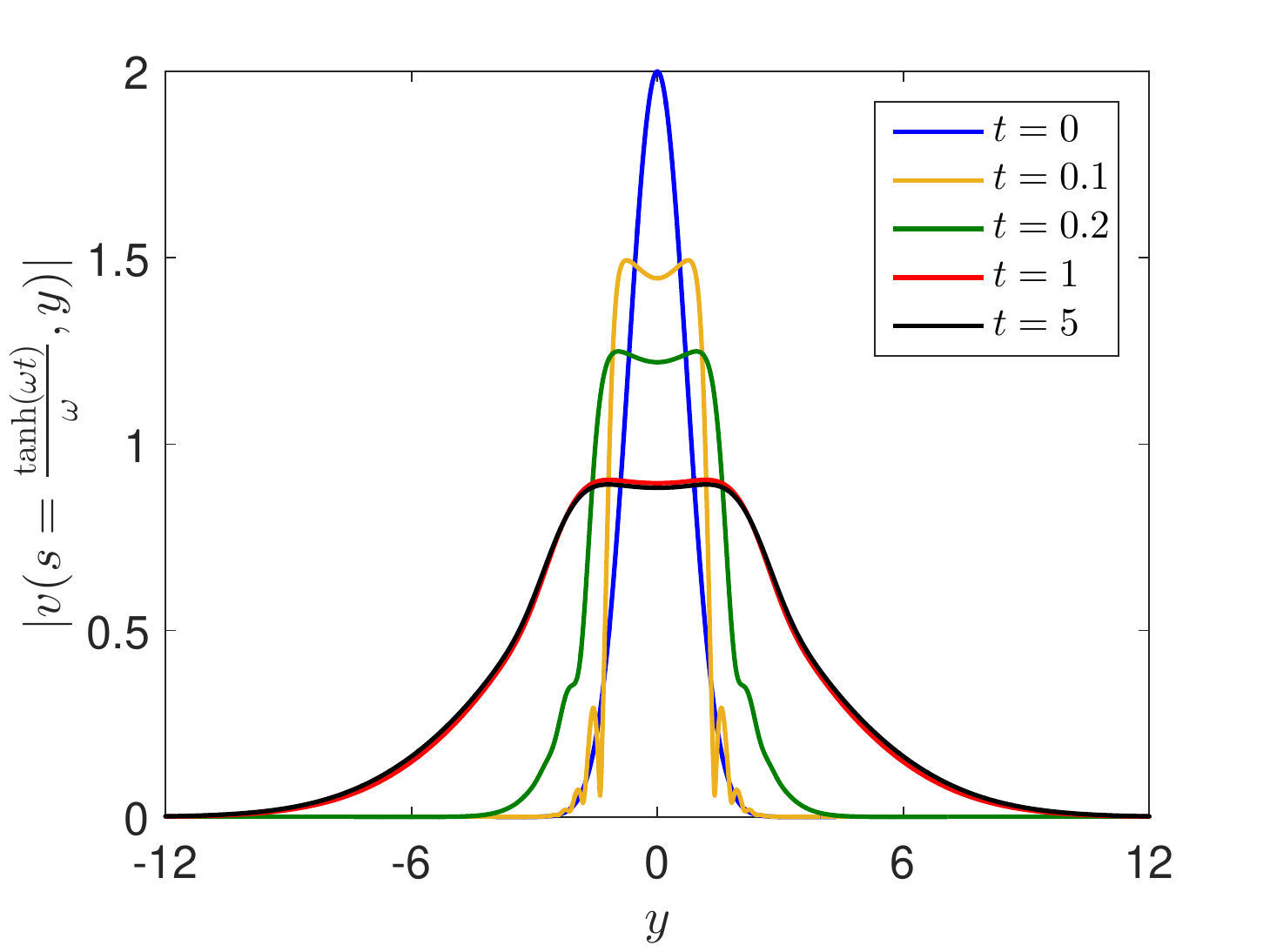}
\includegraphics[width=2.4in,height=1.5in]{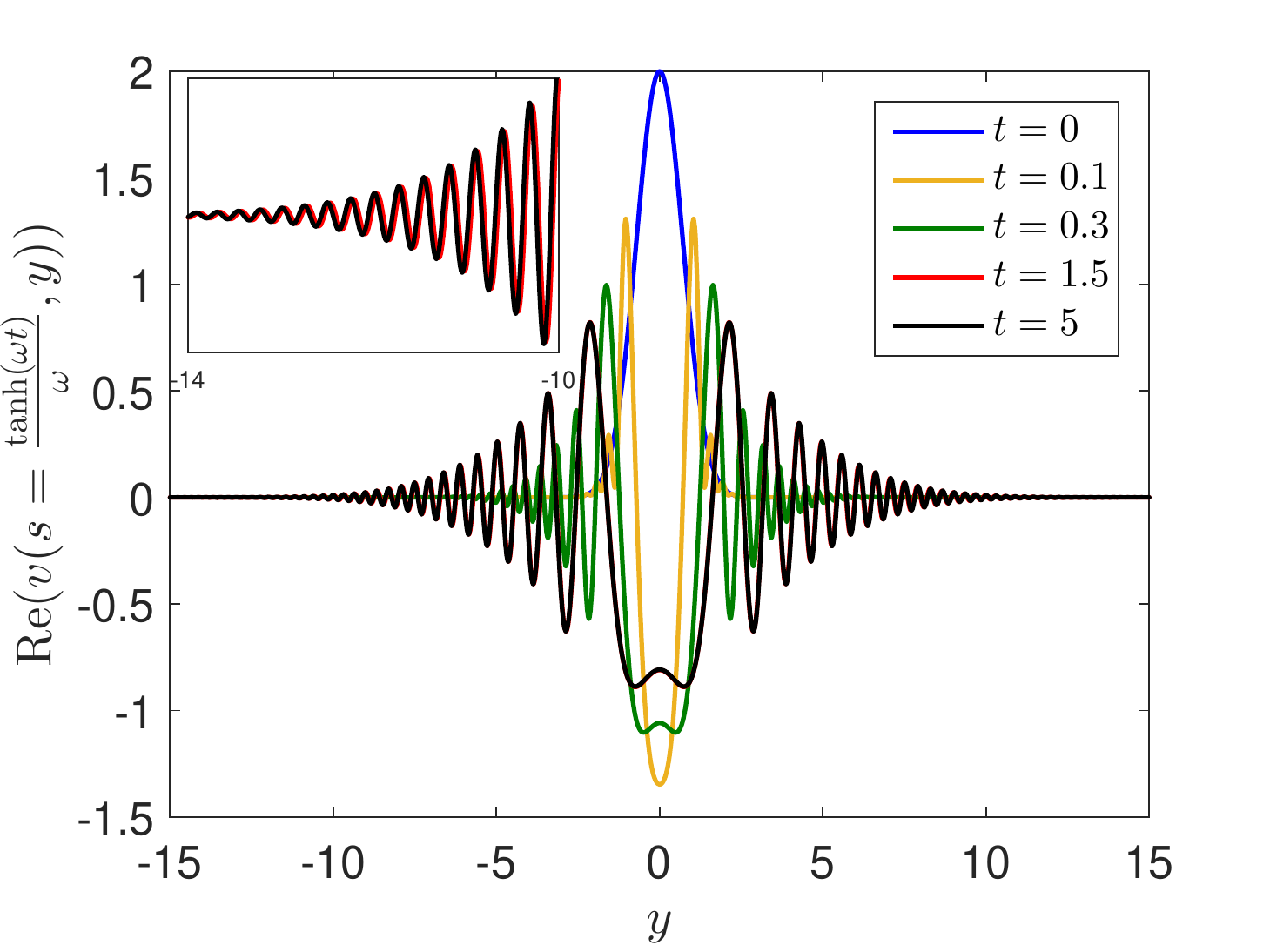}
\end{center}
\caption{Dynamics of $v$ with $\sigma=1, 2, 3$ (up to down) for Case (i) initial data in Example 5.}
\label{Ex5}
\end{figure}

\begin{figure}[htbp]
\begin{center}
\includegraphics[width=2.4in,height=1.5in]{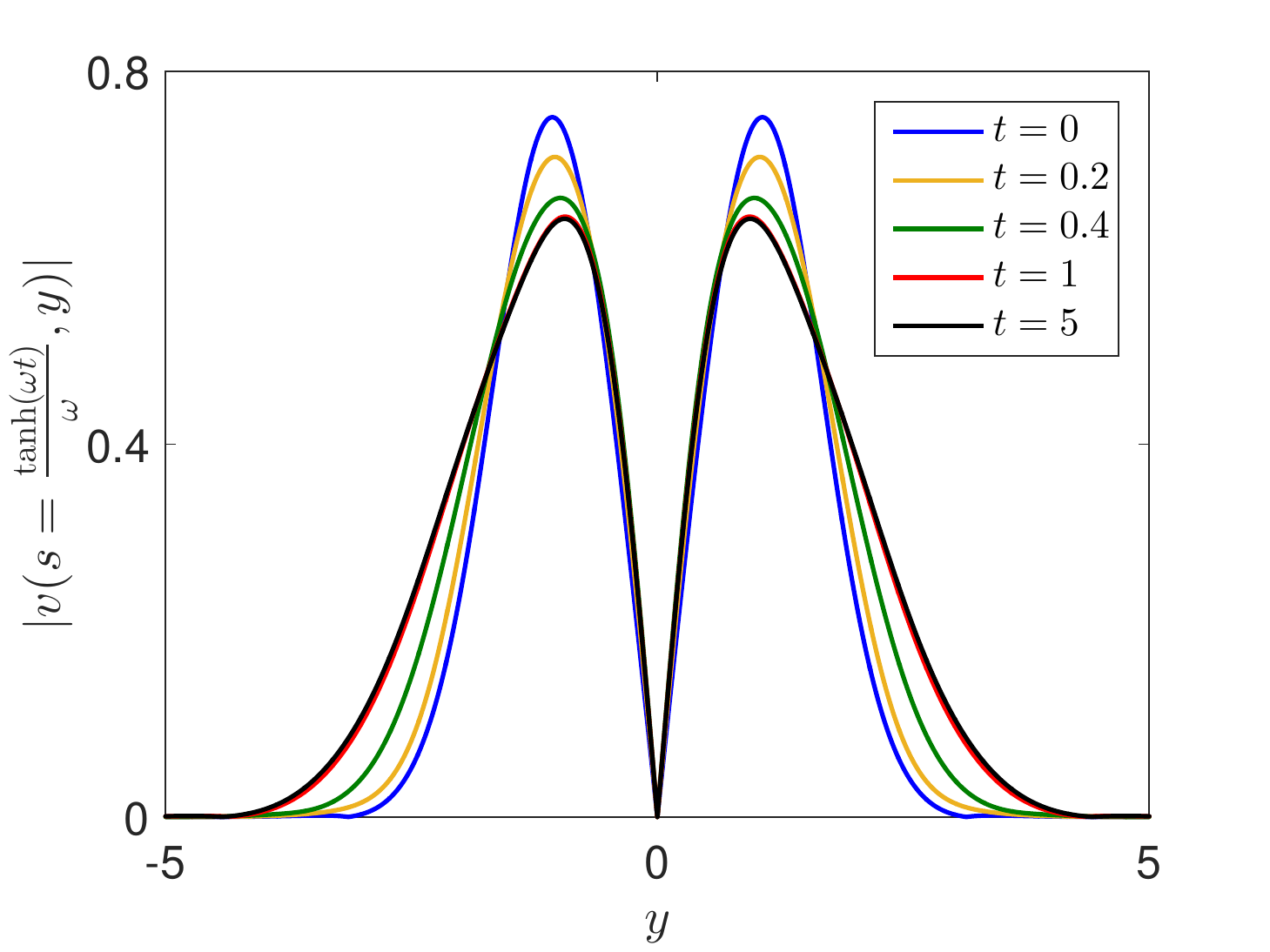}
\includegraphics[width=2.4in,height=1.5in]{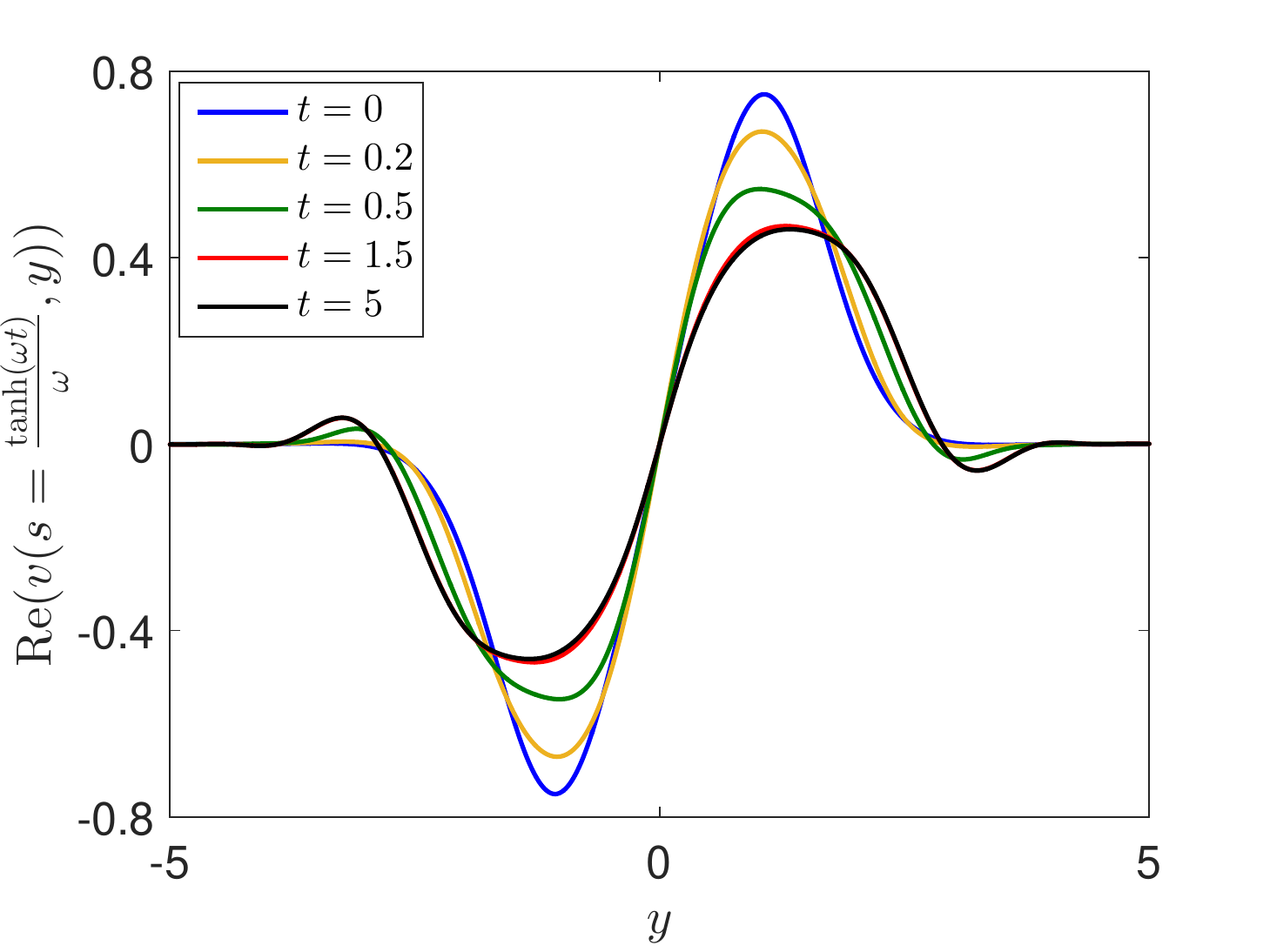}\\
\includegraphics[width=2.4in,height=1.5in]{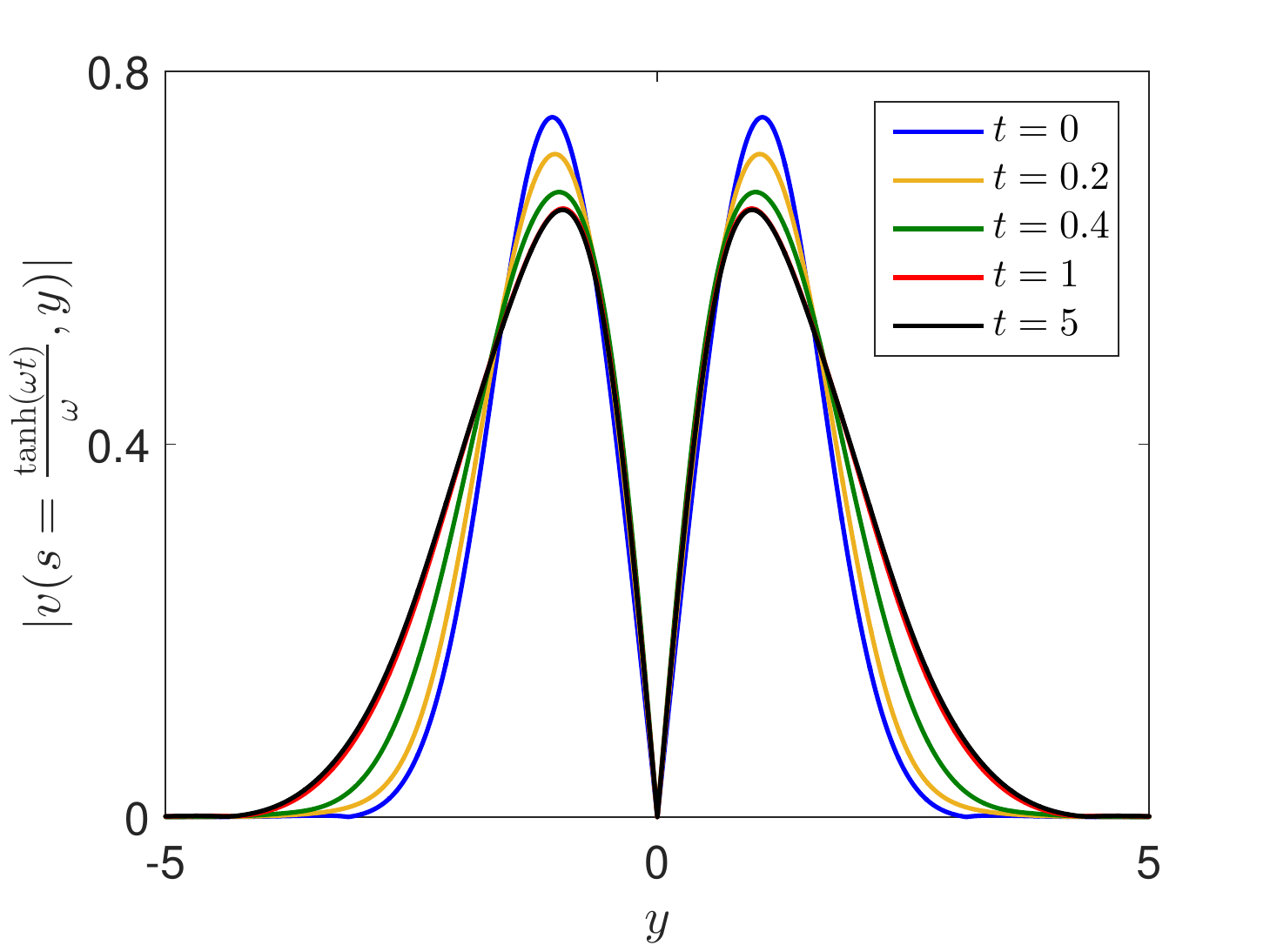}
\includegraphics[width=2.4in,height=1.5in]{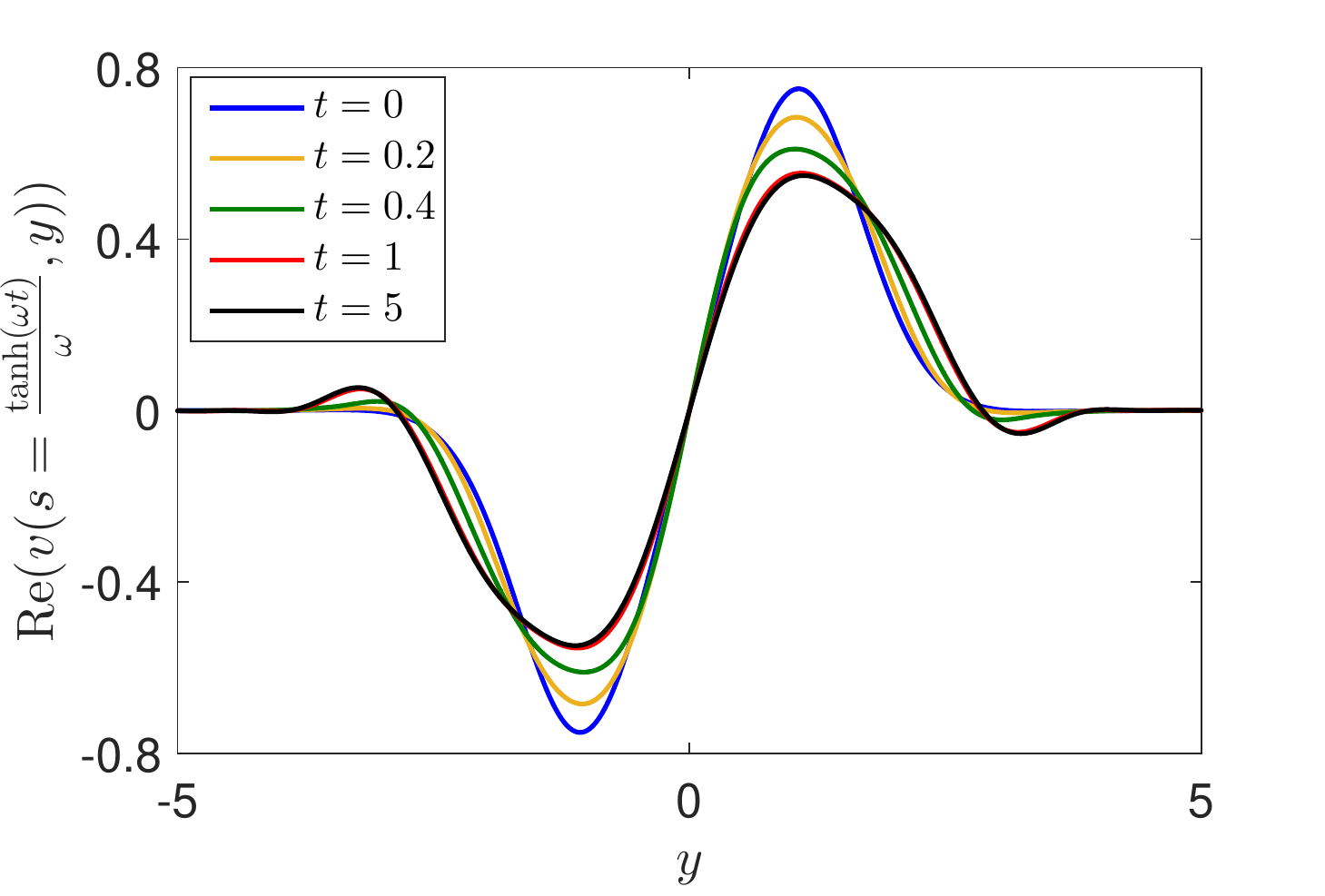}\\
\includegraphics[width=2.4in,height=1.5in]{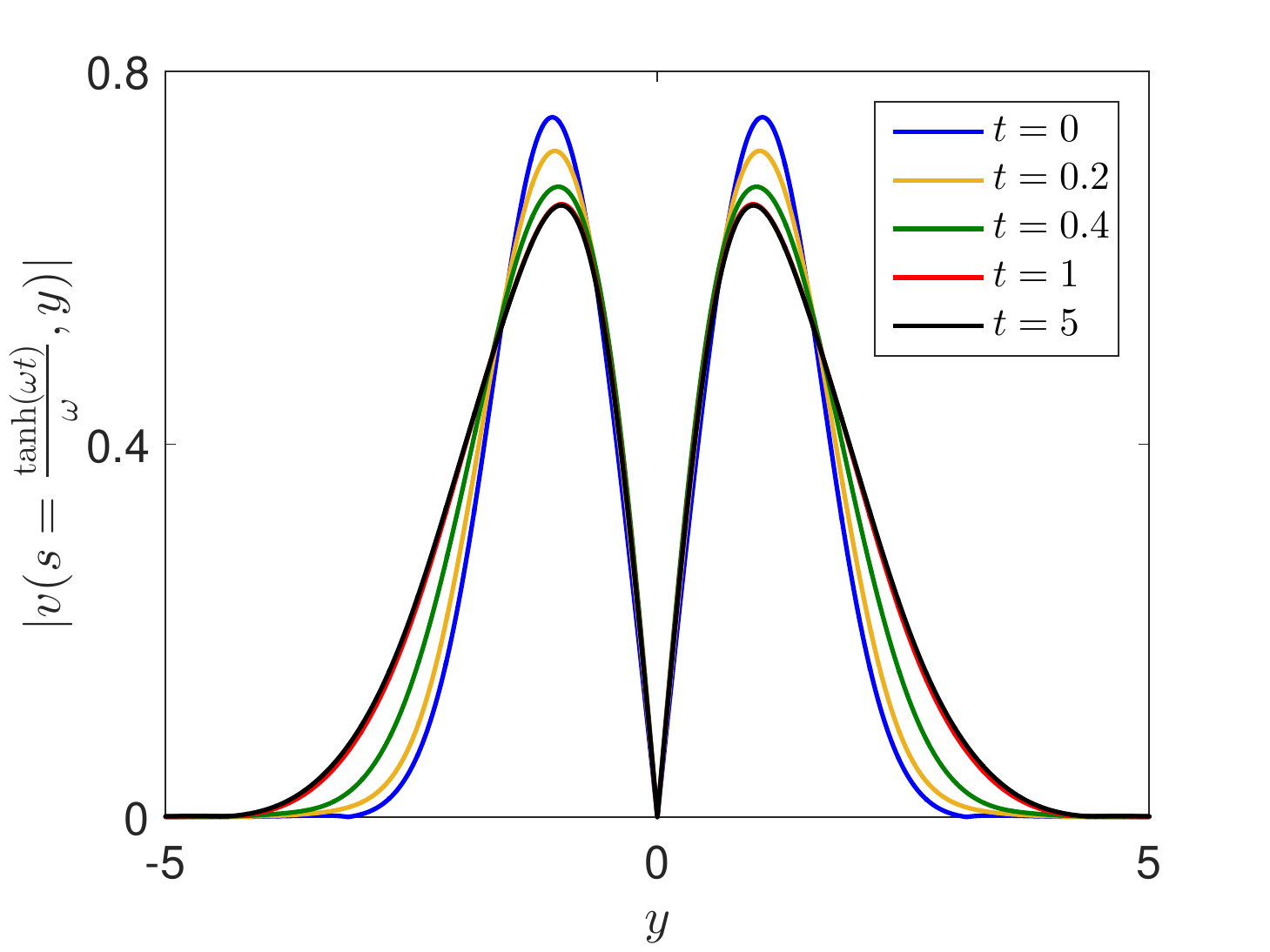}
\includegraphics[width=2.4in,height=1.5in]{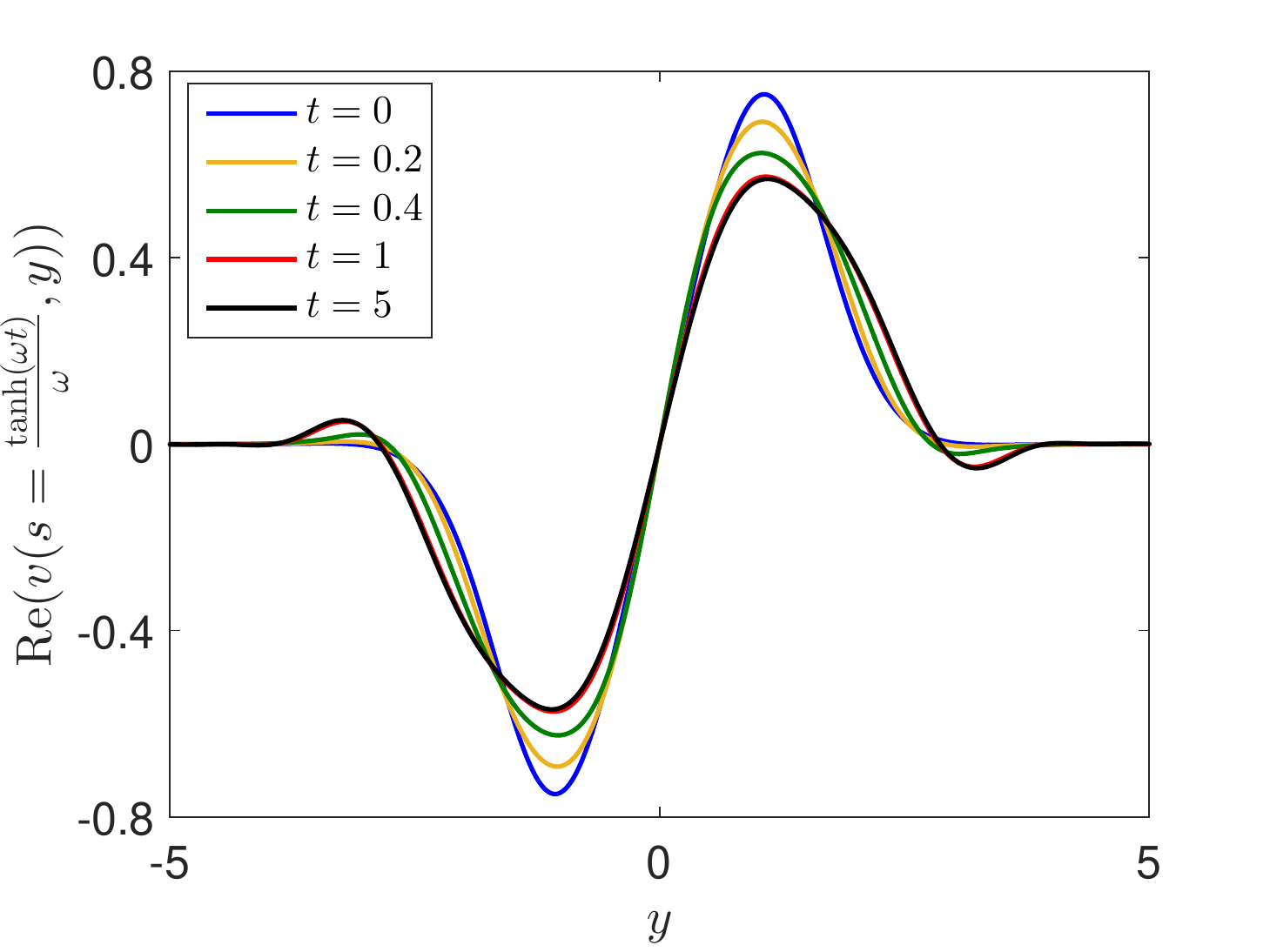}
\end{center}
\caption{Dynamics of $v$ with $\sigma=1, 2, 3$ (up to down) for Case (ii) initial data in Example 5.}
\label{Ex52}
\end{figure}

Figs. \ref{Ex5} and \ref{Ex52} display the dynamics of $v$ for Cases (i) and (ii) initial data, respectively. We observe that for both cases and all chosen $\sigma=1, 2, 3$, both the modulus $|v|$ and the argument varies very slowly (almost invariant) after some time, which is different from the case with logarithmic nonlinearity where the argument varies quickly with respect to time. For Case (i) initial data, the ``essential" support of the (almost) invariant state gets larger (cf. left column in Fig. \ref{Ex5}) and more and more oscillation is created in space (cf. right column in Fig. \ref{Ex5}) when $\sigma$ increases. While for Case (ii) initial data, one can hardly tell any remarkable difference for various choices of $\sigma$, as shown in Fig. \ref{Ex52}.

In the above examples, the numerical solution of \eqref{eq:v}
stays almost invariant after some time. This is understandable by
noticing that the exponent in \eqref{zsol} gets very small as $t_n$
increases, and due to the rapid convergence toward the linear
dynamics, recalled in \eqref{eq:scattering-power}, as well as the
linear behavior, which involves, asymptotically as $t\to +\infty$, the
same rescaling and phase shift as in \eqref{uv}.

\medskip
\noindent{\bf Example 6}.
We set $\lambda=-1$, $\omega=2$ in the equation \eqref{eq:NLSrep} with initial data $u_0(x)=2e^{-x^2}$ and different $\sigma$.

\begin{figure}[htbp]
\begin{center}
\includegraphics[width=2.4in,height=1.5in]{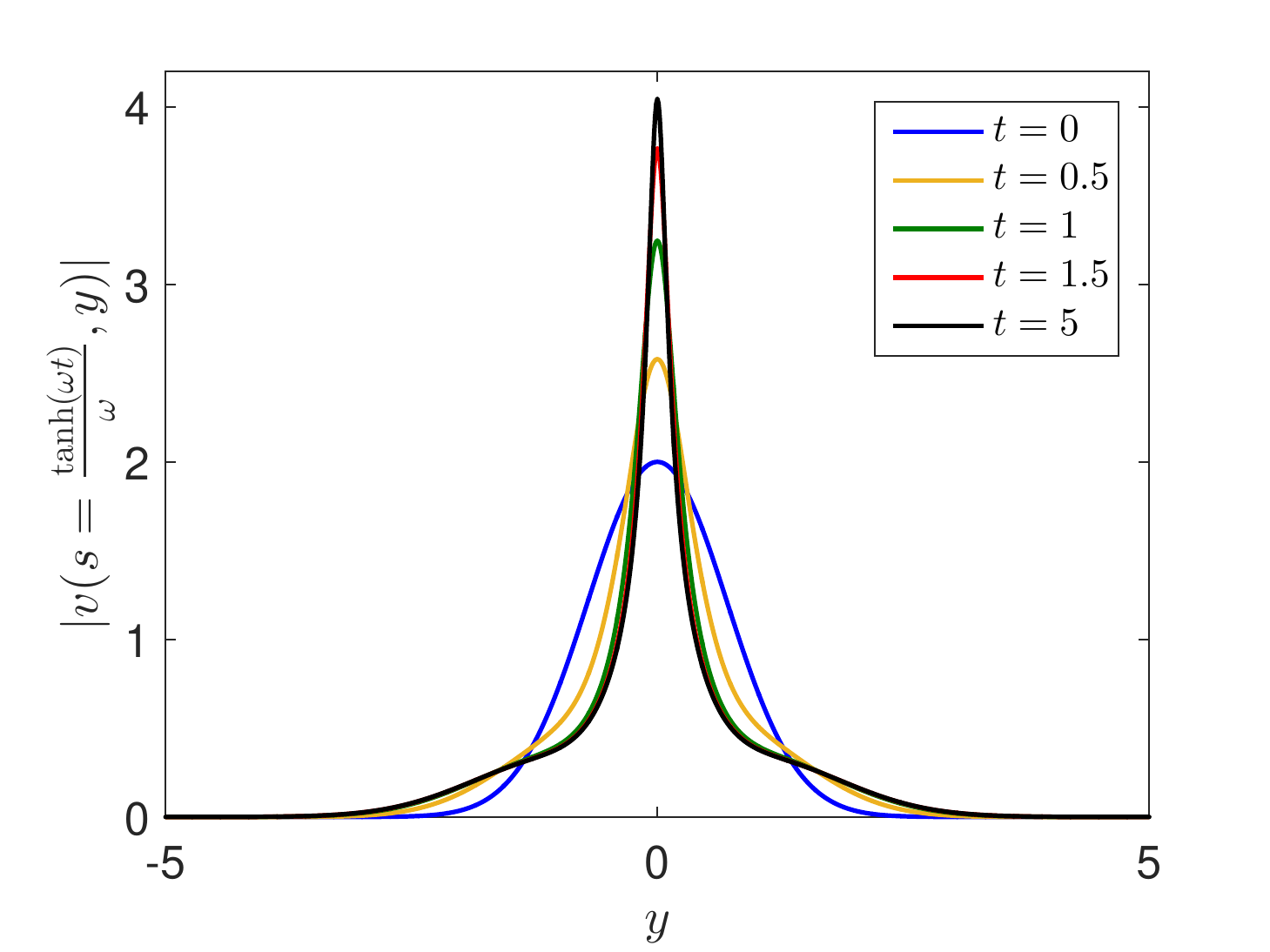}
\includegraphics[width=2.4in,height=1.5in]{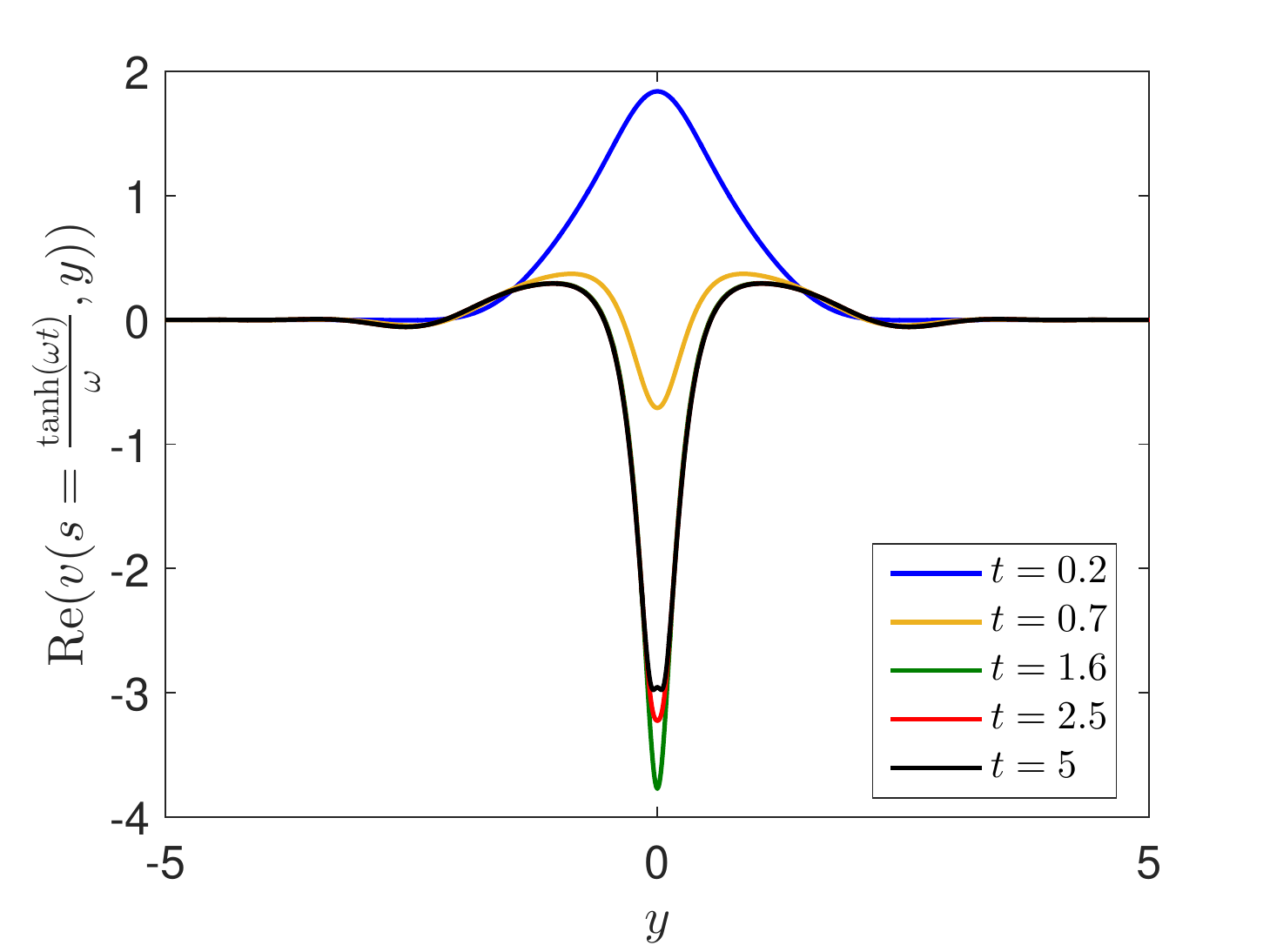}
\end{center}
\caption{Dynamics of $v$ in the equation \eqref{nonl1} with $\sigma=1$ in Example 6.}
\label{Ex6}
\end{figure}
\begin{figure}[htbp]
\begin{center}
\includegraphics[width=2.4in,height=1.5in]{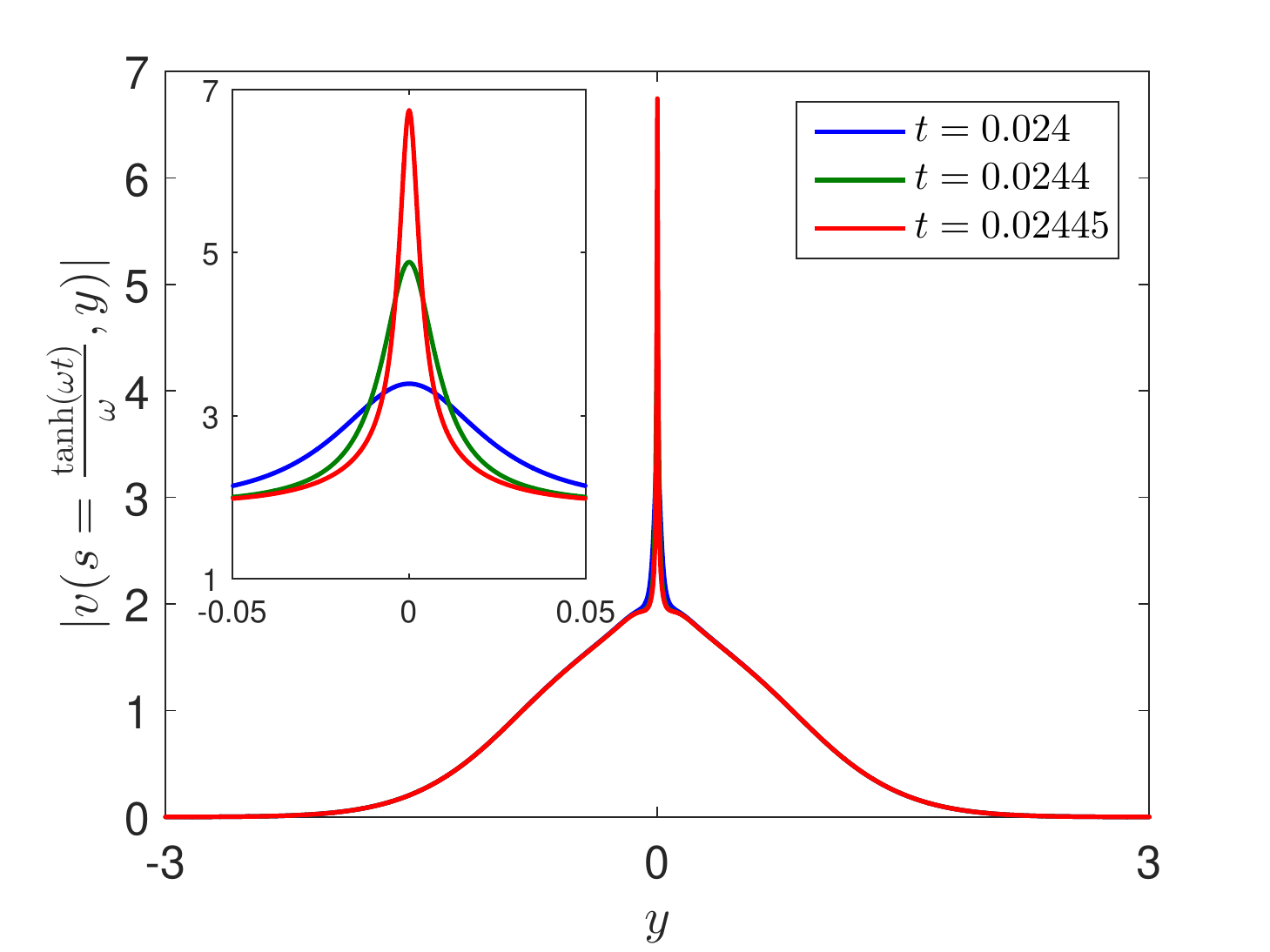}
\includegraphics[width=2.4in,height=1.5in]{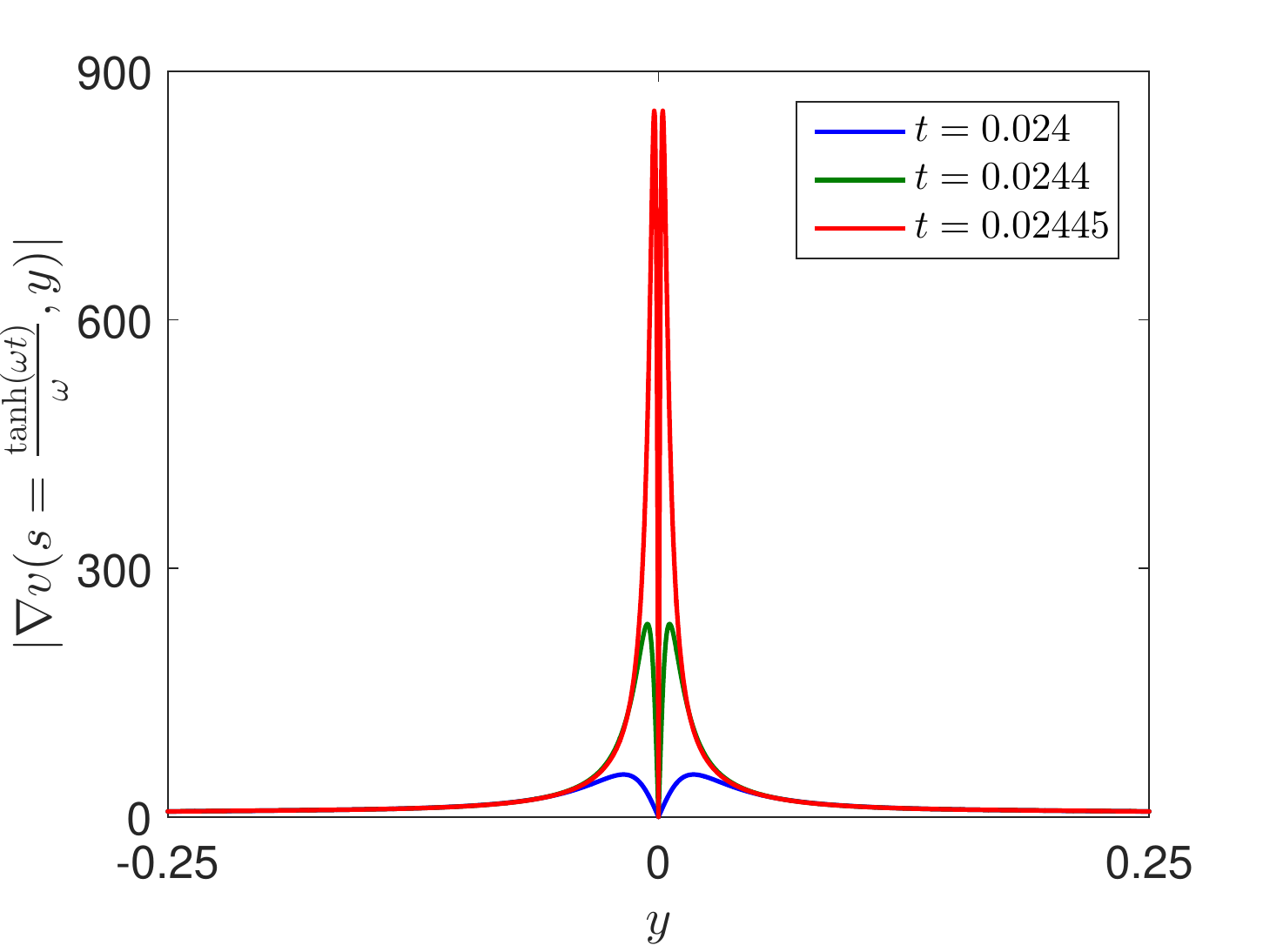}
\end{center}
\caption{Dynamics of $v$ in the equation \eqref{nonl1} with $\sigma=3$ in Example 6.}
\label{Ex62}
\end{figure}

As was shown in \cite{carles2003}, finite-time blow-up might occur
when the parameters $\lambda$, $\omega$, $\sigma$ and $u_0$ satisfy
some condition. Here we choose $\sigma=1, 2, 3$, respectively, and
solve the equation \eqref{nonl1} by Strang I method. For $\sigma=1$,
the dynamics of $v$ is shown in Fig. \ref{Ex6}. Similar results are
obtained as in the case $\lambda>0$, which suggests that in this case,
the solution in $u$ has the same dispersion as before, i.e., it
disperses exponentially in time at the rate $e^{\omega t/2}$. While
for $\sigma=3$, the dynamics of $v$ is displayed in Fig. \ref{Ex62},
from which we observe that the solution concentrates around the origin
quickly. This is in agreement with the virial computation
  leading to \cite[Theorem~1.1, case 3]{carles2003}: in the present case
  we check that
  \begin{equation*}
    \frac{1}{2}\|\nabla u_0\|_{L^2}^2
+\frac{\lambda}{\sigma +1}\|
u_0\|_{L^{2\sigma +2}}^{2\sigma +2}<-\frac{\om^2}{2}\|x u_0\|_{L^2}^2-\om
\left|\operatorname{Im}\int \overline{u_0}x\cdot \nabla_x
u_0 \right|,
  \end{equation*}
as every term can be computed in the Gaussian case,
and so the solution blows up in the future and in the past.
The right plot shows the dynamics of $\nabla v$, which suggests that blow-up \eqref{blowup} occurs at around $T=0.02445$. The situation is similar for $\sigma=2$, which is omitted here for brevity. In one word, the proposed method based on the formulation \eqref{uv} is able to capture the correct dynamics for both dispersive solutions and finite-time blow-up solutions. It is superior than that by solving the original equation
\eqref{eq:NLSrep} directly when the finite-time blow-up does not occur so soon.

\section{Conclusion}
We proposed some time splitting methods for the Schr\"odinger equation with a logarithmic nonlinearity and a repulsive potential based on the generalized lens transform. This transformation can capture the main dispersion and oscillation in the solution. It neutralizes the possible boundary effects and enables the classical splitting methods to work smoothly for the equivalent formulation. This approach was extended to the case with a power nonlinearity. Error estimates for the semi-discrete Lie-Trotter splitting method were established for the equation with both types of nonlinearities. Finally we investigated the dynamics and revealed the dispersion property of the mentioned two types of equations with different parameters by employing the proposed methods.

\bibliographystyle{siam}
\bibliography{biblio}
\end{document}